\documentclass[12pt]{article}
\usepackage{amsmath,amsthm}
\usepackage{amsfonts,amssymb,color}
\usepackage{pstricks}
\usepackage{epsfig}
\usepackage{hyperref}

\newcommand{\dih}[2]{DIH_{#1}(#2)}
\newcommand{\drt}[1]{\frac{1}{\sqrt{2}} #1}
\newcommand{\drtp}[1]{  \frac{1}{\sqrt{2}}\left( #1\right)}
\newcommand{\drtpp}[1]{ \left( \frac{1}{\sqrt{2}}\left( #1\right)\right)}
\newcommand{\mspan}{\mathrm{span}}
\newcommand{\Rspan}[1]{\mathrm{span}_{\mathbb{R}}(#1)}
\newcommand{\Dih}[1]{DIH_{#1}}

\def\refpp#1{(\ref {#1})}

\def\noi{\noindent}

\def\pf{\noi{\bf Proof.\ \,}}
\def\eop{{$\square$}}
\def\qed{{$\square$}}

\def\labtt#1{\label {#1} }
\def\labttr#1{\label {#1} }

\def\a{\alpha}
\def\b{\beta}
\def\g{\gamma}

\def\vep{\varepsilon}

\def\DD{{{\mathbb Z} [D]}}
\def\FF{{\mathbb F}}
\def\QQ{{\mathbb Q}}
\def\RR{{\mathbb R}}

\def\ZZ{{\mathbb Z}}

\def\la{\langle}
\def\ra{\rangle}
\def\<{\langle}
\def\>{\rangle}

\def\bs{\it}            

\def\dim{{\bs dim}}

\def\det{{\bs det}}

\def\half{\frac{1} 2}
\def\fourth{\frac{1} 4}

\def\third{{\frac{1} 3}}
\def\tthird{\frac{2}{3}}

\def\dual#1{#1^*}        

\def\kron#1#2{\delta_{#1#2}}  

\def\mv#1{MinVec(#1)}

\def\ratholoex#1{2^{1+2#1}_+\Omega^+(2#1,2)}

\def\dg#1{{\cal D}({#1})}

\def\bw#1{BW_{{2^{#1}}}}

\def\brw#1{BRW^+(2^{#1})}

\def\ses#1#2#3{0\rightarrow #1 \rightarrow #2 \rightarrow #3 \rightarrow 0} 

\begin{document}

\newtheorem{mainthm}{Main Theorem}[section]
\newtheorem{thm}{Theorem}[section]
\newtheorem{prop}[thm]{Proposition}
\newtheorem{lem}[thm]{Lemma}
\newtheorem{coro}[thm]{Corollary}
\newtheorem{conj}[thm]{Conjecture}
\newtheorem{conv}[thm]{Convention}
\newtheorem{hyp}[thm]{Hypothesis}
\theoremstyle{definition}
\newtheorem{de}[thm]{Definition}
\newtheorem{nota}[thm]{Notation}
\newtheorem{ex}[thm]{Example}
\newtheorem{proc}[thm]{Procedure}
\newtheorem{rem}[thm]{Remark}

\begin{center}
{\Large  $EE_8$-lattices and dihedral groups}

\vspace{10mm}
Robert L.~Griess Jr.
\\[0pt]
Department of Mathematics\\[0pt] University of Michigan\\[0pt]
Ann Arbor, MI 48109  \\[0pt]
USA  \\[0pt]
\vskip 1cm

Ching Hung Lam
\\[0pt]
Department of Mathematics and \\[0pt]
National Center for Theoretical Sciences\\[0pt]
National Cheng Kung University\\[0pt]
Tainan, Taiwan 701\\[0pt]
\vskip 1cm
\end{center}

\abstract{ We classify integral rootless lattices
which are sums of pairs of $EE_8$-lattices (lattices isometric to
$\sqrt 2$ times the $E_8$-lattice) and which define  dihedral groups
of orders less than or equal to $12$.  Most of these may be seen in the Leech
lattice. Our classification may help understand Miyamoto
involutions on lattice type vertex operator algebras and give a
context for the dihedral groups which occur in the Glauberman-Norton
moonshine theory.  }

\tableofcontents

\section{Introduction}

By {\it
lattice}, we mean a finite rank free abelian group with rational
valued, positive definite symmetric bilinear form.  A \textit{root}
in an integral lattice is a norm $2$ vector. An integral lattice is
\textit{rootless} if it has no roots.
The notation $EE_8$ means $\sqrt 2$ times the famous $E_8$ lattice.

In this article, we classify pairs of $EE_8$-lattices
which span an integral and rootless lattice and whose associated involutions
(isometries of order 2) generate a dihedral group of order at most
12.   Examples of such pairs are easy to find within familiar lattices, such as the Barnes-Wall lattices of ranks 16 and 32 and the Leech lattice, which has rank 24.

Our main theorem is as follows.  These results were announced in \cite{gl}.  

\begin{mainthm}\label{mainthm}
Let $M, N\cong EE_8$  be sublattices in a Euclidean space such that
$L=M+N$ is integral and rootless. Suppose that the involutions
associated to $M$ and $N$ \refpp{ssd} generate a dihedral group of
order less than or equal to $12$.   Then the possibilities
for $L$ are
listed in Table \ref{NREE8SUM} and all these possibilities exist.  The lattices in Table \ref{NREE8SUM} are uniquely determined (up to isometry of pairs $M, N$) by the notation in column 1 (see Table \ref{notationandterminology}).  Except for $\dih{4}{15}$, all of them embed as sublattices of the Leech lattice.
\end{mainthm}

\begin{table}[bht]
\caption{ \bf  NREE8SUMs: integral rootless lattices which are sums of
$EE_8$s}
\begin{center}
\begin{tabular}{|c|c|l|c|c|}
\hline
 Name & $\<t_M,t_N\>$ & Isometry type of $L$ (contains)& $\dg L$ & In Leech? \cr
 \hline \hline
 $\dih{4}{12}$ & $Dih_4$ & $\geq  DD_4^{\perp 3}$ & $1^4 2^6 4^2$ & Yes\cr
 \hline
 $\dih{4}{14}$ & $Dih_4$ & $\geq AA_1^{\perp 2} \perp DD_6^{\perp 2}$ & $1^4 2^8 4^2$  &Yes\cr
 \hline
 $\dih{4}{15}$ & $Dih_4$ & $\geq AA_1\perp EE_7^{\perp 2}$ & $1^22^{14}  $  & No\cr
 \hline $\dih{4}{16}$ & $Dih_4$ & $\cong EE_8\perp EE_8$ & $2^{16}$ &   Yes\cr
 \hline $\dih{6}{14}$ & $Dih_6$ & $\geq AA_2\perp A_2\otimes E_6 $ & $1^7 3^3 6^2$ &Yes\cr
 \hline
 $\dih{6}{16}$ & $Dih_6$ & $\cong A_2\otimes E_8 $ & $1^8 3^8 $ &Yes\cr
 \hline
 $\dih{8}{15}$ & $Dih_8$ & $ \geq AA_1^{\perp 7}\perp EE_8$     & $1^{10}4^5 $&   Yes\cr
 \hline
 $\dih{8}{16, DD_4}$ & $Dih_8$ & $\geq DD_4^{\perp 2} \perp  EE_8$  & $1^8 2^4 4^4$  &Yes\cr
 \hline
 $\dih{8}{16, 0}$ & $Dih_8$ & $\cong BW_{16} $ & $1^8 2^8$ &Yes\cr
 \hline $\dih{10}{16}$ & $Dih_{10}$ & $\geq A_4\otimes A_4 $ & $1^{12}5^4$  &Yes\cr
 \hline $\dih{12}{16}$ & $Dih_{12}$ & $\geq AA_2\perp AA_2 $ & &\cr
 & & $\ \perp A_2\otimes E_6 $ & $1^{12}6^4 $  & Yes\cr
 \hline \hline
\end{tabular}
\medskip

$X^{\perp n}$ denotes the orthogonal sum of $n$ copies of the lattice $X$.
\end{center}
\label{NREE8SUM}
\end{table}%

\begin{table}[bht]
\caption{ \bf  Containments of NREE8SUM}
\label{contain}
\begin{center}
\begin{tabular}{|c|cc|}
\hline
 Name &  & Sublattices  \cr \hline
  $\dih{8}{15}$ &\qquad  & $\dih{4}{12}$ \cr
 \hline
 $\dih{8}{16,DD_4}$& & $\dih{4}{12}$ \cr
 \hline
 $\dih{8}{16, 0}$ &  & $\dih{4}{16}$\cr
  \hline
  $\dih{12}{16}$ &  & $\dih{4}{12}$, $\dih{6}{14}$ \cr
 \hline
\end{tabular}
\end{center}
\end{table}

Our methods are probably good enough to determine all the cases where $M+N$ is integral, but such a work would be quite long.

This work may be considered purely  as a study of positive definite integral lattices.
Our real motivation, however, is the evolving theory of vertex operator algebras (VOA) and their automorphism groups, as we shall now explain.

The primary connection between the Monster and vertex operator
algebras was established in \cite{flm}. Miyamoto showed \cite{Mi1}
that there is a bijection between the conjugacy class of $2A$
involutions in the Monster simple group and conformal vectors of
central charge $\half$ in the moonshine vertex operator algebra
$V^\natural$. The bijection between the $2A$-involutions and
conformal vectors offers an opportunity to study, in a VOA context,
the McKay observations linking the extended $E_8$-diagram and pairs
of $2A$-involutions \cite{lyy}.  This McKay theory was originally  described in purely finite group theory terms.

Conformal vectors of central charge $\half$ define automorphisms of
order 1 or 2 on the VOA, called {\it Miyamoto involutions} when they
have order 2.  They were originally defined in \cite{Mi}; see also \cite{Mi1}.  Such conformal vectors are not found in most VOAs but
are common in many VOAs of great interest,   mainly lattice type
VOAs and twisted versions \cite{dmz,dlmn}. Unfortunately, there are
few general, explicit formulas for such conformal vectors in lattice type
VOAs. We know of two.  The first such formula (see \cite{dmz}) is
based on a norm 4 vector in a lattice.  The second such formula (see
\cite{dlmn}, \cite{0+102}) is based on a sublattice which is
isometric to $EE_8$. This latter formula indicates special interest
in $EE_8$ sublattices for the study of VOAs.

We call the dihedral group generated by a pair of Miyamoto
involutions a {\it Miyamoto dihedral group}. Our assumed upper bound
of 12 on the order of a Miyamoto dihedral group is motivated by  the
fact that in the Monster, a pair of $2A$ involutions generates a dihedral
group of order at most 12 \cite{gms}.   Recently, Sakuma \cite{Sa}
announced that 12 is an upper bound for the order of a Miyamoto
dihedral group in an OZVOA (= CFT type with zero degree 1 part)\cite{gnavoa1} with a positive definite invariant form. This broad class of VOAs contains all lattice type VOAs $V_L^+$ such that the even lattice $L$ is rootless, and
the Moonshine VOA $V^\natural$.   If a VOA has nontrivial degree 1
part, the order of a Miyamoto dihedral group may not be bounded in
general (for instance, a Miyamoto involution can invert a nontrivial
torus under conjugation).   See \cite{gnavoa1}.

If $L$ is rootless, it is conjectured \cite{lsy} that the above two
kinds of conformal vectors will exhaust all the conformal vectors of
central charge $1/2$ in $V_L^+$. This conjecture was proved when $L$
is a $\sqrt{2}$ times a root lattice or the Leech lattice
\cite{lsy,ls} but it is still open if $L$ is a general rootless
lattice.   The results of this paper could help settle this
conjecture, as well as provide techniques for more work on the Glauberman-Norton
theory \cite{gn}. 

The first author acknowledges financial support from the National Science Foundation,  National Cheng Kung University where the first author was a visiting distinguished professor, and the University of Michigan.
The second author acknowledges financial support from the National Science
Council of Taiwan (Grant No. 95-2115-M-006-013-MY2).

\section{Background and notations}

\begin{table}[htdp]

\caption{ \bf  Notation and Terminology}
\begin{center}
\begin{tabular}{|c|c|c|}
   \hline
\bf{Notation}& \bf{Explanation} & \bf{Examples in text}  \cr
   \hline\hline
$\Phi_{A_1},  \cdots ,\Phi_{ E_8}$ & root system of the indicated type &  Table \ref{NREE8SUM}, \ref{contain}\cr
   \hline
$A_1,  \cdots , E_8$ & root lattice for root system $\Phi_{A_1}, \dots , \Phi_{E_8}$ &  Table \ref{NREE8SUM}, \ref{contain} \cr
   \hline
$AA_1, \cdots , EE_8$ & lattice isometric to $\sqrt 2$ times the
lattice  & Table \ref{NREE8SUM} \cr
  &$A_1, \cdots,  E_8$ &\cr
    \hline
 $ann_L(S)$ &  $\{v\in L |\ (v,s)=0 \text{ for all } s\in S\}$ & \refpp{ssd}, \refpp{dih12setup}\cr
     \hline
 $\brw d$ & the Bolt-Room-Wall group, a subgroup of  & \cr &
$O(2^d,\QQ )$of shape $\ratholoex d$ & \cr \hline
$BW_{2^n}$ & the
Barnes-Wall lattice of rank $2^n$  & Table \ref{NREE8SUM},
\ref{dih8}, \refpp{dih8FJ} \cr
    \hline
 $\dih{n}{r}$ &  an NREE8SUM $M$, $N$ such that the SSD & Table \ref{NREE8SUM}, Sec. \ref{sec:F3}\cr
 & involutions $t_M, t_N$
 generate a dihedral group & \cr
 & of order $n$ and $M+N$ is of rank $r$ & \cr
       \hline
 $\dih{8}{16,X}$ &  an NREE8SUM $\dih{8}{16}$ such that&  Table \ref{NREE8SUM}, Sec. \ref{sec:F3} \cr
 & $X\cong ann_M(N)\cong ann_N(M)$ & \cr
       \hline
 $DIH_n$-theory &  the theories for $\dih{n}{r}$ for all $r$&  Sec.
\ref{sec:dih8}, \ref{sec:dih12} \cr
       \hline
$\dg L$ & discriminant group of integral lattice $L$: $\dual L/L$ & \refpp{prankdg}, \refpp{e6overaa2+4}\cr
      \hline
$HS_n$ or $D_n^+$ & the half spin lattice of rank $n$, i.e., & \refpp{not3c}\cr
   &the lattice generated by $D_n$ and $\frac{1}2(11\cdots 1)$ & \cr
      \hline
$HHS_n^+$ or $DD_n^+$ & $\sqrt{2}$ times the half spin lattice $HS_n^+$ & \refpp{invoinLambda}\cr
    \hline
IEE8 pair & a pair of $EE_8$ lattices whose sum is integral & \refpp{IEE8pair}, Sec. \ref{sec:F3} \cr
     \hline
IEE8SUM &  the sum of an IEE8 pair  & \cr
      \hline
 NREE8 pair & an IEE8 pair whose sum has no roots &  \refpp{IEE8pair}, Sec. \ref{sec:F3} \cr
       \hline
 NREE8SUM  & the sum of an NREE8 pair &  Table \ref{NREE8SUM} \cr
       \hline
$\dual L$ & the dual of the rational lattice $L$, i.e., those  &
\refpp{ssd}, \refpp{prankdg}\cr
  & elements $u$ of $\QQ \otimes L$ which satisfy $(u,L)\le \ZZ$ &  \cr
    \hline
$\Lambda$ & the Leech lattice & Sec. \ref{sec:F3}\cr
    \hline
$L^+(t), L^-(t)$ & the eigenlattices for the action of $t$& \cr
   & on the lattice $L$: $L^\vep (t):=\{x\in L| xt=\vep x \}$ & \refpp{not3c}\cr
      \hline
$Tel(L,t)$ & total eigenlattice for action of $t$&  \refpp{g1}\cr
& on
$L$;  $L^+(t)\perp L^-(t)$  & \cr
       \hline
$Tel(L,D)$ & total eigenlattice for action of an elementary &  \refpp{f/fmfnelemabel} \cr & abelian
$2$-group $D$ on $L$; $Tel(L,D):=\sum L^{\lambda}$,  &
\cr & where $\lambda \in Hom(D, \{ \pm 1\} )$ and & \cr  &  $L^\lambda=\{\a\in L|\, \a g=\lambda(g)\a \text{ for all } \g\in D\}$ & \cr
       \hline
$m^n$ & the homocyclic group $\ZZ_m^n=\ZZ_m\times \cdots \times
\ZZ_m$, & \refpp{df}, \refpp{allabouta4(1)} \cr
  &$n \text{ times}$&\cr
    \hline
$|g|, |G|$ & order of a group element, order of a group & \refpp{duala2}, Sec. \ref{sec:F3}\cr
      \hline
$O(X)$ or $Aut\, X$ & the isometry group of the lattice $X$&
\refpp{rootless} \cr \hline
\end{tabular}
\end{center}
\label{notationandterminology}
\end{table}%

\begin{table}[htdp]
\caption{ \bf  Notation and Terminology (continued)}
\begin{center}
\begin{tabular}{|c|c|c|}
   \hline
\bf{Notation}& \bf{Explanation} & \bf{Examples in text}  \cr
   \hline\hline
$O_p(G)$& the maximal normal $p$-subgroup of $G$ & \refpp{commutatoreigen}, \refpp{dgxdgy}\cr \hline
$O_{p'}(G)$& the maximal normal $p'$-subgroup of $G$,  & \refpp{elemabelsylow} \cr
\hline
$p$-rank & the rank of the maximal elementary &
\refpp{prankdg}\cr & abelian $p$-subgroup of an abelian group    & \cr
    \hline
 root & a vector of norm 2 & Sec. \ref{NREE8SUM}, \refpp{ee8intensor} \cr
      \hline
rectangular  & a lattice with an orthogonal basis & \refpp{rank3}\cr
     lattice&&\cr
       \hline
square lattice &  a  lattice isometric to some $\sqrt m \, \ZZ^n$
       & \refpp{rank3}
\cr \hline
 $Weyl(E_8)$, $Weyl(F_4)$ & the Weyl group of type $E_8$, $F_4$, etc  & \refpp{ssde8}\cr
     \hline
 $X_n$ & the set of elements of norm $n$  & \refpp{equivreln4}, \refpp{tmaps}\cr
 & in the lattice $X$ &\cr
     \hline
$X^{\perp n}$ & the orthogonal sum of $n$ copies
& Table \ref{NREE8SUM}, \ref{dih4}, \ref{dihedralpairs}\cr
& of the lattice $X$ &\cr
     \hline
$\xi$& an isometry of Leech lattice &
\refpp{defofxi}, Sec. \ref{sec:F3} \cr
& (see Notation \ref{defofxi})&\cr
    \hline
$\ZZ^n$ & rank $n$ lattice with an orthonormal basis &
\refpp{rank4det4}\cr
      \hline
\end{tabular}
\end{center}
\label{notationandterminology2}
\end{table}%

\begin{conv}\labttr{rightaction} Lattices in this article shall be rational and positive definite.  Groups and linear transformations will generally act on the right and $n$-tuples will be row vectors.
\end{conv}

\begin{de} Let $X$ be an integral lattice. For any positive integer $n$, let
$
X_n=\{ x\in L | \, (x,x)=n\}
$
be the set of all norm $n$ elements in $X$.
\end{de}

\begin{de}\labttr{determinesummand}  If $L$ is a lattice,  {\it the summand of
$L$ determined by the subset $S$} of $L$  is the intersection of $L$
with the $\QQ$-span of $S$.
\end{de}

\begin{nota}\labttr{ssd} Let $X$ be a subset of Euclidean space.
Define $t_X$ to be the orthogonal transformation which is $-1$
on $X$ and is 1 on $X^\perp$.
\end{nota}

\begin{de}\labttr{rssd}
A sublattice $M$ of an integral lattice $L$ is {\it RSSD (relatively semiselfdual)}
if and only if $2L\le M + ann_L(M)$. This implies that $t_M$
maps $L$ to $L$ and is equivalent to this property when $M$ is
a direct summand.

The property that $2\dual M \le M$ is called {\it SSD
(semiselfdual)}. It implies the RSSD property, but the RSSD
property is often more useful.  For example, if $M$ is RSSD in $L$ and $M\le J\le L$, then $M$ is RSSD in $J$, whence the involution $t_M$ leaves $J$ invariant.
\end{de}

\begin{ex}\labttr{exssd}
An example of a SSD sublattice is $\sqrt 2 U$, where $U$ is a
unimodular lattice. Another is the family of Barnes-Wall
lattices.
\end{ex}

\begin{lem}
\labtt{ssd=rssd} If the sublattice $M$ is a direct summand of
the integral lattice $L$ and $(det(L),det(M))=1$, then SSD and
RSSD are equivalent properties for $M$.
\end{lem}
\pf It suffices to assume that $M$ is RSSD in $M$ and prove
that it is $SSD$.

Let $V$ be the ambient real vector space for $L$ and define
$A:=ann_V(M)$. Since $(det(L),det(M))=1$, the natural image of
$L$ in $\dg M$ is $\dg M$, i.e., $L+A=\dual M +A$ \refpp{dgxdgy}. We have
$2(L+A)=2(\dual M+A)$, or $2L+A=2\dual M + A=M\perp A$. The left side is
contained in $M+A$, by the RSSD property. So, $2\dual M \le
M+A$. If we intersect both sides with $ann_V(A)$, we get
$2\dual M \le M$. This is the SSD property.
\eop

\begin{lem}\labtt{subrssd}
Suppose that $L$ is an integral lattice and
$N\le M\le L$ and both $M$ and $N$ are RSSD in $L$.   Assume that $M$ is a direct summand of $L$.
Then $ann_M(N)$ is an RSSD sublattice of $L$.
\end{lem}
\pf
This is easy to see on the level of involutions.  Let $t, u$ be  the involutions
associated to $M, N$.   They are in $O(L)$ and they
commute since $u$ is the identity on $ann_L(M)$, where $t$ acts as the scalar 1,
 and since $u$ leaves invariant $M=ann_L(ann_L(M))$, where $t$ acts as the
scalar $-1$.  Therefore, $s:=tu$ is an involution. Its negated sublattice
$L^-(s)$ is RSSD \refpp{ssd}, and this is $ann_M(N)$.
\eop

\begin{de}\labtt{IEE8pair}
An {\it IEE8 pair} is a pair of sublattices $M, N\cong EE_8$ in a
Euclidean space such that $M+N$ is an integral lattice. If $M+N$ has
no roots, then the pair is called an {\it NREE8 pair}. An {\it IEESUM} is  the
sum of an IEE8 pair and an {\it NREE8SUM} is the sum of an NREE8 pair.
\end{de}

\begin{lem}\labttr{nocommonfp}
We use Definition \ref{ssd}. Let $M$ and $N$ be RSSD in  an integral lattice
$L=M+N$. A vector in $L$ fixed by both $t_M$ and $t_N$ is 0.
\end{lem}
\pf
 We use $L=M+N$.
 If we tensor $L$ with $\QQ$,  we have complete reducibility for the action  of
$\la t_M, t_N\ra$.  Let $U$ be the fixed point space for $\la t_M, t_N \ra$ on
$\QQ \otimes L$.  The images of $M$ and $N$ in $U$ are 0, whence $U=0$.
\eop

\section{Tensor products} 

\begin{de}
Let $A$ and $B$ be integral lattices with the inner products $(\ ,\ )_A$  and
$(\ ,\ )_B$, respectively.
\textit{The tensor product of the lattices} $A$ and $B$ is defined to be the
integral lattice which is isomorphic to $A\otimes_\ZZ B$ as a $\ZZ$-module and
has the inner product given by
\[
(\a\otimes \b, \a' \otimes \b') = (\a,\a')_A \cdot (\b,\b')_B, \quad \text{ for any } \a,\a'\in A,\text{ and } \b,\b'\in B.
\]
We simply denote the tensor product of the lattices $A$ and $B$ by $A\otimes B$.
\end{de}

\begin{lem}\labtt{tensorwitha2}  Let $D:=\la t, g\ra$ be a dihedral group of order 6, generated by an involution $t$ and element $g$ of order 3.   Let $R$ be a rational lattice on which $D$ acts such that $g$ acts fixed point freely.  Suppose that $A$ is a sublattice of $R$ which satisfies  $at=-a$ for all $a\in A$.
Then
\begin{enumerate}
\item[(i)] $A\cap Ag=0$; so $A+Ag=A\oplus Ag$ as an abelian group.

\item[(ii)] $A+Ag$ is isometric to
$A\otimes B$, where $B$ has Gram matrix
$\begin{pmatrix}1&-\half \cr -\half & 1\end{pmatrix}$.

\item[(iii)] Furthermore $ann_{A+Ag}(A)= A(g-g^2)
\cong \sqrt 3 A$.
\end{enumerate}
\end{lem}
\pf
(i)
Take $a\in A$ and suppose $a=a'g$. Since $at=-a$, we have
$-a=at=a'gt=a'tg^2=-a'g^2$. That means $a=a'g^2$ and $a=ag$. Thus $a=0$ since $g$ acts fixed point freely on $R$.

For any $x, y \in R$, we have $0=(x,0)=(x,y+yg+yg^2)=(x,y)+(x,yg)+(x,yg^2)$.
Now, take $x, y \in A$.
We have $(x,yg)=(xt, ygt)=(-x,ytg^2)=(-x,-yg^2)=(x,yg^2)$.  We conclude that
$(x,yg)=(x,yg^2)=-\half (x,y)$.

Let bars denote images under the quotient $\ZZ \la g \ra \rightarrow \ZZ \la g \ra /(1+g+g^2)$

We use the linear monomorphism $A\otimes \overline {g^i} \rightarrow R$ where $\ZZ \la g \ra /(1+g+g^2)$ has the bilinear form which take value 1 on a pair $\overline {g^i}, \overline {g^i}$ and value $-\half$ on a pair $\overline {g^i}, \overline {g^j}$ where $j=i\pm 1$.  This proves (ii)

For (iii), note that $\psi: x\mapsto xg-xg^2$ for $x \in A$ is a scaled isometry and $Im\,\psi$ is a direct summand of $Ag\oplus Ag^2$. Note also that $A\oplus Ag= Ag\oplus Ag^2= Im\,\psi \oplus Ag$. Thus we have
\[
Im\, \psi \le ann_{A+Ag}(A) \le  Im\,\psi \oplus Ag.
\]
By Dedekind law, $ann_{A+Ag}(A)=Im\,\psi +(ann_{A+Ag}(A)\cap Ag)$. Since $(x,yg)=-\frac{1}2 (x,y)$, $ann_{A+Ag}(A)\cap Ag=0$ and we have $Im\, \psi = ann_{A+Ag}(A)$ as desired.
\eop

\begin{lem}\labtt{tensorminvec} Suppose that $A, B$ are lattices, where $A\cong A_2$.  The minimal vectors of $A\otimes B$ are just $u \otimes z$, where $u$ is a minimal vector of $A$ and $z$ is a minimal vector of $B$.
\end{lem}
\pf
Let $u$ be a minimal vector of $A$.
The minimal vectors of $\ZZ u \otimes B$ have the above shape.
Let $u'$ span $ann_A(u)$. Then $(u',u')=6$ and $|A: \ZZ u+\ZZ u'|=2$.
The minimal vectors of $(\ZZ u  \perp \ZZ u')\otimes B$ have the above shape.
Now take a vector $w$ in $A\otimes B \setminus (\ZZ u  \perp \ZZ u')\otimes B$.
It has the form $pu \otimes x + qu'\otimes y$, where $p, q \in \half +\ZZ$ and $p+q\in \ZZ$. The norm of this vector is therefore
$2p^2 (x,x)+6q^2 (y,y)$.  A necessary condition that $w$ be a minimal vector in $A\otimes B$ is that each of $x, y$ be minimal in $B$ and $p,q\in \{\pm \half\}$.

Define $v:=\half u+\half u'$.
Then $u', v$ forms a basis for $A$. We have $w= v \otimes x + \half u' \otimes (y-x)$.  Since $w\in A\otimes B$, $y-x\in 2B$. Suppose $y-x=2b$.
If $b=0$, done, so assume that $b\ne 0$.
In case $x, y$ are minimal, $(y-x,y-x)=4(b,b)\geq 4(x,x)$ and thus $-2(x,y) \geq 2(x,x)$. This implies $x=-y$ and then $w= (v-u')\otimes x$ as required.
\eop

\begin{nota}\labtt{minvector}
For a lattice $L$, let $\mv L$ be the set of minimal vectors.
\end{nota}
\begin{lem}\labtt{ee8intensor}
 We use the notations of
\refpp{tensorminvec}.  If $B$ is a root lattice of an indecomposable
root system and $rank(B)\ge 3$, the only sublattices of $A\otimes B$ which are isometric
to $\sqrt 2 B$ are the
$u \otimes B$, for $u$ a minimal vector of $A$.
\end{lem}
\pf
Let $S$ be a sublattice of $A\otimes B$ so that $S\cong \sqrt{2} B$.  Then $S$
is spanned by $\mv S$, which by \refpp{tensorminvec} equals $M_u \cup
M_v \cup M_w$, where $u, v, w$ are pairwise nonproportional minimal
vectors of $A$ which sum to 0 and where $M_t:=(t\otimes B) \cap \mv S$, for
$t=u, v, w$.  Note that $(u,v)=(v,w)=(w,u)=-1$.

We suppose that $M_u$ and $M_v$ are nonempty and seek a contradiction.
Take $b, b'\in B$ so that $u \otimes b \in M_u, v\otimes b' \in M_v$.
Then $(u \otimes b, v\otimes b')=(u,v)(b,b') = -(b,b')$.
Since $S$ is doubly even, all such $(b,b')$ are even.

We claim that all such $(b,b')$ are 0.

Assume that some such $(b,b')\ne 0$.  Then, since $b, b'$ are roots,  $(b,b')$ is $\pm 2$ and $b=\pm b'$.   Then
$u\otimes b, v\otimes b \in S$, whence $w\otimes b \in S$.  In other words, $A\otimes b \le S$.

Since $rank(S)=rank(B)\ge 3$, $S$ properly contains $A\otimes b$.  Since $S$ is generated by its minimal vectors and the root system for $B$ is connected, $S$ contains some $t\otimes d$ where $d \in \mv B$ and $(d,b)\ne 0$.  It follows that $(d,b)=\pm 1$.  Take $t'\in \mv A$ so that $(t,t')=\pm 1$.  Then $(t\otimes d, t'\otimes b)=\pm 1$, whereas $S$ is doubly even, a contradiction.  The claim follows.

The claim implies that $M_u$ and $M_v$ are orthogonal.  Similarly, $M_u, M_v, M_w$ are pairwise orthogonal, and at least two of these are nonempty.   Since $\mv S$ is the disjoint union of $M_u, M_v, M_w$, we have a contradiction to indecomposability of the root system for $B$.
\eop

\section{Uniqueness}

\begin{thm}\labtt{uniquenessofl}
Suppose that $L$ is a free abelian group and that $L_1$ is a
subgroup of finite index.  Suppose that $f:{L}_1 \times
{L}_1\rightarrow K$ is a $K$-valued bilinear form, where $K$ is
an abelian group so that multiplication by $|L:{L}_1|$ is an
invertible map on $K$.  Then $f$ extends uniquely to a $K$-valued
bilinear form $L\times L \rightarrow K$.
\end{thm}
\pf Our statements about bilinear forms are equivalent to statements
about linear maps on tensor products. We define $A:=L_1\otimes L_1, B:=L \otimes L$ and $C:=B/A$.   Then $C$ is finite and is
annihilated by $|L:L_1|^2$.
From $\ses ABC$, we get the long exact sequence $0\rightarrow
Hom(C,K) \rightarrow Hom(B,K) \rightarrow Hom(A,K) \rightarrow Ext^1(C,K)\rightarrow
\cdots $.   Each of the  terms $Hom(C,K)$ and $Ext^1(C,K)$ are 0 because
they are annihilated by $|C|$ and multiplication by $|C|$ on $K$ is an
automorphism. It follows that the restriction map from $B$ to $A$ gives
an isomorphism $Hom(B,K)\cong Hom(A,K)$. \eop

\begin{rem}\labttr{applicuniquenessofl}
We shall apply \refpp{uniquenessofl} to $L=M+N$ when we determine
sufficient information about a pairing ${M}_1\times
 {N}_1 \rightarrow \QQ$, where ${M}_1$ is a finite index
sublattice of $M$ and ${N}_1$ is a finite index sublattice of
$N$. The pairings $M\times M\rightarrow \QQ$ and  $N\times
N\rightarrow \QQ$ are given by the hypotheses $M\cong N \cong EE_8$,
so in the notation of
\refpp{uniquenessofl} we take $L_1=M_1+N_1$.
\end{rem}

\begin{rem}
We can determine all the lattices in the main theorem by explicit glueing. However, it is difficult to prove the rootless property in some of those cases. In Appendix F, we shall show that all the lattices in Table 1 can be embedded into the Leech lattice except $\dih{4}{15}$. The rootless property follows since the Leech lattice has no roots.  The proof that $\dih{4}{15}$ is rootless will be included at the end of Subsection \ref{sec:dih4}.
\end{rem}

\section{$\Dih{4}$ and $\Dih{8}$ theories}

\subsection{$\Dih{4}$: When is $M+N$ rootless?}\label{sec:dih4}

\begin{nota}\labttr{nota-dih4}
Let $M, N$ be $EE_8$ lattices such that the dihedral group $D:=\la t_M, t_N\ra$ has order 4.
Define $F:=M\cap N$,
$P:=ann_M(F)$ and $Q:=ann_N(F)$.
\end{nota}

\begin{rem}\labttr{remdih4}
Since $t_M$ and $t_M$ commute, $D$ fixes each of $F$, $M$, $N$, $ann_M(F)$, $
ann_N(F)$.  Each of these may be interpreted as eigenlattices since
$t_M$ and $t_N$ have common negated space $F$, zero common fixed space, and  $t_M, t_N$ are
respectively $-1, 1$ on $ann_M(F)$ and $t_M, t_N$ are respectively $1,
-1$ on $ann_N(F)$.  Since $L=M+N$, $D$ has only 0 as the fixed point
sublattice (cf. \refpp{nocommonfp}).  Therefore, the elementary abelian group
$D$ has total
eigenlattice $F \perp ann_M(F) \perp ann_N(F)$.
Each of these summands is RSSD as a sublattice of $L$, by
\refpp{subrssd}.
It follows that $\drtp{M\cap N}$ is an RSSD sublattice in $\drt{M}$ and in $\drt{N}$. Since $\drt{M}\cong \drt{N}\cong E_8$, we have that $\drtp{M\cap N}$ is an SSD sublattice in $\drt{M}$ and in $\drt{N}$ (cf. \refpp{ssd=rssd}).
\end{rem}

\begin{prop}\labtt{rootless}
If  $M+N$ is rootless, $F$ is one of 0, $AA_1, AA_1\perp AA_1, DD_4$.
Such sublattices in $M$ are unique up to the action of $O(M)$.
\end{prop}

\pf This can be decided by looking at cosets of $P+Q+F$ in $M+N$.
A glue vector will have nontrivial projection to two or three of
$\Rspan{P}, \Rspan{Q}, \Rspan{F}$.

Since $F$ is a direct summand of $M$ by \refpp{mnfsummands} and $\drt{F}$ is an SSD  by \refpp{remdih4}, we have  $\drt{F}\cong 0, A_1, A_1\perp A_1, A_1\perp A_1\perp A_1, A_1\perp A_1\perp A_1\perp A_1, D_4, D_4\perp A_1, D_6, E_7$ and $ E_8$ by \refpp{ssde8}. If $\drt{F}=E_8$, then $t_M=t_N$ and $D:=\la t_M, t_N\ra$ is only a cyclic group of order $2$. Hence, we can eliminate $\drt{F}=E_8$.

Now we shall note that in each of these cases, $F\perp P$ contains a
sublattice $A\cong AA_1^8$ such that $F\cap A\cong AA_1^{k}$ and $A\cap
P\cong AA_1^{8-k}$, where $k=rank\, F$. We use an orthogonal basis of $A$ to identify $M/A$ with a code.  Since $M\cong EE_8$, this code is the Hamming $[8,4,4]$ binary code
$H_8$.  Let $\varphi: M/A \to H_8$ be such an identification. Then $
\varphi((F\perp P)/A)$ is a linear subcode of $H_8$.

Next, we shall show that $\dual{\left( \drt{F}\right )}$ contains a
vector $v$ of norm $3/2$ if
$\drt{F}\cong A_1\perp A_1\perp A_1, A_1\perp A_1\perp A_1\perp A_1, D_4\perp A_1, D_6$ or $E_7$.

Recall that if $A_1=\ZZ\a, (\a,\a)=2$, then $A_1^*=\frac{1}2 \ZZ\a$ and
$\frac{1}2 \a\in A_1^*$ has norm $\frac{1}2$. Since $(A_1^{\perp k})^* =
(A_1^*)^{\perp k}$, $(A_1^{\perp k})^*$ contains a vector of norm $3/2$ if
$k\geq 3$.

We use the standard model
\[
D_n= \{ (x_1, x_2, \dots, x_n)\in \ZZ^n|\, x_1+\cdots +x_n\equiv 0\mod 2\}.
\]
Then  $\frac{1}2 (1,\dots,1) \in D_n^*$ and its norm is $\frac{1}4 n$. Therefore, there exists vectors of norm $3/2$ in $(D_4\perp A_1)^*= D_4^*\perp A_1^*$ and $D_6^*$. Finally, we recall that
$E_7^*/E_7\cong \ZZ_2$ and the non-trivial coset is represented by a vector
of norm $3/2$.

\medskip

Now suppose $\drt{F}\cong A_1\perp A_1\perp A_1, A_1\perp A_1\perp A_1\perp A_1, D_4\perp A_1, D_6$ or $E_7$  and
let $\gamma\in 2 F^*$ be a vector of norm $3$.  Since $F\cap A\cong
AA_1^k$, $F\cap A$ has a basis $\{\a_1, \dots,\a_k\}$ such that
$(\a_i,\a_j)=4\delta_{i,j}$. Then 
\[
2F^* < 2(F\cap A)^*=\mspan_\ZZ\{\frac{1}{2}\sum_{i=1}^k a_i \a_i|\, a_i\in
\ZZ\}.
\]
Thus, by replacing some basis vectors by their negatives, we have  $\gamma= \frac{1}{2}(\a_{i_1}+\a_{i_2}+\a_{i_3})$ for some
$1\leq i_1, i_2, i_3 \leq k$.

Since the image of the natural map $M\to \dg{F}$ is $2\dg F$, there exists a vector $u\in M$
such that the projection of $u$ to  $\Rspan{F}$ is $\gamma$.

Now consider the image of $u+A$ in $H_8$ and study the projection of the codeword $\varphi(u+A)$ to the first $k$ coordinates. Since $\gamma= \frac{1}{2}(\a_{i_1}+\a_{i_2}+\a_{i_3})$, the projection of $\varphi(u+A)$ to the first $k$ coordinates has weight $3$.

If $k=rank\, F\geq 4$, then the projection of $(1,\cdots,1)$ to the first $k$ coordinates has weight $k\geq 4$. Thus, $\varphi(u+A)\neq (1,\dots,1)$ and hence $\varphi(u+A)$ has weight $4$ since $\varphi(u+A)\in H_8$ .

If $k=3$, then $F\cong AA_1^3$ and
$P\cong AA_1\perp DD_4$.  Let $K\cong DD_4$ be 
an orthogonal direct summand of $P$ and let $\ZZ \a = ann_P(K)\cong AA_1$.  Note that $ F\perp \ZZ\a < ann_ M(K)\cong DD_4$ and $F=\ZZ\a_1\perp \ZZ\a_2\perp \ZZ\a_3$.  Thus, $\frac{1}2(\a_1+\a_2+\a_3+\a)=\frac{1}2(\g+\a)\in ann_ M(K)< M$  and it has norm $4$. Therefore, we may assume $\varphi(u+A)$ has weight $4$  and $u$ is a norm $4$ vector.

Similarly,
there exists a norm $4$ vector $w\in N$ such that the projection of $w$ in
$\Rspan{F}$ is also $\gamma$. Then $u-w\in L=M+N$ but
$(u,w)=(\gamma,\gamma)=3$ and hence $u-w\in L$ is a root, which
contradicts the rootless property of $L$. Therefore, only the remaining cases
occur, i.e., $F\cong 0, AA_1, AA_1\perp AA_1, DD_4$.  \eop

\begin{table}[bht]
\caption{ \bf  $\Dih{4}$: Rootless cases}
\begin{center}
\begin{tabular}{|c|c|c|c|}
 \hline $M\cap N$ & $P\cong Q$ & $dim(M+N)$ & Isometry type of L\cr
 \hline \hline 0&$EE_8$& 16&$\cong EE_8\perp EE_8$\cr
 \hline $DD_4$&$DD_4$&12& $\geq DD_4\perp DD_4\perp DD_4$\cr
 \hline $AA_1$& $EE_7$&15& $\geq  AA_1\perp EE_7\perp EE_7$\cr
 \hline $AA_1\perp AA_1$& $DD_6$&14& $\geq AA_1\perp AA_1\perp DD_6\perp DD_6$ \cr
\hline
\end{tabular}
\end{center}
\label{dih4}
\end{table}%

\begin{rem}\labttr{rootlessverification}
Except for the case $F=M\cap N\cong AA_1$, we shall show in Appendix F that all
cases in Proposition \ref{rootless} occur inside the Leech lattice $\Lambda$ .
The rootless property of $L=M+N$ then follows from the rootless property of
$\Lambda$.  The rootless property for the case $F=M\cap N\cong AA_1$ will be
shown in the next proposition.
\end{rem}

\begin{prop}\labtt{rootlessdih415} [Rootless property for $\dih{4}{15}$]
If $F=M\cap N\cong AA_1$, then $P\cong Q\cong EE_7$ and $L=M+N$ is rootless.
\end{prop}

\pf
We shall use the standard model for the lattice $E_7$, i.e.,
\[
E_7=\left \{(x_1, \dots, x_8)\in \ZZ^8\left|\, {{\text{all } x_i\in \ZZ \text{ or all } x_i\in \frac{1}2+\ZZ}\atop  {\text{ and }
x_1+\cdots+x_8=0}}\right.  \right \}
\]
The dual lattice is $E_7^* = E_7 \cup (\gamma + E_7)$, where $\gamma
= \frac{1}4 (1,1,1, 1,1,1, -3,-3)$. Recall that the minimal weight
of $E_7^*$ is $3/2$ \cite[p.125]{cs}.

If $F=M\cap N=AA_1=\ZZ\a$, then it is clear that $P\cong Q\cong
EE_7$. In this case, $M=\mspan_\ZZ\{ F+P, \frac{1}2 \a + \xi_M\}$
and $N= \mspan_\ZZ\{ F+Q, \frac{1}2 \a + \xi_N\} $ for some
$\xi_M\in P^*$ and $\xi_N\in Q^*$ with $(\xi_M,\xi_M)=(\xi_N,
\xi_N)=3$. Therefore,
\[
L=M+N=\mspan_\ZZ\{ F+P+Q, \frac{1}2\a+\xi_M,\frac{1}2\a+ \xi_N\}.
\]
Take $\b\in L=M+N$. If $\b \in F+P+Q$, then $(\b, \b)\geq 4$. Otherwise, $\b$  will have nontrivial projection to two or three of $\Rspan{P}, \Rspan{Q}, \Rspan{F}$.  Now note that the projection of
$L$ onto $\Rspan{P}$ is $\mspan_\ZZ\{ P, \xi_M\}\cong \sqrt{2}E_7^*$ and the projection of
$L$ onto $\Rspan{Q}$ is $\mspan_\ZZ\{ Q, \xi_N \}\cong \sqrt{2}E_7^*$. Both of them have minimal norm $3$.
On the other hand, the projection of $L$ onto $\Rspan{F}$ is $\ZZ \frac{1}2 \a$, which has minimal norm $1$. Therefore, $(\b,\b)\geq 1 + 3=4$ and so $L$ is rootless.
\eop

\subsection{$\Dih{8}$} \label{sec:dih8}

\begin{nota}\labttr{dih8setup}   Let $t:=t_N, u:=t_N$, and $g:=tu$,  which has order 4.
Define $z:=g^2, t':= tz$ and $u':=uz$.  We define $F:=Ker_L(z-1)$ and $J:=Ker_L(z+1)$.
\end{nota}

By Lemma \ref{involonlattice},  $L/(F\perp J)$ is an elementary
abelian 2-group of rank at most $min\{rank(F),rank(J)\}$. We have two
systems $(M, t, Mg, t')$ and $(N,u,Ng,u')$ for which the $\Dih{4}$ analysis applies.

\begin{nota}\labtt{dih8:ljf}
If $X$ is one of $M, N$, we denote by $L_X, J_X,
F_X$ the lattices $L:=X+ Xg, J, F$ associated to the pair
$X, Xg$, denoted $``M"$ and $``N"$ in the $\Dih{4}$ section.
\end{nota}

\subsubsection{$\Dih{8}$: What is $F$?}\label{sec:dih8F}

We now determine $F$.

\begin{rem}\labttr{rankdeterminesf}
It will turn out that
the two systems $(M, t, Mg, t')$ and $(N,u,Ng,u')$
have the same $\Dih{4}$ types (cf.
Table \ref{dih4}). Also,  we shall prove that $rank(F_X)$ determines $F_X$, hence also
determines $J_X$, for $X=M, N$.
\end{rem}

\begin{lem}\labtt{squareofg-1}
Let $f= g$ or $g^{-1}$. Then
(i) As endomorphisms of $J$, $f^2=-1$, $(f-1)^2=-2f$.  For $x, y \in J$, $(x(f-1),y(f-1))=2(x,y)$.

(ii) $(M\cap J, (M\cap J)f)=0$ and $(N\cap J, (N\cap J)f)=0$.

(iii) For $x, y \in M$ or $x, y \in N$, $(x,y(f-1))=-(x,y)$.
\end{lem}
\pf
(i) As endomorphisms of  $J$, $(f-1)^2=f^2-2f+1=-2f$.

(ii) We take $x, y \in M\cap J$ (the argument for $x, y \in N\cap J$ is similar).

We have
$(x,yf)=(xt,yft)=(-x,ytf^{-1})=(-x,-yfz)=(-x,yf)=-(x,yf)$, whence
$(x,yf)=0$.

(iii) We have $(x,y(f-1))=(x,yf)-(x,y)=-(x,y)$.
\eop

\begin{lem}\labtt{conjcalc}(i) In $\QQ \otimes End(J)$, $(g^{-1}-1)^{-1}t(g^{-1}-1)=u$.

(ii) $(M\cap J)(g^{-1}-1)\le N\cap J$ and
$(N\cap J)2(g^{-1}-1)^{-1}\le M\cap J$.

(iii) $rank(M\cap J)=rank(N\cap J)$.

(iv) $rank(F_M)=rank(F_N)$.
\end{lem}
\pf 
We use the property that $g^{-1}$ acts as $-g$ on $\QQ\otimes J$.  We also abuse notation by identifying elements of $\QQ [D]$ with their images in $End(\QQ \otimes J)$.  For example, $(g^{-1}-1)$ is not an invertible element of $\QQ [D]$, though its image in
$End(\QQ \otimes J)$ is invertible.

For (i), observe
that $(g^{-1}-1)^2=-2g^{-1}$, so that
$g^{-1}-1$ maps $J$ to $J$ and has zero kernel.
Secondly,
$2(g^{-1}-1)^{-1}$ maps $J$ to $J$ and has zero kernel.

The
equation $(g^{-1}-1)^{-1}t(g^{-1}-1)=u$
in $\QQ \otimes End(J)$
is equivalent to
$t(g^{-1}-1)=(g^{-1}-1)u$ which is the same as $(g-1)t=(g^{-1}-1)u$
or $tut-t=-gu-u=-tuu-u=-t-u$, which is true since $tut=-u$.

The statement (ii) follows since in a linear representation of a
group, a group element which conjugates one element to a second one
maps the eigenspaces of the two elements correspondingly.  Here, this means $g^{-1}-1$
conjugates $t$ to $u$, so that $g^{-1}-1$ maps $\QQ \otimes (M\cap J)$ to $\QQ \otimes (N \cap J)$.
Since $g^{-1}-1$ maps $J$ into $J$ (though not onto $J$), $g^{-1}-1$ maps the direct summand $M\cap J$ into the direct summand $N\cap J$.

For (iii), observe that we have monomorphisms  $M\cap
J\rightarrow N\cap J \rightarrow M\cap J$ and $N\cap J\rightarrow
M\cap J \rightarrow N\cap J$ by use of $g^{-1}-1$ and $2(g^{-1}-1)^{-1}$.
Therefore,
(iv) follows from (iii).
\eop

\begin{lem}\labtt{rankjcapmeven}  Suppose that
$det(J\cap M)det(J\cap N)$ is the square of an integer
(equivalently,  that $det(F_M)det(F_N)$ is the square of an
integer). Then $rank(J\cap M)=rank(J\cap N)$ is even.
\end{lem}
\pf Note that $rank(M\cap J)=rank(N\cap J)$ by \refpp{conjcalc}.
 Let $d:=det(J\cap M)$ and $e:=det(J\cap N)$ and let $r$ be  the
common rank of $M\cap J$ and $N\cap J$.  First note that $(M\cap
J)(g^{-1}-1)$ has determinant $2^r d$ and second note that $(M\cap
J)(g^{-1}-1)$ has finite index, say $k$, in $N\cap J$.  It follows
that $2^r d = k^2 e$.  By hypothesis, $de$ is a perfect square.
Consequently, $r$ is even. \eop

\begin{coro}\labtt{rank} $rank(F_M)=rank(F_N)$ is even.
\end{coro}
\pf We have $rank(F)+rank(M\cap J)=rank(M)=8$ and similarly
for $N$. Since $rank\, F_M=rank \, F_N$,  we have  $F_M\cong F_N$ by
\refpp{rootless} and hence $\det F_M \det F_N= (\det F_M)^2$ is a
square. Now use \refpp{rankjcapmeven}. \eop

\begin{prop}\label{possfmfn}
If $L=M+N$ is rootless, then  $F_M\cong F_N \cong 0, AA_1\perp AA_1$ or $DD_4$.
Moreover, $M\cap J \cong N\cap J$.
\end{prop}
\pf
Since by \refpp{rank}, $rank(F_M)=rank(F_N)$ is even, Proposition \ref{rootless} implies that $F_M\cong F_N \cong 0, AA_1\perp AA_1$ or $DD_4$.
It is well-known that there is one orbit of $O(E_8)$ on the family of sublattices which have a given one of the latter isometry types.
It follows that
$ M\cap J=ann_M(F_M)\cong ann_N(F_N)\cong N\cap J$.
\eop

\subsubsection{$\Dih{8}$: Given that $F=0$, what is $J$?}\label{dih8FJ}

By Proposition \ref{possfmfn}, when $L$ is rootless, $F_M\cong F_N \cong 0, AA_1\perp AA_1$ or $DD_4$.
We now consider each case for $F_M$ and $F_N$ and determine the possible pairs $M, N$.
The conclusions
are listed in Table \ref{dih8}.

\begin{table}[htdp]
\caption{ \bf  $\Dih{8}$ which contains a rootless $\Dih{4}$ lattice}
\begin{center}
\begin{tabular}{|c|c|c|c|c|c|}
   \hline $F_M\cong F_N$ &$F_M\cap F_N$ & $rank(M+N)$ & $M+N$  & Isometry type\cr
                         &               &            &integral? roots ?&if rootless \cr
 \hline
 0& 0& 16 &  rootless& $\cong BW_{16}$\cr
 \hline
 \hline
 $AA_1^{\perp 2}$&0&16& non-integral &\cr
 \hline
 $AA_1^{\perp 2}$ &$AA_1$ &15& non-integral&\cr
 \hline
 $AA_1^{\perp 2}$ &$AA_1^{\perp 2}$& 14& non-integral&\cr
 \hline
 $AA_1^{\perp 2}$ &$ 2A_1$& 15& has roots &\cr
 \hline
 \hline
 $DD_4$& $0$ & $16$ &  rootless & $\geq DD_4^{\perp 2}\perp  EE_8$\cr
 \hline
 $DD_4$& $AA_1$ & $15$ &  rootless & $\geq AA_1^{\perp 7}\perp EE_8$\cr
 \hline
 $DD_4$& $2A_1$ & $15$ & non-integral &\cr
 \hline
 $DD_4$& $AA_1^{\perp 2}$ &$14$& has roots &\cr
 \hline
 $DD_4$& $AA_1\perp 2A_1$ & $14$ & non-integral&\cr
 \hline
 $DD_4$& $AA_1^{\perp 3}$ & $13$& has roots&\cr
 \hline
 $DD_4$& $AA_3$ & $13$ & non-integral&\cr
 \hline
 $DD_4$& $DD_4$ & $12$ & has roots &\cr
 \hline
\end{tabular}
\end{center}
\label{dih8}
\end{table}%

\medskip

\begin{prop}
If $F_M= F_N=0$, then $L=M+N$ is isometric to the Barnes-Wall lattice
$BW_{16}$.
\end{prop}
\pf The sublattice $M':=Mt_N$ is the $1$-eigenspace for $t_M$ and so
$M+M'=M\perp M'$.  Consider how $N$ embeds in $\dual {(M+M')}=\half
(M+M')$.  Let $x\in N\setminus (M+M')$ and let $y\in \half  M, y'\in
\half M'$ so that $x=y+ y'$.  We may replace $y, y'$ by members of
$y+M$ and $y'+M'$, respectively, which have least norm. Both $y, y'$
are nonzero. Their norms are therefore one of 1, 2, by a property of
the $E_8$-lattice.
 Since $(x,x)\ge  4$, $y$ and $y'$ each has norm 2.  It
follows that the image of $N$ in $\dg {M}$ is totally singular in the
sense that all norms of representing vectors in $\dual M$ are
integers.  A similar thing is true for the image of $N$ in $\dg {M'}$.
It follows that these images are elementary abelian groups which have
ranks at most 4.  On the other hand, diagonal elements of the
orthogonal direct sum $M\perp M'$ have norms at least 8, which means
that $N\cap (M+M')$ contains no vectors in $N$ of norm 4.   Therefore,
$N/(N\cap (M+M'))$ is elementary abelian of rank at least 4.  These two inequalities imply that the rank is 4.
The action of
$t_N$ on this quotient is trivial.  We may therefore use the uniqueness
theorem of \cite{bwy} to prove that $M+N$ is isometric to the
Barnes-Wall lattice $BW_{16}$. \eop
\medskip

\subsubsection{$\Dih{8}$: Given that $F\cong AA_1\perp AA_1$, what is $J$?}
\begin{prop}
If $F_M\cong F_N\cong AA_1\perp AA_1$, $M+N$ is  non-integral or
has a root.
\end{prop}
\pf
If $F_M\cong F_N\cong AA_1\perp AA_1$, then $M\cap J\cong N\cap J\cong DD_6$.  We shall first
determine the structure of $M\cap J + N\cap J$.

Let $h=g^{-1}$. Then, by Lemma \ref{conjcalc}, we have $(M\cap
J)(h-1)\leq  N\cap J$. Since $( M\cap J, (M\cap J)h) =0$, $(M\cap
J)(h-1)\cong 2D_6$ and $\det ( (M\cap J)(h-1))=2^8$. Therefore, $|N\cap J: (M\cap J)(h-1)|=(2^{rank\, N\cap J})^{1/2}=2^3$.

\medskip

Let $K=(M\cap J)(h-1)$. Then, by \refpp{sqrt2D}, there exists a
subset $\{ \eta_1, \dots, \eta_6\}\subset N\cap J$ with $(\eta_i,
\eta_j)=4 \delta_{i,j}$ such that
\[
K=\mathrm{span}_\ZZ\left \{ (\eta_i\pm \eta_j)\left|\, i,j =1,2,
\dots,6 \right.\right \}
\]
and
\[
N\cap J= \mathrm{span}_\ZZ\left \{ \eta_1, \eta_2,\eta_4, \eta_6,
\frac{1}2(-\eta_1+\eta_2-\eta_3+\eta_4),
\frac{1}2(-\eta_3+\eta_4-\eta_5+\eta_6) \right \}.
\]

By computing the Gram matrix, it is easy to show that $\{
\eta_1+\eta_2, -\eta_1+\eta_2, -\eta_2+\eta_3,-\eta_3+\eta_4,
-\eta_4+\eta_5, -\eta_5+\eta_6\}$ forms a basis of $K=(M\cap
J)(h-1)\cong 2D_6$.   Now let
\[
\begin{split}
\a_1=(\eta_1+\eta_2)(h-1)^{-1},\ \ \quad
\a_2=(-\eta_1+\eta_2)(h-1)^{-1},\quad
\a_3=(-\eta_2+\eta_3)(h-1)^{-1},\\
\a_4=(-\eta_3+\eta_4)(h-1)^{-1},\quad
 \a_5=(-\eta_4+\eta_5)(h-1)^{-1}, \quad
\a_6=(-\eta_5+\eta_6)(h-1)^{-1}.
\end{split}
\]
Then $\{\a_1,\a_2, \a_3, \a_4,  \a_5, \a_6\}$ is a basis of $M\cap
J$. Moreover, $(\a_1, \a_2)=0$,  $(\a_1, \a_3)=-2$, $(\a_i,
\a_{i+1})=-2$ for $i=2,\dots,5$.

By the definition, we have $\eta_2=-\frac{1}2 (\a_1+\a_2)(h-1),
(-\eta_1+\eta_2-\eta_3+\eta_4)= (\a_2+\a_4)(h-1)$ and
$(-\eta_3+\eta_4-\eta_5+\eta_6)= (\a_4+\a_6)(h-1)$. Hence,
\[
 N\cap J= \mathrm{span}_\ZZ\left \{ K, \frac{1}2 (\a_1+\a_2)(h-1), \frac{1}2 (\a_2+\a_4)(h-1),
\frac{1}2 (\a_4+\a_6)(h-1) \right \}.
\]

Since $M \cong EE_8$, $\det( M)=2^8$ and $|M/ (F_M+M\cap J)|=2^2$.
Note that $M$, $F_M$ and $M\cap J$ are all doubly even.
Recall that
$D_6^*/D_6\cong \ZZ_2\times \ZZ_2$. Since $M\cap J$ is a direct
summand of $M$,
the natural map $\drt{M}\to \dg{\drtp{M\cap J}}$ is onto.
Similarly, the natural map $\drt{M}\to \dg{\drtp{F_M}}$ is also
onto.

\medskip

Define $H:=F_M\cap F_N$. Let $H_X:=ann_{F_X}(H)$, for  $X=M, N$.
Since $H$ is the negated sublattice of the involution $t_N$ on
$F_M$, $H$ is isometric to
either $0, AA_1, AA_1\perp AA_1$ or $2A_1$ since $F_M$ and $F_N$ are rectangular.

\medskip

Let $\{\a_M^1, \a_M^2\}$ and $\{\a_N^1, \a_N^2\}$ be bases of $F_M$
and $F_N$ such that $(\a_M^i,\a_M^j)=4\delta_{i,j}$ and
$(\a_N^i,\a_N^j)=4\delta_{i,j}$.  Since  $|M: F_M+M\cap J|=2^2$ and
the natural map $\drt{M} \to \dg{\drt{F_M}}$ is onto, there exist
$\b^1 \in \dual {(M\cap J)}, \b^2\in \dual {(M\cap J)}$ so that
\[
\begin{split}
\xi_M= \frac{1}2 \a_M^1 + \b^1 \quad \text{ and }\quad  \zeta_M=
&\frac{1}2\a_M^2 + \b^2
\end{split}
\]
are glue vectors and  the cosets  $\drtp{ \b^1 +(M\cap J)}, \drtp{\b^2+
(M\cap J)}$ generate the abelian group $\drtpp{M\cap J}^*/ \drtp{M\cap
J}\cong \ZZ_2 \oplus \ZZ_2$. Since $M$ is spanned by
norm 4 vectors,
we may also assume that $\xi_M$ and $\zeta_M$ both have norm $4$ and
thus $\b_1$ and $\b_2$ have norm 3.

Recall that a standard basis for the root lattice $D_6$ is given by
 $\{ (1,1,0,0,0,0)$, $(-1,1,0,0,0,0)$, $(0, -1,1,0,0,0)$,
$(0,0,-1,1,0,0)$, $(0,0,0,-1,1,0)$, $(0,0,0,0,-1,1)\}$ and the
elements of norm $3/2$ in $\dual{(D_6)}$ have the form $\half (\pm 1,
\dots, \pm 1)$ with evenly many $-$ signs or $\half (\pm 1, \dots,
\pm 1)$ with oddly many $-$ signs (cf. \cite[Chapter 5]{cs}). They
are contained in two distinct cosets of $\dual{(D_6)}/D_6$. Note that
$(D_6)^*/D_6$ have 3 nontrivial cosets and their elements have norm
$3/2$, $3/2$, and $1$ modulo $2\ZZ$, respectively.

Now define $\phi: D_6 \to M\cap J$ by
\[
\begin{split}
 (1,1,0,0,0,0)\mapsto \a_1, \quad (-1,1,0,0,0,0)\mapsto \a_2, \quad (0,
-1,1,0,0,0)\mapsto \a_3\ \\
(0,0,-1,1,0,0)\mapsto \a_4, \quad (0,0,0,-1,1,0)\mapsto \a_5,\quad
(0,0,0,0,-1,1)\mapsto \a_6.
\end{split}
\]
A comparison of Gram matrices shows that  $\phi$ is $\sqrt{2}$ times an isometry. Since $\frac{1}2(-1,1,-1,1,-1,1)$ and
$\frac{1}2(1,1,-1,1,-1,1)$ are the representatives of the two cosets
of $(D_6)^*/D_6$ represented by norm $3/2$ vectors,  by
\refpp{dualrescaling},
\[
\phi\left( \frac{1}2(-1,1,-1,1,-1,1)\right)=
\frac{1}2(\a_2+\a_4+\a_6)
\]
and
\[
\phi\left( \frac{1}2(1,1,-1,1,-1,1)\right)=
\frac{1}2(\a_1+\a_4+\a_6)
\]
are the representatives of the two cosets of $2(M\cap J)^*/(M\cap
J)$ represented by norm $3$ vectors. Therefore, without loss, we may
assume
\[
\{\b^1, \b^2\}= \{\frac{1}2 (\a_2+\a_4+\a_6), \frac{1}2 (\a_1+\a_4+\a_6)\}.
\]

Similarly, there exist $\gamma^1, \gamma^2\in (N\cap J)^*$ with
$(\gamma^1, \gamma^1)=(\gamma^2, \gamma^2)=3$ such that
\[
\xi_N= \frac{1}2 \a_N^1 + \gamma^1,\quad \text{ and }\quad
 \zeta_N= \frac{1}2\a_N^2 + \gamma^2
\]
are glue vectors and $N=\mathrm{span}_\ZZ\{ F_N+N\cap J, \xi_N,
\zeta_N\}$. Moreover, $\drtp{ \gamma^1 +N\cap J}$,  $\drtp{\gamma^2+
N\cap J}$ generate the group $\drtpp{N\cap J}^*/ \drtp{N\cap J}$.

We shall prove that $(\b, \g)\equiv \half (mod \ \ZZ)$, resulting in a contradiction.

Define $\varphi: D_6 \to N\cap J$ by
\begin{align*}
 (1,1,0,0,0,0)&\mapsto \eta_1, &&(-1,1,0,0,0,0)\mapsto \eta_2, \\
 (0, -1,1,0,0,0)&\mapsto
 \frac{1}2(\eta_1+\eta_2+\eta_3+\eta_4), &&
(0,0,-1,1,0,0)\mapsto \eta_4,\\
(0,0,0,-1,1,0)&\mapsto \frac{1}2(-\eta_3+\eta_4-\eta_5+\eta_6), &&
(0,0,0,0,-1,1)\mapsto \eta_6.
\end{align*}

By comparing the Gram matrices, it is easy to show that $\varphi$ is
a $\sqrt{2}$ times an isometry. Thus, we may choose $\gamma^1,
\gamma^2$ such that
\[
\begin{split}
\{\gamma^1, \gamma^2\}&=\left\{\varphi\left(
\frac{1}2(-1,1,-1,1,-1,1)\right),  \varphi\left(
\frac{1}2(1,1,-1,1,-1,1)\right)\right\}\\
&=\left\{\half( \eta_2+\eta_4+\eta_6), \half(
\eta_1+\eta_4+\eta_6)\right\}.
\end{split}
\]
By the definition of $\a_1, \dots, \a_6$, we have
\[
\begin{split}
\eta_1&= \half(\a_1-\a_2)(h-1),\quad  \eta_2= \half(\a_1+\a_2)(h-1), \\
\eta_4&=[\half(\a_1+\a_2)+(\a_3+\a_4)](h-1),\\
\eta_6&=[\half(\a_1+\a_2)+(\a_3+\a_4+\a_5+\a_6)](h-1).
\end{split}
\]
Thus,
\[
(\a_2+\a_4+\a_6, \eta_2+\eta_4+\eta_6)= -6 \quad \text{ and }\quad
(\a_2+\a_4+\a_6, \eta_1+\eta_4+\eta_6)= -2.
\]
Therefore, $(\b^1, \gamma^1)\equiv 1/2 \mod \ZZ$.

\medskip

{\sl Subcase 1.} $H\cong 0, AA_1$ or $AA_1\perp AA_1$.  In this case,
we may choose the bases $\{\a_M^1, \a_M^2\}$ and $\{\a_N^1,
\a_N^2\}$ of  $F_M$ and $F_N$ such that $(\a_M^i, \a_N^j) \in \{0,
4\}$, for all $i, j$. Hence,
\[
(\xi_M, \xi_N)= (\frac{1}2\a_M, \frac{1}2 \a_N)+ (\b^1, \gamma^1) \equiv 1/2 \mod \ZZ
\]
and  $L=M+N$ is non-integral.

\medskip

{\sl Subcase 2.} $H\cong 2A_1$. Then $H_M\cong H_N \cong 2A_1$, also.
By replacing $\a_M^i$ by $-\a_M^i$ and $\a_N^i$ by $-\a_N^i$ for
$i=1,2$ if necessary, $\a_M^1+\a_M^2=\a_N^1+\a_N^2 \in H$.  Write
$\rho := \a_M^1+\a_M^2=\a_N^1+\a_N^2$.  Then  we calculate the
difference of the glue vectors
\[
\begin{split}
\eta_M-\zeta_M&= \frac{1}2(\a_M^1-\a_M^2)+ \half(\a_2+\a_4+\a_6)-\half(\a_1+\a_4+\a_6)\\
&\equiv \frac{1}2\rho + \half(-\a_1+\a_2) \mod (F_M+M\cap J).
\end{split}
\]
Similarly,
\[
\begin{split}
\eta_N-\zeta_N& = \frac{1}2(\a_N^1-\a_N^2)+ \half(\eta_2+\eta_4+\a_6)-\half(\eta_1+\eta_4+\eta_6)\\
&\equiv \frac{1}2\rho - \half(-\eta_1+\eta_2) \mod (F_N + N\cap
J).
\end{split}
\]
Let $\nu_M=\frac{1}2\rho + \half(-\a_1+\a_2)$ and
$\nu_N=\frac{1}2\rho - \half(-\eta_1+\eta_2)$. Then $\nu_M$ and
$\nu_N$ are both norm $4$ vectors in $L$. Recall that
$(-\eta_1+\eta_2)=\a_2(h-1)$.  Since $(\a_i, \a_M^j)=(\a_i, \a_N^j)=0$ for all $i, j$, we have $(\rho, \a_i)=0$ for all $i$.
\[
\begin{split}
(\nu_M, \nu_N)=&(\frac{1}2\rho + \half(-\a_1+\a_2),
\frac{1}2\rho - \half(-\eta_1+\eta_2))\\
= &\frac{1}4\big[ (\rho, \rho) - (-\a_1+\a_2,
\a_2(h-1))\big]\\
\end{split}
\]
Recall that $(\a_M^i, \a_M^j)=4 \delta_{i,j}$ and $(\a_i, \a_j)=4
\delta_{i,j}$ for $i,j=1,2$. Moreover, $(x,y h)=0$ for all $x,y \in M\cap
J$ by (ii) of \refpp{squareofg-1}. Thus, $(\rho,\rho)=(\a_M^1+\a_M^2,
\a_M^1+\a_M^2)=4+4=8$  and $(-\a_1+\a_2,
\a_2(h-1))= (-\a_1+\a_2, -\a_2)=-4$.

Therefore,
$
(\nu_M, \nu_N)= \frac{1}4(8-(-4))=3$ and hence
$\nu_M-\nu_N$ is a root in $L$. \eop

\medskip

\subsubsection{$\Dih{8}$: Given that $F_M\cong F_N\cong DD_4$, what is $J$?}

Next we shall consider the case $F_M\cong F_N\cong DD_4$.
In this case, $M\cap J\cong N\cap J\cong DD_4$.

\begin{nota}\labtt{aa1indd4}
Let $\{\a_1,\a_2,\a_3,\a_4\}\subset M\cap J$ such that $(\a_i,\a_j)=
4\delta_{i,j}, i,j=1,2,3,4$.   Then $ M\cap J= \mspan_\ZZ\{\a_1,
\a_2, \a_3, \frac{1}2(\a_1+\a_2+\a_3+\a_4)\}. $  In this case, the
norm $8$ elements of $M\cap J$  are given by $\pm \a_i\pm \a_j$ for
$i\neq j$.
\end{nota}

\begin{lem}\labtt{dih8:ncapj}
Let $h=g^{-1}$.  By rearranging the subscripts if necessary,
we have
\[
N\cap J =\mspan_\ZZ\{ (M\cap J)(h-1), \half(\a_1+\a_2)(h-1),
\half(\a_1+\a_3)(h-1)\}.
\]
\end{lem}
\pf Let $K: =(M\cap J)(h-1)$.  Then by (ii) of Lemma \ref{conjcalc},
we have $K \leq N\cap J$.  Since $( M\cap J, (M\cap J)h ) =0$ by
\refpp{squareofg-1}, $K= (M\cap J)(h-1)\cong 2D_4$. Therefore, by
\refpp{dd4ind4}, we have  $K\leq N\cap J \leq \half K$.

Note that, by determinants, $|N\cap J: K|=\sqrt{2^4}=2^2$. Therefore,
there exists two glue vectors $\b_1,\b_2\in (N\cap J)\setminus K$
such
that $ N\cap J =\mspan_\ZZ\{ (M\cap J)(h-1), \b_1, \b_2\}.$

Since $K$ has minimal norm $8$ and $N\cap J$ is generated by
norm $4$ elements, we may choose $\b_1,\b_2$ such that both are of
norm $4$. On the other hand, elements of norm $4$ in $N\cap J$ are
given by $\frac{1}2 \g(h-1)$, where $\g\in M\cap J$ with norm $8$,
i.e., $\g=\pm \a_i\pm \a_j$ for some $i\neq j$. Since
$\a_1(h-1),\a_2(h-1),\a_3(h-1),\a_4(h-1)\in (M\cap J)(h-1) \leq
N\cap J$, we may assume
\[
\b_1=\half(\a_i+\a_j)(h-1)\quad  \text{ and } \quad
\b_2=\half(\a_k+\a_\ell)(h-1)
\]
for some $i,j,k,\ell\in \{1,2,3,4\}$. Note that $|\{i,j\}\cap
\{k,\ell\}|=1$ because $\b_1+\b_2\notin K$. Therefore, by rearranging
the indices if necessary,  we may assume
$\b_1=\half(\a_1+\a_2)(h-1), \b_2=\half(\a_1+\a_3)(h-1)$ and
\[
N\cap J =\mspan_\ZZ\{ (M\cap J)(h-1), \half(\a_1+\a_2)(h-1),
\half(\a_1+\a_3)(h-1)\}.
\]
as desired. \eop

\begin{prop}
If $F_M\cong F_N\cong DD_4$, then $M\cap J+ N\cap J \cong EE_8$.
\end{prop}

\pf First we shall note that $ (M\cap J)+ (M\cap J)(h-1)=(M\cap J)
\perp (M\cap J)h \cong DD_4\perp DD_4. $ Moreover, we have $|N\cap J+
M\cap J: (M\cap J)+ (M\cap J)(h-1)|=|N\cap J: (M\cap
J)(h-1)|=\sqrt{(2^8\cdot 4)/(2^4\cdot 4)}= 4$, by determinants.
Therefore, $\det( M\cap J +N\cap J)= (2^4 \cdot 4)^2/ 4^2=2^8$.

Now by \refpp{dih8:ncapj}, we have
\[
N\cap J =\mspan_\ZZ\{ (M\cap J)(h-1), \half(\a_1+\a_2)(h-1),
\half(\a_1+\a_3)(h-1)\}.
\]
 Next we shall show that $(M\cap J, N\cap J)\subset 2\ZZ$.
Since $(M\cap J, (M\cap J)h))=0$ and $M\cap J$ is doubly even, it is
clear that $(M\cap J, (M\cap J)(h-1))\subset 2\ZZ$. Moreover, for
any $ i,j\in \{1,2,3,4\}, i\neq j$,
\[
\left(\a_k, \frac{1}2(\a_i+\a_j)(h-1)\right)=
\begin{cases}
0 &\text{ if } k\notin \{i,j\},\\
-2&\text{ if } k\in \{i,j\},
\end{cases}
\]
and
\[
\left (\frac{1}2(\a_1+\a_2+\a_3+\a_4), \frac{1}2(\a_i+\a_j)(h-1)\right)=-2.
\]
Since $M\cap J$ is spanned by $\a_1, \a_2, \a_3$ and
$\frac{1}2(\a_1+\a_2+\a_3+\a_4)$ and $N\cap J =\mspan_\ZZ\{ (M\cap
J)(h-1), \half(\a_1+\a_2)(h-1), \half(\a_1+\a_3)(h-1)\}$, we have
$(M\cap J, N\cap J)\subset 2\ZZ$ as required. Therefore,
$\drtp{M\cap J+ N\cap J}$ is an integral lattice and has determinant
$1$ and thus $M\cap J+ N\cap J\cong EE_8$, by the classification of
unimodular even lattices of rank $8$. \eop

\medskip

\begin{lem}\labtt{dualmjnj}
Let $\a_1, \a_2, \a_3,\a_4\in M\cap J$ be as in Notation
\ref{aa1indd4}. Then
\[
(M\cap J)^*= \frac{1}4 \mspan_\ZZ\left \{ \a_1- \a_2, \a_1+\a_2,
\a_1+\a_3, \a_1+\a_4\right\}.
\]
and
\[
\begin{split}
\dual{(N\cap J)}&= \frac{1}4\mspan_\ZZ \left \{ \a_1(h-1),
\a_2(h-1),\a_3(h-1), \frac{1}2(\a_1+\a_2+\a_3+\a_4)(h-1)\right\}.
\end{split}
\]
\end{lem}
\pf
Since $\dual{(\ZZ\a_i)}=\frac{1}4\ZZ\a_i$ and  $M\cap J=\mspan_\ZZ\{ \a_1,\a_2,\a_3, \frac{1}2(\a_1+\a_2+\a_3+\a_4)\}$,
\[
\begin{split}
&(M\cap J)^*\\
= &\left \{\left.\b=\frac{1}4(a_1\a_1+a_2\a_2+a_3\a_3+a_4\a_4)
\right |a_i\in \QQ, (\b,\gamma)\in
\ZZ\text{ for all }\gamma\in M\cap J\right \}\\
= &\left \{ \frac{1}{4}(a_1\a_1+ a_2\a_2+a_3\a_3+a_4\a_4)\left |\,
a_i\in \ZZ \text{ and } \sum_{i=1}^4 a_i\in 2\ZZ, i=1,2,3,4\right.
\right\}\\
 =&\frac{1}4 \mspan_\ZZ\left \{ \a_1- \a_2, \a_1+\a_2,
\a_1+\a_3, \a_3+\a_4\right\}.
\end{split}
\]

Now by \refpp{dih8:ncapj}, we have
\[
N\cap J= \mspan_\ZZ\left \{ (M\cap J)(h-1),
\frac{1}2(\a_1+\a_2)(h-1),\frac{1}2(\a_1+\a_3)(h-1)\right \}.
\]

Let $ \b_1=\frac{1}2(\a_1+\a_2)(h-1)$,
$\b_2=\frac{1}2(\a_1-\a_2)(h-1)$, $\b_3=\frac{1}2(\a_3+\a_4)(h-1)$,
and $\b_4=\frac{1}2(\a_3-\a_4)(h-1)$. Then $\{\b_1,
\b_2,\b_3,\b_4\}$  forms an orthogonal subset of $N\cap J$ with
$(\b_i,\b_j)=4\delta_{i,j}$. Note that
$\frac{1}2(\b_1+\b_2+\b_3+\b_4)=\frac{1}2(\a_1+\a_3)(h-1)\in N\cap
J$. Thus, $ N\cap J= \mspan_\ZZ\{ \b_1, \b_2,\b_3,
\frac{1}2(\b_1+\b_2+\b_3+\b_4)\} $  since both of them are isomorphic to $DD_4$. Hence we have
\[
\begin{split}
(N\cap J)^* & =\frac{1}4 \mspan_\ZZ\left \{ \b_1- \b_2, \b_1+\b_2,
\b_1+\b_3, \b_3+\b_4\right\}\\
&= \frac{1}4\mspan_\ZZ \left \{ \a_1(h-1), \a_2(h-1),\a_3(h-1),
\frac{1}2(\a_1+\a_2+\a_3+\a_4)(h-1)\right\}.
\end{split}
\]
as desired. \eop

\begin{lem}\labtt{cosetindd4}
We shall use the same notation  as in \refpp{aa1indd4}. Then the
cosets of $2\dual{(M\cap J)}/(M\cap J)$ are represented by
\[
0,\ \frac{1}{2}(\a_1+\a_2),\ \frac{1}{2}(\a_1+\a_3),\
\frac{1}{2}(\a_2+\a_3),
\]
and the cosets of $2\dual{(N\cap J)}/(N\cap J)$ are represented by
\[
0,\ \frac{1}{2}\a_1(h-1),\ \frac{1}{4}(\a_1+\a_2+\a_3+\a_4)(h-1),\
\frac{1}{4}(\a_1+\a_2+\a_3-\a_4)(h-1).
\]
Moreover, $(2(M\cap J)^*,2 \dual{(N\cap J)})\subset \ZZ $.
\end{lem}
\pf Since $X:=M\cap J\cong DD_4$, it is clear that $2\dual X /X$ is a
four-group.  The three nonzero vectors in the list
\[
0,\ \frac{1}{2}(\a_1+\a_2),\ \frac{1}{2}(\a_1+\a_3),\
\frac{1}{2}(\a_2+\a_3),
\]
have norms two, so all are in $2\dual X \setminus X$.  Since the difference of any two has norm 2, no two are congruent modulo $X$.   A similar argument proves the second statement.

For the third statement, we calculate the following inner products.

For any $i,j, k\in \{1,2,3,4\}$ with $i\neq j$,
\[
(\a_i\pm \a_j, \a_k(h-1))=
\begin{cases}
0 & \text{ if } k\notin \{i,j\},\\
\pm 4 &\text{ if } k\in \{i,j\},
\end{cases}
\]
and
\[
(\a_i\pm \a_j, \frac{1}2(\a_1+\a_2+\a_3+\a_4)(h-1)) = 0 \text{ or }
-4.
\]
Since $(M\cap J)^*= \frac{1}4 \mspan_\ZZ\left \{ \a_1- \a_2,
\a_1+\a_2, \a_1+\a_3, \a_1+\a_4\right\} $ and $ \dual{(N\cap J)}=
\frac{1}4\mspan_\ZZ \left \{ \a_1(h-1), \a_2(h-1),\a_3(h-1),
\frac{1}2(\a_1+\a_2+\a_3+\a_4)(h-1)\right\}$ by \refpp{dualmjnj}, we
have $((M\cap J)^*, \dual{(N\cap J)})\subset \frac{1}4 \ZZ $ and
hence $(2(M\cap J)^*,2 \dual{(N\cap J)})\subset \ZZ $ as desired.
\eop

\begin{rem}\labttr{sevenconjclasses}
Note that the lattice $D_4$ is $\bw 2$, so the involutions in its
isometry group $\brw 2\cong Weyl(F_4)$  may be deduced from the theory in
\cite{ibw1}, especially Lemma 9.14 (with $d=2$).  The results are in Table \refpp{involsond4}.
\end{rem}

\begin{nota}\labttr{defofh}
Define $H:=F_M\cap F_N$ and let $H_X:=ann_{F_X}(H)$ for $X=M, N$.
Since $H$ is the negated sublattice of the involution $t_N$ on
$F_M$, we have the possibilities listed in Table \ref{involsond4}.
We label the case for $\drt{H}$ by the corresponding involution $2A,
\cdots , 2G$. (i.e., the involution whose negated space is
$\drt{H}$)
\end{nota}

\begin{table}[bht]
\caption{ \bf  \bf The seven conjugacy classes of involutions in
$\brw 2\cong Weyl(F_4)$}
\begin{center}
\begin{tabular}{|c|c|c|}
\hline Involution & Multiplicity of $-1$ & Isometry type of  \cr
   &  & negated sublattice \cr
\hline
   2A&1&$A_1$\cr \hline
   2B&1&$AA_1$\cr \hline
   2C&2&$A_1\perp A_1$\cr \hline
   2D&2&$A_1\perp AA_1$\cr
 \hline 2E&3&$A_1\perp A_1\perp A_1$\cr \hline
   2F&3&$A_3$\cr \hline 2G&4&$D_4$\cr \hline
\end{tabular}
\end{center}
\label{involsond4}
\end{table}%

We shall prove the main result of this section, Theorem
\ref{dih8:rootlesscases} in several steps.

\begin{lem}\labtt{dih8lem1}
Suppose $F_M\cong F_N\cong DD_4$. If $\drt{H}\cong AA_1, A_1\perp AA_1$
or $A_3$ (i.e., the cases for $2B$, $2D$ and $2F$), then the lattice
$L$ is non-integral.
\end{lem}

\pf  We shall divide the proof into 3 cases.  Recall notations \refpp{defofh}.

\textsl{Case 2B.} In this case, $\drt{H}\cong AA_1$ and
$\drt{H_M}\cong \drt{H_N}\cong A_3$. Then $F_M\geq H\perp H_M$ and $
M\geq H\perp H_M\perp M\cap J$.

Let $\alpha \in H$ with $( \alpha , \alpha)=8$. Then $H=\ZZ\a$ and
$H^*=\frac{1}8 \ZZ\a$.

By \refpp{dualrescaling},  we have
$\dual{(\drt{F_M})}/\drt{F_M}\cong 2(F_M)^*/F_M$. Thus, by
\refpp{aa1ind4}, the natural map $2(F_M)^* \to
2(H^*)=\frac{1}4\ZZ\a$ is onto. Therefore, there exists $\delta_M\in
H_M^*$ with $(\delta_M, \delta_M)= 3/2$ such that $\frac{1}4
\a+\delta_M\in 2(F_M)^*$. Note that the natural map $\drt{M}\to
\dg{\drt{F_M}}$ is also onto since $\drt{M}$ is unimodular and $F_M$ is a
direct summand of $M$. Therefore, there exists $\gamma_M\in
2\dual{(M\cap J)}$ with $(\gamma_M,\gamma_M)=2$ such that
\[
\xi_M= \frac{1}4\a +\delta_M+\gamma_M
\]
is a glue vector for $H\perp H_M\perp M\cap J$ in $M$. Similarly,
there exists $\delta_N\in H_N^*$ with $(\delta_N, \delta_N)= 3/2$
and $\gamma_N\in 2\dual{(N\cap J)}$ with $(\gamma_N,\gamma_N)=2$
such that
\[
\xi_N=\frac{1}4\a +\delta_N+\gamma_N
\]
is a glue vector for $H\perp N_N \perp N\cap J$ in $N$.

Since $(\gamma_M, \gamma_N)\in \ZZ$ by \refpp{cosetindd4} and $H_M\perp H_N$,
\[
(\xi_M, \xi_N) = \frac{1}{16}(\alpha, \alpha) + (\delta_M,
\delta_N)+(\gamma_M, \gamma_N) \equiv \half \mod \ZZ,
\]
which is not an integer.

\medskip

\textsl{Case 2D.} In this case, $\drt{H}\cong A_1\perp AA_1$ and
$\drt{H_M}\cong \drt{H_N}\cong A_1\perp AA_1$, also.  Take $\alpha,
\beta\in H$ such that $(\alpha ,\alpha) =8$, $(\beta ,\beta ) =4$
and $(\alpha ,\beta) =0$. Similarly, there exist $\alpha_M, \beta_M\in
H_M$ and $\alpha_N, \beta_N\in H_N$ such that $( \alpha_M ,\alpha_M
) =8$, $(\beta_M ,\beta_M)=4$ and $(\alpha_M ,\beta_M ) =0$ and
$(\alpha_N ,\alpha_N) =8$, $( \beta_N ,\beta_N) =4$ and $(\alpha_N
,\beta_N ) =0$.

Note that $\ZZ\b\perp \ZZ\b_M\perp \ZZ\a_M\leq ann_{F_M}(\a)\cong
AA_3$. Set $A:= ann_{F_M}(\a)\cong AA_3$. Then $|A: \ZZ\b\perp
\ZZ\b_M\perp \ZZ\a_M|= \sqrt{ (4\times 4\times 8)/ (2^3\times
4)}=2.$ Therefore, there exists a $\mu \in (\ZZ\b\perp \ZZ\b_M\perp
\ZZ\a_M)^*= \frac{1}4 \ZZ\b\perp \frac{1}4 \ZZ\b_M\perp
\frac{1}8\ZZ\a_M$ such that $A=\mspan_\ZZ\{ \b, \b_M, \a_M, \mu\} $
and $2\mu \in \ZZ\b\perp \ZZ\b_M\perp \ZZ\a_M$.

Since $A$ is generated by norm 4 vectors, we may choose $\mu$ so
that $\mu$ has norm 4. The only possibility is $\mu =\frac{1}2 ( \pm
\b\pm \b_M\pm \a_M)$. Therefore, $A=\mspan_\ZZ\{ \b, \b_M, \a_M,
\frac{1}2(\b+\b_M+\a_M)\}$ and $2A^*/A\cong \ZZ_4$ is generated by
$\frac{1}2 \b +\frac{1}{4}\a_M+ A$. Note also that $\frac{1}2 \b
+\frac{1}{4}\a_M$ has norm $3/2$.

Now recall that $\ZZ\a\cong 2A_1$ and $A\cong AA_3$.
Thus, by
\refpp{aa1ind4}, the natural map $2(F_M)^* \to
2(H^*)=\frac{1}4\ZZ\a$ is onto. Thus, there exists a $\delta\in
2A^*$ with $(\delta, \delta)=3/2$ such that $\frac{1}4\a +\delta\in
2(F_M)^*$.  By the previous paragraph, we may assume $\delta=
\frac{1}2 \b +\frac{1}{4}\a_M$. Since the natural map $\drt{M}\to
\dg{\drt{F_M}}$ is onto, there exists $\gamma_M\in 2(M\cap J)^*$
such that
\[
\xi_M= \frac{1}4 \alpha + \frac{1}2 \beta + \frac{1}4 \alpha_M
+\gamma_M
\]
is a glue vector for $H\perp H_M\perp M\cap J$ in  $M$. Similarly,
there exists $\gamma_N\in 2(N\cap J)^*$ such that
\[
\xi_N= \frac{1}4 \alpha + \frac{1}2 \beta + \frac{1}4 \alpha_N
+\gamma_N
\]
is a glue vector for $H\perp H_N\perp N\cap J$ in $N$. Then
\[
(\xi_M, \xi_N)= \frac{1}{16}(\alpha,\alpha) + \frac{1}4( \beta,
\beta) + \frac{1}{16} (\alpha_M, \alpha_N) + (\gamma_M, \gamma_N
)\equiv 1/2 \mod \ZZ,
\]
since $(\a_M,\a_N)=0$ and $(\gamma_M, \gamma_N
)\in \ZZ$ by \refpp{cosetindd4}.   Therefore, $L$ is not integral.

\medskip

\textsl{Case 2F.} In this case, $\drt{H}\cong A_3$ and
$\drt{H_M}\cong \drt{H_N} \cong AA_1$. Then $F_M\geq H\perp H_M$ and
$F_N\geq H\perp H_N$. Let $\delta\in H^*$, $\alpha_M \in H_M$ and
$\alpha_N \in H_N$ such that $(\delta, \delta)= 3/2$,  $(\alpha_M,
\a_M)=8$ and $(\alpha_N, \a_N)=8$.

Recall in Case 2B that $\drt{H}\cong AA_1$ and $\drt{H_M}\cong \drt{H_N}\cong A_3$. Now by exchanging the role of $H$ with $H_M$ (or $H_N$) and using the same argument as in Case 2B, we may show that there exist $\gamma_M\in
2(M\cap J)^*$ and $\gamma_N\in 2(N\cap J)^*$ such that
$\xi_M=\delta+\frac{1}4\a_M+\gamma_M$ is a glue vector for $H\perp
H_M\perp M\cap J$ in  $M$ and $\xi_N=\delta+\frac{1}4\a_N+\gamma_N$
is a glue vector for $H\perp H_N\perp (N\cap J)$ in  $N$. However,
\[
(\xi_M,\xi_N)=(\delta, \delta)+(\gamma_M, \gamma_N) \equiv 1/2 \mod
\ZZ,
\]
since $(\gamma_M, \gamma_N
)\in \ZZ$ by \refpp{cosetindd4}. Again, $L$ is not integral. \eop
\medskip

\begin{lem}\labtt{innerproduct=-1}
Let $\gamma_M$ be any norm $2$ vector in $2(M\cap J)^*$. Then for each non-zero coset
$\gamma_N + (N\cap J)$ in $2(N\cap J)^*/ (N\cap J)$, there exists a norm $2$ vector $\gamma\in \gamma_N + (N\cap J)$ such that $ (\gamma_M, \gamma)=-1$.
\end{lem}

\pf Recall from \refpp{dualmjnj} that
\[
2(M\cap J)^* =\frac{1}2 \mspan_\ZZ\left \{ \a_1- \a_2, \a_1+\a_2,
\a_1+\a_3, \a_3+\a_4\right\}.
\]
Thus, all norm $2$ vectors in $2(M\cap J)^*$ have the form $\frac{1}2(\pm \a_i\pm \a_j)$ for some $i\neq j$.
Without loss, we may assume $\gamma_M= \frac{1}2( \a_i+\a_j)$ by replacing $\a_i$, $\a_j$ by $-\a_i$, $-\a_j$ if necessary.

Now by \refpp{cosetindd4}, the non-zero cosets of $2(N\cap J)^*/ (N\cap J)$ are represented by
$\frac{1}2 \a_1(h-1)$, $\frac{1}4(\a_1+\a_2+\a_3+\a_4)(h-1)$ and $\frac{1}4(
\a_1+\a_2+\a_3-\a_4)(h-1)$. Moreover, by \refpp{dih8:ncapj},
\[
N\cap J= \frac{1}2 \mspan_\ZZ\{\, (\a_i\pm \a_j)(h-1)\, |\ 1\leq i< j\leq 4\}.
\]

If $\gamma_N + (N\cap J)=\frac{1}2 \a_1(h-1) + (N\cap J)$,  we take
\[
\gamma= \frac{1}2 \a_i(h-1) = \frac{1}2 \a_1(h-1) + \frac{1}2(-\a_1+ \a_i)(h-1) \in \frac{1}2 \a_1(h-1) + (N\cap J).
\]

Recall from \refpp{aa1indd4} that $\{\a_1, \a_2,\a_3,\a_4\}\in M\cap J$
and
$(\a_i, \a_j)=4\delta_{i,j}$ for $i,j=1, \dots, 4$. Moreover, $(x,yh)=0$ for all $x,y \in M \cap J$ by
\refpp{squareofg-1}.

\medskip

Thus, $(\gamma, \gamma)=(\frac{1}2 \a_i(h-1), \frac{1}2 \a_i(h-1))=\frac{1}4 [(\a_i h, \a_i h)+ (\a_i, \a_i)]=2 $ and $
(\gamma_M^1, \gamma)= (\frac{1}2(\a_i+\a_j),
\frac{1}2\a_i(h-1))=-\frac{1}{4}(\a_i+\a_j, \a_i)=-1.
$

\medskip

If $\gamma_N + (N\cap J)=\frac{1}4(\a_1+\a_2+\a_3+\a_4)(h-1) + (N\cap J)$, we simply take $\gamma= \frac{1}4(\a_1+\a_2+\a_3+\a_4)(h-1)$. Then $(\gamma, \gamma)=2$ and
\[
\begin{split}
(\gamma_M,\gamma)= & \
(\frac{1}2(\a_i+\a_j),
\frac{1}4(\a_1+\a_2+\a_3+\a_4)(h-1))\\
= &
-\frac{1}8(\a_i+\a_j,\a_1+\a_2+\a_3+\a_4)\\
=&-\frac{1}8(4+4)=-1.
\end{split}
\]

Finally we consider the case $\gamma_N + N\cap J=\frac{1}4(\a_1+\a_2+\a_3-\a_4)(h-1) + (N\cap J)$.
Let $\{k,\ell\}=\{1,2,3,4\}-\{i,j\}$ and take
\[
\begin{split}
\gamma & = \frac{1}4(\a_i+\a_j+\a_k-\a_\ell)(h-1)= \frac{1}4(\a_1+\a_2+\a_3-\a_4)(h-1)+ \frac{1}2 (\a_4 - \a_\ell)\\
&\in \frac{1}4(\a_1+\a_2+\a_3-\a_4)(h-1)+ N\cap J.
\end{split}
\]
Then $(\gamma, \gamma)=2$ and
\[
\begin{split}
(\gamma_M,\gamma)= & \
(\frac{1}2(\a_i+\a_j),
\frac{1}4(\a_i+\a_j+\a_k-\a_\ell)(h-1))\\
= &
-\frac{1}8(\a_i+\a_j, \a_i+\a_j+\a_k-\a_\ell)\\
=&-\frac{1}8(4+4)=-1
\end{split}
\]
as desired.
\eop

\begin{lem}\labtt{dih8lem2}
If $\drt{H}\cong  A_1\perp A_1,  A_1\perp A_1\perp A_1$ or $D_4$ (i.e.,
the cases for $2C$, $2E$ and $2G$), then the lattice $L$ has roots.
\end{lem}

\pf  We continue to use the notations \refpp{defofh}.
First, we shall note that the natural maps $\drt{M}\to
\dg{\drt{F_M}}$,  $\drt{M}\to \dg{\drtp{M\cap J}}$ and $\drt{N}\to
\dg{\drt{F_N}}$, and $\drt{N}\to \dg{\drtp{N\cap J}}$ are all onto.

\textsl{Case 2C.} In this case, $\drt{H}\cong A_1\perp A_1$ and
$\drt{H_M}\cong \drt{H_N}\cong A_1\perp A_1$. Let $\mu^1, \mu^2$ be a
basis of $H$ such that $(\mu^i,\mu^j)=4\delta_{i,j}$. Let
$\mu_{M}^{1}, \mu_{M}^{2}$ and $\mu_{N}^{1}, \mu_{N}^{2}$ 
be
bases of $H_M$ and $H_N$ which consist of norm 4 vectors. Then,
$F_M=\mspan_\ZZ\{ \mu^1, \mu^2, \mu^1_M,
\frac{1}2(\mu^1+\mu^2+\mu^1_M+\mu^2_M)\}$ and $F_N=\mspan_\ZZ\{
\mu^1, \mu^2, \mu^1_N, \frac{1}2(\mu^1+\mu^2+\mu^1_N+\mu^2_N)\}$.
Therefore, by the same arguments as in Lemma \ref{cosetindd4}, the
cosets representatives of $(2F_M^*)/F_M$ are given by
\[
0,\quad  \frac{1}2(\mu^1+\mu^2),\quad \frac{1}2(\mu^1+\mu^1_M),\quad
\frac{1}2(\mu^2+\mu^1_M),
\]
and the cosets representatives of $(2F^*_N)/F_N$ are given by
\[
0,\quad  \frac{1}2(\mu^1+\mu^2),\quad \frac{1}2(\mu^1+\mu^1_N),\quad
\frac{1}2(\mu^2+\mu^1_N).
\]
Therefore, there exist $\g_M^1, \g_M^2 \in \dual {(M\cap J)}$ so that
\[
\xi_M= \frac{1}2(\mu_1+\mu_2)+\gamma^1_M\quad \text{ and }\quad
\zeta_M = \frac{1}2(\mu_1+\mu_{M}^{1})+\gamma^2_M,
\]
are glue vectors for $F_M+ (J\cap M)$ in $M$ and such that
$\gamma^1_M+
(M\cap J)$, $\gamma^2_M+ (M\cap J)$ generate $2(M\cap J)^*/ (M\cap
J)$.

Similarly, there exist $\g_N^1, \g_N^2 \in \dual {(N\cap J)}$ so that
\[
\xi_N= -\frac{1}2(\mu_1+\mu_2)+\gamma^1_N\quad \text{ and }\quad
\zeta_N = -\frac{1}2(\mu_1+\mu_{N}^{1})+\gamma^2_N,
\]
are glue vectors for $F_N+ N\cap J$ in $N$, and such that $\gamma^1_N+
(N\cap J)$, $\gamma^2_N+ (N\cap J)$ generate $2(N\cap J)^*/ (N\cap
J)$.

By Lemma \refpp{innerproduct=-1}, we may assume $(\gamma^1_N,
\gamma^1_M)=-1$. Then
\[
\begin{split}
(\xi_M,\xi_N)&= (\frac{1}2(\mu_1+\mu_2)+\gamma_M^1,
-\frac{1}2(\mu_1+\mu_2)+\gamma_N^1)\\
&=-\frac{1}4((\mu_1, \mu_1)+(\mu_2,\mu_2))+(\gamma_M^1,
\gamma_N^1)\\
&=-\frac{1}4(4+4)-1=-3
\end{split}
\]
and hence $\xi_M+\xi_N$ is a root.

\medskip

\textsl{Case 2E.} In this case, $\drt{H}\cong A_1\perp A_1\perp A_1$
and $\drt{H_M}\cong \drt{H_N}\cong A_1$.

Let $\mu_1,\mu_2,\mu_3\in H$ be such that $(\mu_i,\mu_j)=4\delta_{i,j}$.  Let
$\mu_M\in H_M$ and $\mu_N\in H_N$ be norm 4 vectors. Then
$H_M=\ZZ\mu_M$ and $H_N=\ZZ\mu_N$. Moreover,
\[
F_M=\mspan_\ZZ \{ \mu_1, \mu_2,\mu_3,
\frac{1}2(\mu_1+\mu_2+\mu_3+\mu_M)\}\cong DD_4
\]
and
\[
F_N=\mspan_\ZZ\{ \mu_1, \mu_2,\mu_3,
\frac{1}2(\mu_1+\mu_2+\mu_3+\mu_N)\}\cong DD_4.
\]
Then, by \refpp{cosetindd4}, $\frac{1}2(\mu_1+\mu_2)$ is in both $2(F_M)^*$ and $2(F_N)^*$. Therefore, there exist $\gamma_M\in 2(M\cap J)^*$ and $\gamma_N\in 2(N\cap J)^*$ such that
\[
\xi_M= \frac{1}{2}(\mu_1+\mu_2)+ \gamma_M \in M \quad
\text{and}\quad \xi_N=- \frac{1}{2}(\mu_1+\mu_2)+ \gamma_N \in N
\]
are  norm 4 glue
vectors for $F_M\perp M\cap J$ in $M$ and $F_N\perp N\cap J$ in $N$, respectively.
By Lemma \refpp{innerproduct=-1}, we may assume $(\gamma_M,
\gamma_N)=-1$.  Then,
\[
\begin{split}
(\xi_M, \xi_N)=&(\frac{1}{2}(\mu_1+\mu_2)+ \gamma_M, -
\frac{1}{2}(\mu_1+\mu_2)+ \gamma_N)\\
= & -\frac{1}4(\mu_1+\mu_2, \mu_1+\mu_2)+(\gamma_M,
\gamma_N)=-2-1=-3
\end{split}
\]
and $\xi_M+\xi_N$ is a root.

\medskip

\textsl{Case 2G.} In this case, $\drt{H} \cong D_4$ and
$\drt{H_M}=\drt{H_N}=0$. Recall that $2(DD_4)^*/DD_4\cong (D_4)^*/D_4$
by \refpp{dualrescaling}  and all non-trivial cosets of $(D_4)^*/D_4$
can be represented by norm $1$ vectors \cite[p. 117]{cs}. Therefore,
the non-trivial cosets of $2(DD_4)^*/DD_4$ can be represented by norm
$2$ vectors.  Thus we can find vectors $\gamma\in 2 H^*$  with
$(\gamma, \gamma)=2$ and $\gamma_M\in 2(M\cap J)^*$, $\gamma_N\in
2(N\cap J)^*$ such that
\[
\xi_M=\gamma+\gamma_M \in M \quad \text{ and }\quad  \xi_N=
-\gamma+\gamma_N \in N
\]
are norm 4 glue vectors for $F_M\perp M\cap J$ in $M$ and $F_N\perp N\cap J$ in $N$, respectively. Again, we may assume
$(\gamma_M, \gamma_N)=-1$ by \refpp{innerproduct=-1} and thus $(\xi_M,\xi_N)=-3$ and there are
roots. \eop

\medskip

\begin{thm}\labtt{dih8:rootlesscases}
Suppose $F_M\cong F_N\cong DD_4$. If $L=M+N$ is integral and
rootless, then $H=F_M\cap F_N=0$ or $\cong AA_1$.
\end{thm}

The proof of Theorem \ref{dih8:rootlesscases} now follows from Lemmas \ref{dih8lem1} and \ref{dih8lem2}.
\section{$\Dih{6}$ and $\Dih{12}$ theories}

We shall study the cases when $D=\la t_M, t_N\ra \cong Dih_6$ or $Dih_{12}$.  The following is our main theorem in this section. We refer to the notation table (Table \ref{notationandterminology}) for the definition of $\dih{6}{14}$, $\dih{6}{16}$ and $\dih{12}{16}$.

\begin{thm}\labtt{thmDih612}
Let $L$ be a rootless integral lattice which is a sum of sublattices $M$ and $N$ isometric to $EE_8$.  If the associated dihedral group has order 6 or 12, the
possibilities for $L+M+N, M, N$ are listed in Table
\refpp{dihedralpairs}.
\end{thm}

\begin{table}[bht]
\caption{ \bf  $\Dih{6}$ and $\Dih{12}$: Rootless cases}
\begin{center}
\begin{tabular}{|c|c|c|c|c|}
\hline  Name &  $F\cong$ & $L$ contains $\dots$&with index \dots & $\dg L$ \cr
\hline
\hline
$\dih{6}{14}$&  $AA_2$ & $\ge A_2\otimes E_6 \perp AA_2$ & $3^2$ &$1^93^36^2$ \cr
\hline
$\dih{6}{16}$ & 0 & $A_2\otimes E_8$  & 1 & $3^8$ \cr
\hline
\hline
$\dih{12}{16}$& $AA_2\perp AA_2$ &  $\ge A_2\otimes E_6 \perp AA_2^{\perp 2}$&$3^2$  & $1^{12}6^4$  \cr
\hline
\end{tabular}
\end{center}
\label{dihedralpairs}
\end{table}%

\subsection{$\Dih{6}$}

\begin{nota}\labttr{dih6setup}
Define $t:=t_M$, $h:=t_Mt_N$.  We suppose $h$ has order $3$. Then,
$N=Mg$, where $g=h^2$.  The third lattice in the orbit of $D:=\la g, t \ra$ is
$Mg^2$, but we shall not refer to it explicitly henceforth. Define
$F:=M\cap N$,  $J:=ann_L(F)$.  Note that $F$ is the common negated
lattice for $t_M$ and $t_N$ in $L$, so is the fixed point sublattice
for $g$ and is a direct summand of $L$ (cf. \refpp{mnfsummands}).
\end{nota}

\begin{lem}\labtt{intertrivial} Let $X=M$ or $N$.
Two of the sublattices $\{(J\cap X)g^i|\,i\in \ZZ \}$ are equal or meet trivially.
\end{lem}
\pf  We may assume $X=M$.
Suppose that $0\ne U=(J\cap M)g^i \cap (J\cap M)g^j$ for $i, j$ not congruent modulo 3.
Then $U$ is negated by two distinct involutions $t^{g^i}$ and $t^{g^j}$, hence is centralized by $g$, a contradiction.
\eop

\begin{lem}\labtt{tensor}
If $F=0$, $J\cong A_2 \otimes E_8$.
\end{lem}
\pf  Use \refpp{tensorwitha2}.
\eop

\begin{hyp} \labttr{s}
We assume $F\neq 0$ and define the integer $s$ by $3^s:=|L/(J+F)|$.
\end{hyp}

\begin{lem} \labtt{lmodjf} $L/(J\perp F)$ is an elementary abelian group, of order $3^s$ where $s\le \half rank(J)$.
\end{lem}
\pf
Note that $g$ acts trivially on both $F$ and $L/J$ since $L/J$ embeds in $\dual F$.
Observe that $g-1$ induces an embedding $L/F \rightarrow J$.  Furthermore, $g-1$ induces an embedding $L/(J+F) \rightarrow J/J(g-1)$, which is an elementary abelian 3-group whose rank is at most  $\half  rank(J)$ since $(g-1)^2$ induces the map $-3g$ on $J$.
\eop

\begin{lem}\labtt{sle4} $s\le rank(F)$ and $s\in \{1,2,3\}$.
\end{lem}
\pf
 If $s$ were 0, $L=J+F$ and $M$ would  be orthogonally decomposable, a contradiction. Therefore, $s\geq 1$.
The two natural maps $L\rightarrow \dg F$ and $L\rightarrow \dg J$
have common kernel $J\perp F$.
Their images are therefore  elementary abelian group of rank $s$ at most $rank(F)$ and at most $rank(J)$.  In \refpp{lmodjf}, we observed the stronger statement that $s\le \half  rank(J)$.   Since $rank(J)\ge 1$, $8=rank(J)+rank(F) >   rank(J)\ge 2s$ implies that $s\le 3$.
\eop

\begin{lem}\labtt{rank=s}
$M/((M\cap J) + F)\cong L/(J\perp F)$ is an elementary abelian
3-group of order $3^s$.
\end{lem}
\pf The quotient $L/(J+F)$ is elementary abelian by \refpp{lmodjf}.
Since $L=M+N$ and $N=Mg$, $M$ covers $L/L(g-1)$.
Since $L(g-1)\le J$, $M+J=L$  Therefore,
$L/(J\perp F)\cong (M+J)/(J+F)=(M+(J+F))/(J+F)\cong M/(M\cap (J+F))$.
The last denominator is $(M\cap J)+F$ since $F\le M$.
\eop

\begin{lem}\labtt{df}
$\dg F\cong 3^s \times 2^{rank(F)}$.
\end{lem}
\pf
Since $\drt{M}\cong E_8$ and the natural map of
$\drt{M}$ to $\dg {\drt{F}}$ is onto and has kernel $\drt{(M\cap {J} \perp {F})}$, $\dg{\drt{ F}}\cong 3^s$ is elementary abelian.
\eop


\medskip

\begin{nota}\labtt{xyk}
Let $X=M\cap J$, $Y=N\cap J$ and $K=X+Y$.  Note that $Y=Xg$ and thus
by Lemma \ref{tensorwitha2}, we have $ K\cong A_2\otimes
\left(\drt{X}\right).$

Let $\{\a, \a'\}$ be a set of fundamental roots for $A_2$ and denote
$\a''=-(\a+\a')$. Let $g'$ be the isometry of $A_2$ which is induced
by the map $\a\to \a'\to \a'' \to \a$.

By identifying $K$ with $A_2\otimes \left(\drt{X}\right)$, we may
assume $X=M\cap J= \ZZ\alpha\otimes \left(\drt{X}\right)$. Recall
that  $(x,x'g)=-\frac{1}2(x,x')$ for any $x,x'\in K$ (cf.
\refpp{tensorwitha2}). Therefore, for any $\b\in \drt{X}$, we may
identify $(\a\otimes \beta)g$  with $\a'\otimes \b = \a g' \otimes
\beta$ and  identify $Y=Xg$ with  $\a g'\otimes
\left(\drt{X}\right)$.
\end{nota}

\begin{lem}\labtt{l/fj}
We have $J=L(g-1)+K$, where $K=J\cap M+J\cap N$ as in \refpp{xyk}.
The map $g-1$ takes
$L$ onto $J$ and induces an isomorphism of $L/(J+F)$ and $J/K$, as abelian groups.
In particular, both quotients have order $3^s$.
\end{lem}
\pf
{\it Part 1: The map $g-1$ induces a monomorphism.  }
Clearly, $L(g-1)\le J$ and $g$ acts trivially on $L/L(g-1)$.
Obviously, $F(g-1)=0$. 
We also have $L\ge J+F\ge K+F$.
Since $M\cap J \le K$, $t$ acts trivially on $J/K$.
Therefore, so does $g$, whence $J(g-1)\le K$.  Since $g$ acts trivially on $L/J$, $L(g-1)^2\le K$.

Furthermore, $(g-1)^2$ annihilates $L/(F+K)$, which is a quotient of $L/F$, where the action of $g$ has minimal polynomial $x^2+x+1$.   Therefore
$L/(F+K)$ is annihilated by $3g$, so is
an elementary abelian 3-group.
We have
$3L\le F+K$.

Let $P:=\{x\in L | x(g-1)\in K\}$.   Then $P$ is a sublattice and $F+J\le P \le L$.  By coprimeness,
there are sublattices  $P^+, P^-$ so that $P^+\cap P^-=F+K$ and $P^+ + P^-=P$ and
$t$ acts on $P^\vep /(F+K)$ as the scalar $\vep = \pm 1$.
We shall prove now that $P^-=F+K$ and $P^+=F+J$.
We already know that $P^- \ge F+K$ and $P^+ \ge F+J$.

Let $v\in  P^-$ and suppose that $v(g-1)\in K$.
Then $v(g^2-g)\in K$ and this element is fixed by $t$.
Therefore, $v(g^2-g)\in ann_K( M\cap J)$.
By \refpp{tensorwitha2}, there is $u\in M\cap J$ so that  $u(g^2-g)=v(g^2-g)$.  Then $u-v \in L$ is fixed by $g$ and so $u-v\in F$.  Since $u\in K$, $v\in F+K$.
We have proved that $P^-=F+K$.

Now let $v\in P^+$.
Assume that $v\not \in F+J$.
Since $D$ acts on $L/J$ such that  $g$ acts trivially, coprimeness of $|L/J|$ and $|D/\la g\ra |$ implies that $L$ has a quotient of order 3 on which $t$ and $u$ act trivially.  Since $L=M+N$, this is not possible.   We conclude that
$F+J=P^+$.

We conclude that
$P=P^-+P^+=F+J$ and so
$g-1$ gives an embedding of  $L/(J+F)$ into
$J/K$.

{\it Part 2: The map $g-1$ induces an epimorphism.  }
We know that $L/(F+J)\cong 3^s$ and this quotient injects into $J/K$.   We now prove that
$J/K$ has  order bounded by $3^s$.

Consider the possibility that $t$ negates a nontrivial element $x+K$ of $J/K$.  By \refpp{lifteigvec}, we may assume that $xt=-x$.  But then $x\in M\cap J \le K$, a contradiction.  Therefore, $t$ acts trivially on the quotient $J/K$.
It follows that the quotient $J/K$
is covered by
$J^+(t)$.   Therefore $J/K$
embeds in the discriminant group of $K^+(t)$,
which by \refpp{tensorwitha2} is isometric to $\sqrt 3 (M\cap J)$.
Since $J/K$ is an elementary abelian $3$-group and $ \dg {M\cap J}\cong 3^s \times 2^{rank(M \cap J)}$,
the embedding takes $J/K$ to the Sylow 3-group of $\dg {M\cap J}$, which is isomorphic to $3^s$ (see \refpp{df} and use \refpp{dgxdgy}, applied to $\frac 1{\sqrt 2}M$ and the sublattices
$\frac 1{\sqrt 2}F$ and
$\frac 1{\sqrt 2}(M\cap J)$.
\eop

\begin{prop}\labtt{rootlessdih6} If $L$ is rootless and $F\ne 0$, then $s=1$ and $F\cong AA_2$.  Also, $L/(J\perp F)\cong 3$.
\end{prop}
\pf We have $s\leq 3$, so by Proposition \ref{det3ine8}, $F\cong
AA_2, EE_6$ or $AA_2\perp AA_2$ and $s=1,1$ or $2$, respectively. Note
$X=M\cap J$ is the sublattice of $M$ which is orthogonal to $F$.
Since $M\cong EE_8$, $X\cong EE_6, AA_2$ and $AA_2\perp AA_2$ if $F\cong
AA_2, EE_6$ and $AA_2\perp AA_2$, respectively.

\medskip

We shall show that $L$ has roots if $F\cong EE_6$ or $AA_2\perp AA_2$.  The conclusion in the surviving case follows from \refpp{df}.

\medskip

\textsl{Case 1:}  $F= EE_6$ and $s=1$. In this case, $X\cong Y\cong
AA_2$. Hence $K\cong A_2\otimes A_2$. As in Notation \ref{xyk}, we shall identify $X$ with $\ZZ\a\otimes A_2$ and $Y$ with $\ZZ\a g'\otimes A_2$. Then $F\perp X\cong EE_6\perp
\ZZ(\alpha \otimes A_2)$.

In this case, $|M/(F+X)|=3$ and there exist $\gamma\in (EE_6)^*$ and
$\gamma'\in (A_2)^*$ with $(\gamma,\gamma)=8/3$ and $(\gamma',
\gamma')=2/3$ such that $M=\mathrm{span}_\ZZ\{ F+X, \gamma+\alpha\otimes
\gamma'\}.$ Then $N=Mg=\mathrm{span}_\ZZ\{ F+Y, \gamma+\alpha g'\otimes
\gamma'\}$ and we have
\[
L=M+N\cong \mathrm{span}_\ZZ\{ EE_6\perp (A_2\otimes A_2),
\gamma+(\alpha\otimes \gamma'), \gamma+(\alpha g'\otimes \gamma')\}
\]

Let $\beta :=(\gamma+(\alpha\otimes \gamma')) -
(\gamma+(\alpha g'\otimes \gamma'))= (\alpha-\alpha g')\otimes
\gamma'$. Then $(\b, \b)= (\a-\a g',\a-\a g')\cdot (\g ',\g ')=6\cdot 2/3=4$.

Let $\alpha_1$ be a root of $A_2$ such that $(\alpha_1,\gamma')=-1$. Then $\a\otimes \a_1\in A_2\otimes A_2$, where $\a$ is in the first tensor factor and $\a_1$ is in the second tensor factor.
Then $(\beta, \alpha\otimes \alpha_1)= (\alpha-\alpha g',
\alpha)\cdot (\gamma', \alpha_1)= (2+1)\cdot (-1)=-3$ and the norm
of $\beta+(\alpha\otimes \alpha_1)$ is given by
\[
(\beta+(\alpha\otimes \alpha_1), \beta+(\alpha\otimes \alpha_1))=
(\beta,\beta)+(\alpha\otimes \alpha_1,\alpha\otimes
\alpha_1)+2(\beta, \alpha\otimes \alpha_1)= 4+4-6=2.
\]

Thus, $a_1=\beta+(\alpha\otimes \alpha_1)$ is a root in $J$.  So, $L$ has roots.   In
fact, we can say more.
If we take $a_2= \beta g+ \alpha g'\otimes \alpha_1$, then
$a_2$ is also a root and
\[
\begin{split}
(a_1, a_2)&=(\beta+\alpha\otimes \alpha_1, \beta g+ \alpha g'\otimes \alpha_1)\\
&=(\b, \b g)+(\b, \a g'\otimes \a_1) + (\a\otimes \a_1, \b g)+  (\a\otimes \a_1, \a g'\otimes \a_1)\\
&= -\frac{1}2(4-3-3+4)=-1.
\end{split}
\]
Thus, $a_1, a_2$ spans a sublattice $A$ isometric
to $A_2$.

\medskip

\textsl{Case 2:}  $F=AA_2\perp AA_2$ and $s=2$. In this case, $X\cong
Y\cong AA_2\perp AA_2$. Hence, $K\cong A_2\otimes(A_2\perp A_2)$.
Again, we shall identify $X$ with $\ZZ\a\otimes (A_2\perp A_2)$ and
$F+X\cong AA_2\perp AA_2 \perp \ZZ\alpha \otimes (A_2\perp A_2)\cong
AA_2\perp AA_2\perp AA_2\perp AA_2$. For convenience, we shall use a
$4$-tuple $(\xi_1, \xi_2,\xi_3,\xi_4)$ to denote an element in
$(F+K)^*\cong (AA_2)^*\perp (AA_2)^*\perp (A_2\otimes A_2)^* \perp
(A_2\otimes A_2)^*$, where $\xi_1,\xi_2\in (AA_2)^*$ and
$\xi_3,\xi_4\in (A_2\otimes A_2)^*$.

Recall that $|M/(F+X)|=3^2$ and the cosets of $M/(F+X)$ can be
parametrized by the tetracode (cf. \cite[p. 200]{cs}) whose
generating matrix is given by
\[
\begin{pmatrix}
1&1 &1 &0 \\
1&-1&0 &1
\end{pmatrix}.
\]

Hence, there exists a element $\gamma\in (AA_2)^*$ with
$(\gamma,\gamma)=4/3$ and $\gamma'\in (A_2)^*$ with $(\gamma',
\gamma')=2/3$ such that
\[
M=\mathrm{span}_\ZZ\left \{ {AA_2\perp AA_2 \perp \ZZ\alpha \otimes
(A_2\perp A_2),\atop (\gamma, \gamma, \alpha\otimes \gamma',0), (
\gamma, -\gamma, 0, \alpha\otimes \gamma')}\right \}.
\]
Therefore, we also have
\[
N=Mg= \mathrm{span}_\ZZ\left \{ {{AA_2\perp AA_2 \perp \ZZ\alpha g'
\otimes (A_2\perp A_2),}\atop{ (\gamma, \gamma, \alpha g'\otimes
\gamma',0), ( \gamma, -\gamma, 0, \alpha g'\otimes \gamma')}}\right\}
\]
and
\[
\begin{split}
&L=M+N=\\
&\mathrm{span}_\ZZ\left \{ AA_2\perp AA_2 \perp A_2\otimes (A_2\perp
A_2),(\gamma, \gamma, \alpha\otimes \gamma',0), \atop{( \gamma,
-\gamma, 0, \alpha\otimes \gamma'), (\gamma, \gamma, \alpha g'\otimes
\gamma',0), ( \gamma, -\gamma, 0, \alpha g'\otimes \gamma')}\right
\}.
\end{split}
\]

\medskip

Let $\beta_1= (\gamma, \gamma, \alpha\otimes \gamma',0) - (\gamma,
\gamma, \alpha g'\otimes \gamma',0)$ and $\beta_2=( \gamma, -\gamma,
0, \alpha\otimes \gamma')- ( \gamma, -\gamma, 0, \alpha g'\otimes
\gamma')$. Then $\b_1,\b_2\in L(g-1)\leq J$ and both have norm $4$.

Let $\alpha_1$ and $\alpha_2$ be roots in $A_2$ such that $(\alpha_1,
\gamma')=(\alpha_2,\gamma')=-1$ and $(\alpha_1,\alpha_2)=1$. Denote
\[
\begin{split}
a^1_1= \beta_1+ (0,0, \alpha \otimes \alpha_1,0),\quad  a^1_2= \beta_1 g'+
(0,0,
\alpha g'\otimes \alpha_1,0),\\
a^2_1= \beta_1+ (0,0, \alpha \otimes \alpha_2,0),\quad  a^2_2= \beta_1 g'+
(0,0,
\alpha g'\otimes \alpha_2,0),\\
a^3_1= \beta_2+ (0,0, 0, \alpha \otimes \alpha_1),\quad a^3_2= \beta_2 g'+
(0,0,0,
\alpha g'\otimes \alpha_1),\\
a^4_1= \beta_2+ (0,0, 0, \alpha \otimes \alpha_2),\quad  a^4_2= \beta_2 g'+
(0,0,0 , \alpha g'\otimes \alpha_2).\\
\end{split}
\]
Then similar to Case $1$,  we have the inner products
\[
(a^i_1, a^i_1)=(a^i_2, a^i_2)=2,\quad \text{ and }\quad  (a^i_1, a^i_2)= -1,
\]
for any $i=1,2,3,4$. Thus,
each pair $\{a^i_1,a^i_2\}$, for $i=1,2,3,4$,
spans a sublattice isometric to $A_2$. Moreover, $(a^i_k, a^j_\ell)=0$ for any $i\neq j$ and $k,\ell\in \{1,2\}$. Therefore, $J\geq A_2\perp
A_2\perp A_2\perp A_2$. Moreover, $|\mspan_\ZZ\{a^i_1,a^i_2|\,
i=1,2,3,4\}/K|=3^2$ and hence $J\cong A_2\perp
A_2\perp A_2\perp A_2$.  Again, $L$ has roots.  \eop

\begin{coro}\labtt{indexjfinl}
$|J:M\cap J + N \cap J|=3$ and $M\cap J=ann_M(F)\cong EE_6$.
\end{coro}
\pf \refpp{l/fj} and \refpp{rootlessdih6}. \eop

\begin{coro}\labtt{mjnj} (i) $M\cap J + N \cap J$ is isometric to $A_2\otimes E_6$.\\
(ii)  $L=M+N$ is unique
up to isometry.
\end{coro}
\pf For (i), use \refpp{tensorwitha2} and for (ii), use
\refpp{uniquenessofl}. \eop

\begin{lem}\labtt{glue8/3} If $v=v_1+v_2$ with $v_1\in \dual J$ and $v_2\in \dual F$, then $v_2$ has norm $\frac 43$ and $v_1$ has norm in $\frac 83+2\ZZ$.
\end{lem}
\pf
Since $3v_2\in F$, we may assume that $3v_2$ has norm 12 by \refpp{duala2} so that $v_2$ has norm $\frac 43$.  It follows that $v_1$ has norm in $\frac 83+2\ZZ$.
\eop

\subsubsection{$\Dih{6}$: Explicit glueing}
In this subsection, we shall describe the explicit glueing from
$F+M\cap J +N\cap J$ to $L$.  As in Notation \ref{xyk}, $X=M\cap J$,
$Y=N\cap J$ and $K=X+Y$. Since $F\cong AA_2$, we have $X\cong Y \cong
EE_6$ and $K\cong A_2\otimes E_6$. We also identify $X$ with $\ZZ\a
\otimes \left(\drt{X}\right )$ and $Y$ with $\ZZ\a g' \otimes
\left(\drt{X}\right )$, where $\a$ is a root of $A_2$ and $g'$ is a
fixed point free isometry of $A_2$ such that $\a g'\otimes \b=
(\a\otimes \b)g$ as described  in \refpp{xyk}. Then $F\perp X\cong AA_2
\perp \ZZ \alpha \otimes E_6\cong AA_2\perp EE_6$.

\medskip

Recall that $(AA_2)^*/AA_2 \cong \ZZ_2^2\times \ZZ_3$. Therefore,
$2(AA_2)^*/AA_2$ is the unique subgroup of order $3$ in $(AA_2)^*/AA_2$.
Similarly, $2(\ZZ \alpha\otimes E_6)^*/(\ZZ \alpha\otimes E_6)$ is the
unique subgroup of order $3$ in $(\ZZ \alpha\otimes E_6)^*/(\ZZ
\alpha\otimes E_6)\cong \ZZ_2^6\times \ZZ_3$.

\begin{nota}\labtt{gammaetc}
Since $F\perp X \leq M$ and $|M: F\perp X|=3$, there exists an
element $\mu\in F^*\perp X^*$ such that $3\mu\in F+X$ and
$M=\mathrm{span}_\ZZ\{F+X, \mu\}$.
Let $\gamma\in (AA_2)^*$ be a representative of the generator of the
order $3$ subgroup in $2(AA_2)^*/AA_2$ and $\gamma'$ a representative
of the generator of the order $3$ subgroup in $(E_6)^*/E_6 $. Without
loss, we may choose $\gamma$ and $\gamma'$ so that $(\gamma,
\gamma)=4/3$ and $(\gamma',\gamma')=4/3$.
Since the image of $\mu$ in $M/(F\perp X)$ is of order $3$, it is easy
to see that
\[
\mu \equiv \pm (\gamma + \alpha\otimes \gamma') \quad   \text{ or }\quad   \mu \equiv \pm (\gamma -
\alpha\otimes \gamma') \quad \text{ modulo } F\perp X.
\]
By replacing $\mu$ by $-\mu$ and $\gamma'$ by $-\gamma'$ if
necessary, we may assume $\mu= \gamma + \alpha\otimes \gamma'$.
Then $\nu:=\mu g= \gamma+ \alpha g' \otimes \gamma'$ and $N=
\mathrm{span}_\ZZ\{ F+Y, \nu\}$.
\end{nota}

\begin{prop}\labtt{dih6(12):m+n}
With the notation as in \refpp{gammaetc}, $ L=M+N\cong \ \mathrm{span}_\ZZ\{
AA_2\perp A_2\otimes E_6 , \gamma + \alpha\otimes \gamma',\gamma +
\alpha g' \otimes \gamma' \}. $
\end{prop}

\begin{rem} \labtt{Jindih6}
Let $\beta= (\alpha -\alpha g')\otimes \gamma'= (\gamma + \alpha\otimes
\gamma')- (\gamma + \alpha g' \otimes \gamma')$. Then $\beta \in
L(g-1)= J=ann_L(F)$ but $ \beta=(\alpha -\alpha g')\otimes \gamma'
\notin K\cong A_2\otimes E_6$. Hence $J= \mathrm{span}_\ZZ\{ \b, K
\}$ as $|J:K|=3$. Note also that $(\beta,\beta)=6\cdot 4/3=8$.
\end{rem}

\begin{lem}\labtt{jplustsqrt6duale6}
$J^+(t)=ann_J(M\cap J) \cong \sqrt{6} \dual{E_6}$.
\end{lem}

\pf By Remark \refpp{Jindih6}, we have $J= \mathrm{span}_\ZZ\{ \b, K
\}$, where $ \beta=(\alpha -\alpha g')\otimes \gamma'$ and $K= M\cap J +N\cap J\cong A_2\otimes E_6$.  Recall that $M\cap J$ is identified with $\ZZ\a\otimes E_6$ and  $N\cap J$ is identified with
$\ZZ  \a g' \otimes E_6$. Thus, by \refpp{tensorwitha2}, $ann_K(M\cap J) =\ZZ (\a g'- \a g'^2)\otimes E_6\cong \sqrt{6} E_6$. Since
$(\a, \a g'- \a g'^2)=0$,  $\b g = (\alpha g' -\alpha g'^2)\otimes \gamma'$ also annihilates $M\cap J$.
Therefore,
$
J^+(t)=ann_J(M\cap J) \geq  \mathrm{span}_\ZZ\{ ann_K(M\cap J), \b g \}.$
Since $\g' +E_6$ is a generator of $E_6^*/E_6$, we have $\mathrm{span}_\ZZ\{ E_6, \g'\}= E_6^*$ and hence
\[
\mathrm{span}_\ZZ\{ ann_K(M\cap J), \b g \}= \ZZ (\a g'- \a g'^2)\otimes \mathrm{span}_\ZZ\{ E_6, \g'\} \cong \sqrt{6}E_6^*.
\]
Note that $(\a g'- \a g'^2)$ has norm $6$. Now by the index formula, we have $\det(J^+(t))= 2^6\times 3^5 =\det (\sqrt{6}\dual{E_6})$ and thus, we have
$J^+(t)=ann_J(M\cap J) \cong \sqrt{6} \dual{E_6}$. \eop

\begin{coro}\labtt{jct} $J$ is isometric to the Coxeter-Todd lattice.  Each of these is not properly contained in an integral, rootless lattice.
\end{coro}
\pf This is an extension of the result \refpp{ctmaximal}.
Embed $J$ in $J'$, a lattice satisfying \refpp{12+12inleech2} and
embed the Coxeter-Todd lattice $P$ in a lattice $Q$ satisfying
\refpp{12+12inleech2}.
Then
both $J'$ and $Q$ satisfy the hypotheses of \refpp{12+12inleech2}, so are isometric.   Since $det(J)=det(P)=3^6$, $J\cong P$.
\eop

\subsection{$\Dih{12}$} \label{sec:dih12}

\begin{nota}\labttr{dih12setup}
Let $M, N$ be lattices isometric to $EE_8$ such that their respective
associated involutions $t_M, t_N$ generate $D\cong Dih_{12}$. Let
$h:=t_Mt_N$ and $g:=h^2$.  Let $z:=h^3$, the central involution of
$D$.  We shall make use of the $\Dih{6}$ results by working with the pair
of distinct subgroups $D_M:=D_{t_M}:=\la t_M, t_M^g \ra$ and
$D_N:=D_{t_N}:=\la t_N, t_N^g \ra$.
Note that each of these groups
is normal in $D$ since each has index 2. Define $\widetilde{M}:=Mg,
\widetilde{N}:=Ng$. If $X$ is one of $M, N$, we denote by $L_X, J_X,
F_X$ the lattices $L:=X+\widetilde{X}, J, F$ associated to the pair
$X, \widetilde{X}$, denoted $``M"$ and $``N"$ in the $\Dih{6}$ section.
We define $K_X:=(X\cap J)+(X\cap J)g$, a $D_X$-submodule of $J_X$.
Finally, we define $F:=F_D$ to be $\{x\in L| xg=g\}$ and
$J:=J_D:=ann_L(F)$. We assume $L$ is rootless.
\end{nota}

\begin{lem}\labtt{not3c}
An element of order 3 in $D$ has commutator space of dimension 12.
\end{lem}
\pf
The analysis of $\Dih{6}$ shows just two possibilities in case of no
roots (cf. Lemma \refpp{tensor} and Proposition \refpp{rootlessdih6}).  We suppose that $g-1$ has
rank 16, then derive a contradiction.

From
\refpp{tensor},
$M+Mg
\cong A_2\otimes E_8$.
From \refpp{ee8intensor}, there are only three involutions in $O(A_2\otimes E_8)$ which have negated space isometric to $EE_8$.
Therefore, $L > M+Mg$.
Since $\dg {M+Mg}\cong 3^8$, $L/(M+Mg)$ is an elementary abelian 3-group
of rank at most 4, which is totally singular in the natural $\third \ZZ
/ \ZZ$-valued bilinear form.  Note also that $L$ is invariant under the isometry
of order 3 on
$A_2\otimes E_8$ coming from the natural action of $O(A_2) \times \{ 1 \}$
on the tensor product.   We now obtain a contradiction from \refpp{overa2timese8}  since $L$ is rootless.
\eop

\begin{coro}\labtt{fmfnaa2}
$J_M\cong J_N$ has rank 12 and $F_M\cong F_N\cong AA_2$.
\end{coro}
\pf Use  \refpp{rootlessdih6} and \refpp{not3c}.   \qed

\begin{coro}\labtt{jrank6} The lattice $J$ has rank 12 and contains each of $J_M$ and $J_N$ with finite index.
The lattice $F=F_D$ has rank 2, 3 or 4 and $F/(F_M+F_N)$ is an elementary
abelian $2$-group.
\end{coro}
\pf Use \refpp{fmfnaa2} and \refpp{f/fmfnelemabel}, which implies that
$F/(F_M+F_N)$ is elementary abelian.  \eop

\medskip

\begin{nota}\labtt{utgz}
Set $t:=t_M$ and $u:=t_N$.
\end{nota}

\begin{lem}\labtt{commuting} The involutions $t^g$ and $u$ commute, and in fact $t^gu=z$.
\end{lem}
\pf
This is a calculation in the dihedral group of order 12.  See \refpp{dih12setup}, \refpp{utgz}.
We have $h=tu$  and $h^3=z$, so $t^g=t^{h^2}=utut\cdot t \cdot tutu=utututu=uh^3=uz$.
\eop

\medskip

We now study how $t_N$ acts on the lattice $J$.

\begin{lem}\labtt{submodules}
For $X=M$ or $N$, $J_X$ and $K_X$ are $D$-submodules.
\end{lem}
\pf
Clearly, $t^g$ fixes $L_M=M+Mg$, $F_M$, $J_M$
and $K_M=(M\cap J_M)+(M\cap J_M)g$, the $D_M$-submodule of $J$
generated by the negated spaces of all involutions of $D_M$.
Since $t_N=u=t^g z$, it suffices to show that the central involution $z$ fixes all these sublattices, but that is trivial.
\eop

\begin{lem}\labtt{tonj/jm}
The action of $t$ on $J/J_M$ is trivial.
\end{lem}
\pf  Use \refpp{lifteigvec} and the fact that $M\cap J \le J_M$.
\eop

\subsubsection{$\Dih{12}$: Study of $J_M$ and $J_N$}

We work out some general points about $F_M, F_N, J_M, J_N, K_M$ and $K_N$.  We continue to use the hypothesis that $L$ has no roots.

\begin{lem}\labtt{tustabilize}  As in \refpp{utgz}, $t=t_M, u=t_N$.

(i) $J(g-1)\le J_M\cap J_N$ and $J/(J_M\cap J_N)$ is an elementary
abelian 3-group of rank at most $\half rank(J)=6$.
Also, $J/(J_M\cap J_N)$ is a trivial $D$-module.

(ii) $J=J_M+J_N=L(g-1)$.
\end{lem}
\pf
(i) Observe that $g$ acts on $J/J_M$ and that $t$ acts trivially on
this quotient \refpp{tonj/jm}.  Since $g$ is
inverted by $t$, $g$ acts trivially on $J/J_M$. A similar argument
with $u$ proves $g$ acts trivially on $J/J_N$.  Therefore,
$J(g-1)\le J_M\cap J_N$.  Since $(g-1)^2$ acts on $J$ as $-3g$, (i) follows.

(ii)
We observe that $L_X(g-1)\le J_X$,  for $X=M, N$.  Since $L=M+N$,
$L(g-1)\le J_M+J_N$.   The right side is contained in $J=ann_L(F)$.
Suppose that $L(g-1) \lneq J$.
Then $J/L(g-1)$ is a nonzero 3-group which is a trivial module for $D$.
It follows that $L/(F+L(g-1))$ is an elementary abelian 3-group which
has a quotient of order 3
and  is a trivial $D$-module.  This is impossible since $L=M+N$.
So, $J=J_M+J_N=L(g-1)$.
\eop

\begin{lem}\labtt{g1}
$F=F_M+F_N$.
\end{lem}

\pf Since $F$ is the sublattice of fixed points for $g$.  Then $F$ is a direct summand
of $L$ and is $D$-invariant.
Also, $D$ acts on $F$ as a four-group and $Tel(F,D)$ has finite index in $F$.
If $E$ is any $D$-invariant 1-space in $F$, $t$ or $u$ negates $E$
(because $L=M+N$).  Therefore, $F_M+F_N$ has finite index in
$F$ and is in fact 2-coelementary abelian \refpp{f/fmfnelemabel}.

Consider the possibility that $F_M+F_N\lneq F$, i.e., that $F_M+F_N$ is not a direct summand.  Since $F_M$
and $F_N$ are direct summands, there are $\a\in F_M, \b\in F_N$ so that $\half
(\a+\b)\in
L$ but $\half \a$ and $\half \b$ are not in  $L$.   Since by \refpp{fmfnaa2} $F_M\cong
F_N\cong AA_2$, we may assume that $\a, \b$ each have norm 4.   Then by
Cauchy-Schwartz, $\half (\a+\b)$ has norm at most 4 and, if equal to 4,
$\a$ and $\b$ are equal.  But then, $\half (\a+\b)=\a\in L$, a
contradiction. \eop

\subsubsection{$\Dih{12}$: the structure of $F$}

\begin{lem}\labtt{fmcapfn}
Suppose that $F_M \ne F_N$.
Then $F_M\cap F_N$ is 0 or has rank 1 and is spanned by a vector of norm 4 or norm 12.
\end{lem}
\pf    Since $F_M$ and $F_N$ are summands, $F_M\neq F_N$ implies
$rank(F_M+F_N)\geq 3$, so $F_M + F_N$ has rank 3 or 4. Assume that
$F_M\cap F_N$ has  neither rank 0 or rank 2. Since $O(F_M) \cong
Sym_3 \times 2$ and $F_M\cap F_N$ is an eigenlattice in $F_M$ for
the involution $t_N$, it must be spanned by a norm 4 vector or is
the annihilator of a norm 4 vector. \eop

\begin{lem}\labtt{aboutgluevectors}
Suppose that the distinct vectors
$u, v$ have norm 4 and $u, v$ are in
$L\setminus (J+F)$.
We suppose that
$u\in L_M \setminus (J+F)$ and
$v\in L_N\setminus (J+F)$.
Write $u=u_1+u_2$ and $v=v_1+v_2$,
where $u_1, v_1\in \Rspan{J}$ and $u_2, v_2\in \Rspan{F}$.
We suppose that each of $u_1, u_2, v_1, v_2$ are nonzero.
Then

(i) $u_1$ and $v_1$ have norm $\frac 83$;

(ii) $u_2$ and $v_2$ have norm $\frac 43$.
\end{lem}
\pf Since $4=(u,u)=(u_1,u_1)+(u_2,u_2)$, (i) follows from (ii).
To prove (ii), use the fact that $L/(J\perp F)$ is a 3-group, a rescaled version of \refpp{duala2} and $(u_2,u_2)<(u,u)$.
\eop

\begin{lem}\labtt{aboutgluing1}
(i) Suppose that $F_M\cap F_N=\ZZ u$, where $u\ne 0$.
Let $v$ span $ann_{F_M}(u)$ and let $w$ span $ann_{F_N}(u)$.
Then
 $F=span\{u, v, w, \half (u+v), \half (u+w), \half (v+w)\}$ and $(u, u)=4$, $(v,v)=12=(w,w)$.

\medskip

(ii) If the rank of $O_3(\dg F)$ is at least 2, then

(a) $F=F_M\perp F_N$ (rank 4); or

(b) $F$ has rank 3 and the number of nontrivial  cosets of $O_3(\dg F)\cong 3 \times 3$ which have representing elements whose norms lie in  $\frac 43 + 2\ZZ$, respectively $\frac 83 +2\ZZ$
are 4 and 4, respectively.
\end{lem}
\pf
(i)   By \refpp{fmcapfn}, $u$ has norm 4 or 12.
The listed generators span $F=F_M+F_N$ since $F_M=span\{u, v, \half (u+v)\}$ and
 $F_N=span\{u, w, \half (u+w)\}$.
 If $(u,u)=12$, then $v$ and $w$ have norm $4$ and $\half (v+w)$ is a root in $F$, whereas $L$ is rootless.
 So, $(u,u)=4$ and $(v,v)=12=(w,w)$.

Now we prove (ii).
Since we have already discussed the case of $rank(F)$
equal to 2 and 4, we assume $rank(F)=3$, for which we may use the earlier results.
In the above notation, we may assume that
$F=span\{u, v, w, \half (u+v), \half (u+w), \half (v+w)\}$ and
$(u,u)=4$ and $(v,v)=12=(w,w)$.
Then
$O_3(\dg F )\cong 3^2$ and
the pair of elements $\third v, \third w$ map modulo $F$ to generators of $O_3(\dg F )$.  Since
their norms are $\frac 43, \frac 43$ and they are orthogonal, the rest of (ii) follows.
\eop

\subsubsection{$\Dih{12}$: A comparison of eigenlattices}

\begin{nota}\labttr{dih12conc}
Let $\nu$ be the usual additive $3$-adic valuation on $\QQ$, with $\nu (3^k)=k$.
Set  $P:=Mg\cap J$,  $K:=P+Pg=\ZZ [\la g \ra ] P$,
$R:=ann_{K}(P)$. Note that $P$ is the $(-1)$-eigenlattice of $t^g$ in both $J$ and $K$ while $R$ is the $(+1)$-eigenlattice of $t^g$ in $K$.
\end{nota}

We study the actions of $u$ on $J$, $J_M$,  $P, K$ and $R$.

\begin{lem}\labtt{detj=729}
$J=J_M=J_N$ and $\dg J\cong 3^6$.
\end{lem}
\pf  Since $J$ contains $J_M$ with finite index, we may use \refpp{jct}.
\eop

\begin{nota}\labttr{eigenspaces1}
Define integers $r:=rank(P^+(u))$, $s:=rank(P^-(u))$.   We have  $r=rank(R^-(u)), s=rank(R^+(u))$ and $r+s=6$.
\end{nota}

\begin{coro}\labtt{reven}  $(r, s)\in \{(2,4), (4,2), (6,0)\}$.
\end{coro}
\pf
Since $P^-+R^-$ has finite index $3^p$ in $J^-(u)=N\cap J\cong EE_6$ and $\det(P^-+R^-)=2^m3^{1+r}$, for some $m\ge 6$,
the determinant index formula implies that $r$ is even.
Similarly, we get $s$ is even.
If $r=0$, then $s=6$ and $J\cap Mg=J\cap N$, and so $z:=t^g u$ is the identity on $J$.
Since $t^g\ne u$, we have a contradiction to the $DIH_4$ theory since the common negated space for $t^g$ and $u$ is at least 6-dimensional.
So, $r\ne 0$.
\eop

\subsection{$s=0$}

\begin{lem}\labtt{ifs=0}
If $s=0$, then the pair $Mg, N$ is in case $DIH_4(15)$ or $DIH_4(16)$.
\end{lem}
\pf
In this case, $Mg\cap N$ is RSSD in $F_N$ so is isometric to $AA_2$ or $AA_1$ or $\sqrt 6 A_1$.  Then $DIH_4$ theory implies that $Mg\cap N$ is isometric to $AA_1$ or 0.
\eop

\begin{lem}\labtt{sneq0}
$s\ne 0$.
\end{lem}
\pf
Suppose that $s=0$.  Then $u$ acts as 1 on $Mg\cap J$.
The sublattice $(Mg\cap J) \perp (N\cap J)$ has determinant $2^{12}3^2$ is contained in $Tel(J, u)$, which has determinant $2^{12}det(J)=2^{12}3^6$.
Since $Mg\cap J=J^-(t^g)$, it follows from determinant considerations that $N\cap J$ is contained with index $3^2$ in $J^+(t^g)$.   Since $s=0$, $J^+(t^g) \le J^-(u)=N\cap J$ and we have a contradiction.
\eop

\subsection{$s\in \{2, 4\}$}

\begin{lem}\labtt{eigenspaces3}
If $s>0$ and $det(P^-(u))$ is not a power of 2, then $s=2$ and $P^-(u)\cong 2A_2$  and the pair $Mg, N$ is in case $DIH_4(12)$.
\end{lem}
\pf
Since $P^-(u)$ is RSSD in $P$ and $det(P^-(u))$ is not a power of $2$, \refpp{involsoe6} implies that
$P^-(u) \cong AA_5$ or $2A_2$ or $(2A_2)(AA_1)^m$.
Since $s=rank(P^-(u))$ is 2 or 4,
we have $P^-(u) \cong 2A_2$ or $(2A_2)(AA_1)^2$.
By $DIH_4$ theory, $Mg\cap N\cong DD_4$.
Since $P^-(u)$ is contained in $Mg\cap N$,
$P^-(u)\cong 2A_2$.
\eop

\begin{lem}\labtt{eigenspaces4}
(i)
If $s>0$ and $det(P^-(u))$ is a power of 2, then $s=2$ and  $P^-(u)\cong AA_1\perp AA_1$ or $s=4$ and  $P^-(u) \cong DD_4$.

(ii) If $P^-(u)\cong DD_4$, the pair $Mg, N$ is in case $DIH_4(12)$.

(iii) If $P^-(u)\cong AA_1\perp AA_1$, the pair  $Mg, N$ is in case $DIH_4(14)$.  In particular, $F_M\cap F_N=0$ and so $F=F_M \perp F_N$.
\end{lem}
\pf (i)
Use
\refpp{involsoe6} and evenness of $s$.

(ii) This follows from $DIH_4$ theory since $Mg\cap N$ contains a copy of $DD_4$ and $Mg\ne N$.

(iii) Since $dim(P^-(u))=2$, it suffices by $DIH_4$ theory to prove that $dim(Mg\cap N)\ne 4$.  Assume by way of contradiction that $dim(Mg\cap N) =  4$.   Then $Mg\cap N\cong DD_4$ and $rank(Mg\cap N\cap F)=2$.   This means that $F=F_M=F_N\cong AA2$.
Therefore, $Mg\cap N\cong DD_4$ contains the sublattice
$P^-(u)\perp F\cong  AA_1\perp AA_1\perp AA_2$, which is impossible.
\eop

\begin{lem}\labtt{int2a2}
Suppose that $P^-(u)$ is isometric to $2A_2$.   Then $rank(L)=14$ and $L$ has roots.
\end{lem}
\pf
We have $det(P^-(u)\perp R^-(u))=2^23\cdot 2^63^4=2^8 3^5$.   Since $N\cap J$ covers $J/K_M$, the determinant formula implies that $|J:K_M|=3^2$ and so $det(J)=3^4$.
Now use \refpp{jct}.
\eop

\begin{lem}\labtt{pminusunota1a1}
$P^-(u)$ is not isometric to $AA_1\perp AA_1$.
\end{lem}
\pf
Suppose $P^-(u)\cong AA_1\perp AA_1$.
Then, $P^+(u)\cong \sqrt 2 Q$,
where $Q$ is the rank 4 lattice which is described in \refpp{xq}.
Also, $R^+\cong \sqrt 6 A_1^2$ and $R^-\cong \sqrt 6 Q$.
Then $P^-(u)\perp R^-(u)$ embeds in $EE_6$.   Since sublattices of $E_6$ which are isometric to $A_1^2$ are in a single orbit under $O(E_6)$, it follows that $\sqrt 3 Q$ embeds in $Q$.  However, this is in contradiction with \refpp{transxq2}.
\eop

\medskip

To summarizes our conclusion, we have the proposition.

\begin{prop}
$P^-(u)=Mg\cap N\cong DD4$ and the pair is in case $\dih{4}{12}$.
\end{prop}

\subsection{Uniqueness of the case $DIH_{12}(16)$}

As in other sections, we aim to use \refpp{uniquenessofl} for the case \refpp{eigenspaces4}(ii).

The input $M, N$ determines the dihedral group $\la t, u \ra$ and therefore
$Mg$ and $Mg+N$.  By $DIH_4$ theory, the isometry type of $Mg+N$ is determined up to isometry.  Since $Mg+N$ has finite index in $M+N$, $M+N$ is determined  by \refpp{uniquenessofl}.
Thus, Theorem \refpp{thmDih612} is proved.

\section{$\Dih{10}$ theory}

\begin{nota}\labttr{dih10setup}  Define $t:=t_M$, $h:=t_Mt_N$.  We suppose
$h$ has order 5. Let $g:= h^3$. Then $g$ also has order $5$ and $D:=< t_M, t_M>= <t, g>$. In addition, we have  $N=Mg$. Define $F:=M\cap N$,  $J:=ann_L(F)$.  Note that
$F$ is the common negated lattice for $t_M$ and $t_N$ in $L$, so is
the fixed point sublattice for $g$ and is a direct summand of $L$
\refpp{mnfsummands}. \end{nota}

\begin{de} \labttr{s:dih10} Define the integer $s$ by $5^s:=|L/(J+F)|$.
\end{de}

\begin{lem}\labttr{interprets=0}
Equivalent are (i) $L=J+F$; (ii) $s=0$; (iii) $F=0$; (iv) $J=L$.
\end{lem}

\pf Trivially,  $L=J+F$ and $s=0$ are equivalent.  These
conditions follow if $F=0$ or if $J=0$ (but the latter does not happen
since $g$ has order 5).  If $L=J+F$ holds, then $M=(J\cap M) \perp F$
which implies that $F=0$ and $J=L$ since $M\cong EE_8$ is orthogonally
decomposable.  \eop

\begin{lem}\labtt{g-1}
(i) $g$ acts trivially on both $F$ and $L/J$.

(ii) $g-1$ induces an embedding $L/F\rightarrow J$.

(iii)  $g-1$ induces an embedding $L/(J+F) \rightarrow J/J(g-1)$,
whose rank is  at most $\fourth rank(J)$ since $(g-1)^4$ induces the
map $5w$ on $J$, where $w=-g^2+g-1$ induces an invertible linear
map on $J$.

(iv)  $s \le \fourth rank(J)$, so that $s=0$ and $F=0$ or $s\in
\{1,2,3,4\}$ and $F\ne 0$.

(v) The inclusion $M \le L$ induces an isomorphism $M/((M\cap
J)+F)\cong L/(J+F)\cong 5^s$, an elementary abelian group.
\end{lem}
\pf (i) and
(ii) are trivial.

(iii) This is equivalent to some known behavior in the ring of
integers $\ZZ [e^{2\pi i /5}]$, but we give a self-contained proof
here. We calculate
$(g-1)^4=g^4-4g^3+6g^2-4g+1=(g^4+g^3+g^2+g+1)+5w$, which in $End(J)$
is congruent to $5w$.   Note that the images of $g+1$ and $g^3+1$
are non zero-divisors (e.g., because
$(g+1)(g^4-g^3+g^2-g+1)=g^5+1=2$, and 2 is a non zero-divisor) and
are associates in $End(J)$ so that their ratio $w$ is a unit. For
background, we mention \cite{gh}.

(iv) The Jordan canonical form for the action of $g-1$ on $J/5J$ is
a direct sum of degree 4 indecomposable blocks, by (iii), since $(g-1)^4$ has determinant $5^{rank(J)}$.
Since the
action of $g$ on $L/J$ is trivial, $s\le \fourth rank(J)$. Since
$rank(J)\le 16$, $s \le 4$.  For the case $s=0$, see \refpp{interprets=0}.

(v) Since $N=Mg$, $N$ and $M$ are congruent modulo $J$.  Therefore
$L=M+N=M+J$ and so $5^s\cong L/(J+F)=(M+J)/(J+F)=(M+(J+F))/(J+F)
\cong M/(M\cap (J+F))$  (by a basic isomorphism theorem) and this
equals $M/((M\cap J)+F)$   (by the Dedekind law).

Since $L(g-1)\le J$, $(g-1)$ annihilates $L/(J+F)$. Since
$(g-1)^4$ takes $(L/(J+F))$ to $5(L/(J+F))$, it follows that $5L\le
J+F$.   That is, $L/(J+F)$ is an elementary abelian 5-group. \eop

\begin{lem}\labtt{fdih10}
$s=0$, $1$, $2$ or $3$ and $F=M\cap N\cong 0$, $AA_4$, $\sqrt{2}\mathcal{M}(4,25) $
or $\sqrt{2}A_4(1)$.
\end{lem}
\pf We have that  $\drt{M}\cong E_8$.  The natural map of $\drt{M}$
to $\dg {\drt{F}}$ is onto and has kernel $\drtp{(M\cap J) \perp
F)}$. Therefore, $\dg {\drt{F}} \cong 5^s$ is elementary abelian.
Now apply \refpp{det5ine8} to get the possibilities for $\drt{F}$
and hence for $F$. Note that $M=N$ is impossible here, since $t_M\ne
t_N$. \eop

\subsection{$\Dih{10}$: Which ones are rootless?}

From Lemma \ref{fdih10}, $s=0$, $1$, $2$ or $3$. We shall eliminate the
case $s=1$, $s=2$ and $s=3$, proving that $s=0$ and $F=0$.

\begin{lem}
If $L=M+N$ is integral and rootless, then $F=M\cap N = 0$.
\end{lem}
\pf   By Lemma \ref{fdih10}, we know that $M\cap N\cong 0$, $AA_4$
$\sqrt{2} \mathcal{M}(4,25) $ or $\sqrt{2}A_4(1)$ since $M\neq N$. We shall eliminate the cases
$M\cap N\cong AA_4$, $\sqrt{2} \mathcal{M}(4,25) $ and $\sqrt{2}A_4(1)$.

\medskip

\emph{Case:} $F= M\cap N \cong AA_4$.  In this case, $M\cap J\cong
N\cap J\cong AA_4$. Therefore, there exist $\a\in F^*$ and $\b\in
(M\cap J)^*$ such that $M=\mathrm{span}_\ZZ\{ F+(M\cap J), \a+\b\}.$
Without loss, we may assume $(\a,\a)=12/5$ and $(\b,\b)=8/5$. Let
$\g=\b g$. Then, $(\a+\b)g=\a+\g \in N$ and we have
$N=\mathrm{span}_\ZZ\{ F+(N\cap J), \a+\gamma\}.$  Since $L$ is
integral and rootless and since $\a+\b\in L$ has norm $4$, by
\refpp{norm4gorbit},
\[
 0 \geq (\a+\b, (\a+\b)g) = (\a+\b,
 \a+\g)=(\a,\a)+(\b,\g)=\frac{12}5 +(\b,\g).
 \]
Thus, we have $(\b,\g)\leq -\frac{12}5$.  However, by the Schwartz
inequality, $$|(\b,\g)|\leq \sqrt{(\b,\b)(\g,\g)} =\frac{8}5, $$ which
is a contradiction.

\medskip

\emph{Case:} $F= M\cap N \cong \sqrt{2} \mathcal{M}(4,25) $.   In this case,
$M\cap J\cong N\cap J\cong \sqrt{2} \mathcal{M}(4,25) $, also. Let
$\sqrt{2}u,\sqrt{2}v,\sqrt{2}w,\sqrt{2}x$ be a set of orthogonal
elements in $F\cong \sqrt{2} \mathcal{M}(4,25) $ such that their norms are
$4,8,20,40$, respectively (cf. \refpp{m(4,25)}).  Let $\sqrt{2}u',
\sqrt{2}v',\sqrt{2}w'$, $\sqrt{2}x'$ be a sequence of pairwise orthogonal elements in
$M\cap J$ such that their norms are $4,8,20,40$, respectively.
By the construction in
\refpp{rank4det25} and the uniqueness assertion, we may assume that  the element $\g=\frac{\sqrt{2}}5(w+x+x')$ is in
$M$. Since $\gamma$ has norm $4$, by \refpp{norm4gorbit},
\[
 0 \geq (\g, \g g)= \frac{2}{25} (w+x+x', w+x+x'g)=
 \frac{60}{25}  +\frac{2}{25}(x', x'g).
\]
Thus, we have $(x', x'g)\leq - 30$.  By the Schwartz inequality,
\[
|(x',x' g)|\leq \sqrt{(x',x')(x' g,x' g)} =20,
\]
which is again a  contradiction.

\medskip

\emph{Case:} $F= M\cap N \cong \sqrt{2}A_4(1)$.   Since $F$ is a direct summand of $M$ and $N$, we have $M\cap J\cong N\cap J\cong \sqrt{2}A_4(1)$  by \refpp{a41ine8}.  Recall that $\dual{(A_4(1))} \cong \frac{1}{\sqrt{5}}A_4$.

By the construction in \refpp{a41ine82}, there exists $\a\in 2\dual{F}$ with $(\a,\a)=2\times 8/5=16/5$ and
$\a_M\in 2\dual{(M\cap J)}$ with $(\a_M, \a_M)=2\times 2/5=4/5$ such that $y=\a+\a_M\in M$.

Since $(y,y)=4$, by \refpp{norm4gorbit},
\[
0\geq (y, y g)= (\a+\a_M, \a+\a_M g) = (\a, \a)+(\a_M, \a_M g)
\]
and we have $(\a_M, \a_M g) \leq -(\a, \a)=-16/5$. However, by the Schwartz inequality,
$$|(\a_M, \a_M g)| \leq \sqrt{ (\a_M, \a_M)(\a_M g, \a_M g)} = 4/5, $$ which is a contradiction.
\eop

\subsection{$\Dih{10}$: An orthogonal direct sum}

For background, we refer to \refpp{rank4det4}, \refpp{rank4det125}, \refpp{a4(1)} -- \refpp{allabouta4(1)}.
Our goal here is to build up an orthogonal direct sum of four copies
of $AA_4$ inside $L$.  We do so one summand at a time.  This direct
sum shall determine $L$ (see the following subsection).

\begin{nota}\labttr{orbitdotproducts}
Define $Z(i):=\{x\in M | (x,x)=4, (x,xg)=i\}$. Note that $(x,xg)=-3,-2,-1, 0$, or $1$ by Lemma \ref{norm4gorbit}.
\end{nota}

\begin{lem}\labtt{uvg=ugv}
For $u, v\in M$, $(u,vg)=(u,vg^{-1})=(ug,v)$.
\end{lem}
\pf Since $t$ preserves the form,
$(u,vg)=(ut,vgt)=(-u,vtg^{-1})=(-u,-vg^{-1}) =(u,vg^{-1})$.  This
equals $(ug,v)$ since $g$ preserves the form. \eop

\begin{lem}\labtt{z-1nonempty}   If $u, v, w$ is any set of norm 4 vectors so that $u+v+w=0$, then
one or three of $u, v, w$  lies in $Z(-2)\cup Z(0)$.  In particular, $Z(-2)\cup Z(0)\ne
\emptyset$.
\end{lem}
\pf Suppose that we have norm 4 vectors $u, v, w$ so that $u+v+w=0$.
Then
$0=(u+v+w,ug+vg+wg)=(u,ug)+(v,vg)+(w,wg)+(u,vg)+(ug,v)+(u,wg)+(ug,w)+(v,wg)+(vg,w)\equiv
(u,ug)+(v,vg)+(w,wg) (mod\, 2)$, by \refpp{uvg=ugv}, whence evenly
many of $(u,ug)$, $(v,vg)$, $(w,wg)$ are odd. \eop

\medskip

Now we look at $D$-submodules of $L$ and decompositions.

\begin{de}\labttr{equivreln4}
Let $M_4=\{\a\in M|\ (\a,\a)=4\}$. Define a partition of $M_4$ into sets $M_4^1:=\{ \a \in M_4 |\, \a \DD
\cong A_4(1)\}$ and $M_4^2:=\{ \a \in M_4 |\, \a \DD \cong AA_4 \}$ (cf.
\refpp{norm4gorbit}). For $\a, \b \in M_4^2$, say that  $\a$ and $\b$ are equivalent
if and only if $\a \DD = \b \DD$. Define the partition
$N_4=N_4^1\cup N_4^2$ and equivalence relation on $N_4^2$ similarly.
\end{de}

\begin{rem}\labttr{tmaps}
The linear maps $g^i+g^{-i}$ take $M$ into $M$ since they commute
with $t$.  Also, $g^2+g^3$ and $g+g^4$ are linear isomorphisms of
$M$ onto $M$ since their product is $-1$.  Note that they may not
preserve inner products.
\end{rem}

\begin{lem}\labttr{m41empty}
$M_4=M_4^2$ and $M_4^1=\emptyset$.
\end{lem}
\pf Supposing the lemma to be false, we take $\a \in M_4^1$.   Then
the norm of $\a (g^2+g^3)$ is $4+2(\a g^2, \a g^3)+4=8-2=6$ (cf.
\refpp{norm4gorbit}), which is impossible since $M\cong EE_8$ is
doubly even. \eop

\begin{lem}
 Let $\a\in M_4$. Then $M\cap \a\ZZ[D]\cong AA_1^2$.
\end{lem}

\pf Let $\a\in M_4$. Then by \refpp{m41empty}, $\a\ZZ[D]\cong AA_4$.
In this case, we have either (1) $(\a, \a g)=-2$ and $(\a, \a g^2)=0$ or (2)
$(\a, \a g)=0$ and $(\a, \a g^2)=-2$.

\medskip

In case (1), we have $\a(g^2+g^3)\in M_4$ and $\a$ and  $\a(g^2+g^3)$ generate a sublattice of type $AA_1^2$ in $M\cap \a\ZZ[D]$.  Similarly, $\a g$, $\a(g^2+g^3)g$ generate a sublattice of type $AA_1^2$ in $N\cap \a\ZZ[D]$. Since $M\cap N=0$, we have $rank(M\cap \a\ZZ[D])=rank(N\cap \a\ZZ[D])=2$. Moreover, $\{\a, \a(g^2+g^3), \a g, \a(g^3+g^4)\}$ forms a $\ZZ$-basis of an $AA_4$-sublattice of $\a\ZZ[D]\cong AA_4$. Thus, $\{\a, \a(g^2+g^3), \a g, \a(g^3+g^4)\}$ is also a basis of $\a \ZZ[D]$ and $\mspan_\ZZ\{\a, \a(g^2+g^3) \}$ is summand of $\a \ZZ[D]$. Hence, $M\cap \a\ZZ[D]=\mspan_\ZZ\{\a, \a(g^2+g^3) \}\cong AA_1^2$ as desired.

\medskip

In case (2), we have  $\a(g+g^4)\in M_4$ and thus $M\cap \a\ZZ[D]=\mspan_\ZZ\{\a, \a(g+g^4) \}\cong AA_1^2$ by an argument as in case (1). \eop

\begin{lem}\labtt{gramaa4}
Suppose that $\a \in M_4^2$, $\b \in N_4^2$ and $\a \DD = \b \DD$.
Let the equivalence class of $\a$ be $\{ \pm \a, \pm \a' \}$ and let
the equivalence class of $\b$ be $\{ \pm \b, \pm \b' \}$.  After
interchanging $\b$ and one of $\pm \b'$ if necessary, the Gram matrix of
$\a, \a', \b, \b'$ is $$2\begin{pmatrix} 2&0&0&1\cr 0&2&1&-1\cr
0&1&2&0\cr 1&-1&0&2
\end{pmatrix}=\begin{pmatrix} 4&0&0&2\cr 0&4&2&-2\cr
0&2&4&0\cr 2&-2&0&4
\end{pmatrix}.$$
\end{lem}
\pf We think of $A_4$ as the 5-tuples in $\ZZ^5$ with zero coordinate
sum. Index coordinates with integers mod 5: 0, 1, 2, 3, 4.  Consider
$g$ as addition by 1 mod 5 and $t$ as negating indices modulo 5. We
may take $\a :=\sqrt 2 (0,1,0,0,-1), \a' :=\sqrt 2 (0,0,1,-1,0)$.  We define $\b:=\a
g, \b':=\a' g$.  The computation of the Gram matrix is
straightforward. \eop

\begin{lem}\labtt{predirectsum4}   Let $m \ge 1$.  Suppose that $U$ is a rank $4m$ $\DD$-invariant sublattice of $L$ which is generated as a $\DD$-module by $S$, a sublattice of $U\cap M$ which is isometric to $AA_1^{2m}$.
Then $ann_M(U)$ contains a sublattice of type $AA_1^{8-2m}$.
\end{lem}
\pf 
We may assume that $m\le 3$.
Since $t$ inverts $g$ and $g$ is fixed point free on $L$, $U^-(t)=U\cap M$ has rank $2m$.
Let $S$ be a sublattice of $U\cap M$ of type $AA_1^{2m}$.  Then $S$ has finite index in $U \cap M$.
Let $W:=ann_L(U)$, a direct summand of $L$ of rank $16-4m$.
The action of $g$ on $W$ is fixed point free and $t$ inverts $g$ under conjugation,
so $W\cap M=ann_M(U)$ is a direct summand of $M$ of rank $8-2m$.
It is contained in hence is equal to the annihilator in $M$ of $S$, by rank considerations,
so is isometric to $DD_6, DD_4, AA_1^4$ or $AA_1^2$.
Each of these lattices contains a sublattice of type $AA_1^{8-2m}$.
\eop

\begin{coro}\labtt{directsum4}
$L$ contains an orthogonal direct sum of four $D$-invariant lattices, each isometric to $AA_4$.
\end{coro}
\pf
We prove by induction that for $k=0,1,2,3,4$,
$L$ contains an orthogonal direct sum of $k$ $D$-invariant lattices, each isometric to $AA_4$.
This is trivial for $k=0$.  If $0\le k \le 3$, let $U$ be such an orthogonal direct sum of $k$ copies of $AA_4$.   Then $M\cap U\cong AA_1^{2k}$ and thus $ann_M(U)$ contains a norm 4 vector, say $\a$.  By \refpp{m41empty}, $\a \DD \cong AA_4$.   So,
$U\perp  \a \DD$ is an orthogonal direct sum of $(k+1)$ $D$-invariant lattices, each isometric to $AA_4$.  \eop

\begin{coro}
$L=M+N$ is unique up to isometry.
\end{coro}

\pf Uniqueness follows from the isometry type of $U$ (finite index
in $L$) and \refpp{uniquenessofl}. We take the finite index sublattices $M_1:=M\cap U$ and $N_1:=N\cap U$ and use \refpp{gramaa4}. An alternate proof is given by the glueing in \refpp{gluedih10} \eop

\subsection{$\Dih{10}$: From $AA_4^4$ to $L$}

We discuss the glueing from  a sublattice $U=U_1\perp U_2\perp U_3
\perp U_4$, as in \refpp{directsum4} to $L$.  We assume that each
$U_i$ is invariant under $D$.

By construction, $M/(M\cap U)\cong 2^4$, $M\cap U\cong AA_1^8$.  A
similar statement is true with $N$ in place of $M$. Since $L=M+N$,
it follows that $L/U$ is a 2-group.  Since $g$ acts fixed point
freely on $L/U$, $L/U$ is elementary abelian of order $2^4$ or
$2^8$.  Also, $L/U$ is the direct sum of $C_{L/U}(t)$ and
$C_{L/U}(u)$, and each of the latter groups is elementary abelian of
order $|L:U|^\half$.  So $|L:U|^\half = |M:M\cap U|^2= (2^4)^2=2^8$.
Therefore,  $det(L)=5^4$ and the
Smith invariant sequence of $L$ is $1^{12}5^4$.

\begin{prop}\labtt{gluedih10}
The glueing from $U$ to $L$ may be identified with the direct sums
of these two glueings from $U\cap M$ to $M$ and $U\cap N$ to $N$. Each
glueing is based on the extended Hamming code with parameters
[8,4,4] with respect to the orthogonal frame.
\end{prop}

\subsection{$\Dih{10}$: Explicit glueing and tensor products}

In this section, we shall give the glue vectors from $U=U_1\perp U_2\perp U_3
\perp U_4$ to $L$ explicitly in Proposition \ref{dih10glue} (cf. Proposition \ref{gluedih10}). We  also show that $L$ contains a sublattice
isomorphic to a tensor product $A_4 \otimes A_4$.

\begin{nota}\labtt{aa'} 
Recall that $M\cap U_i\cong AA_1\perp AA_1$ for $ i=1,2,3,4$.  Let $\a_i\in M_4 \cap
U_i, i=1,2,3,4,$ such that $(\a_i,\a_ig)=-2$. Note that such $\a_i$
exists because if $(\a_i, \a_i g)\neq -2$, then $(\a_i, \a_i g)=0$
and $(\a_i, \a_i g^2)=-2$. In this case, $\tilde{\a}_i=
\a_i(g+g^4)\in M_4\cap U_i$  and $(\tilde{\a}_i, \tilde{\a}_i
g)=-2$ \refpp{tmaps}.

Set $\a_i':= \a_i(g^2+g^3)$ for $i=1,2,3,4$. Then $\a_i'\in M_4 \cap
U_i$ and  $M\cap U_i =\mathrm{span}_\ZZ\{\a_i, \a_i'\}$ \refpp{tmaps}.
\end{nota}

\begin{lem}\labtt{inn}
Use the same notation as in \refpp{aa'}. Then for all $i=1,2,3,4,$ we have
$ (\a_i, \a_ig)=-2$, $(\a_i, \a_ig^2)=0$, $(\a_i, \a_i')=0$, $(\a_i', \a_i'g)=0$  and $(\a_i, \a_i'g)=-2.$
\end{lem}

\pf By definition, $(\a_i, \a_i g)=-2$ and $(\a_i, \a_i g^2)=0$. Thus, we have
$(\a_i', \a_i)=(\a_i(g^2+g^3), \a_i)=0$. Also,
\[
\begin{split}
(\a'_i, \a'_i g)& =(\a_i(g^2+g^3), \a_i(g^2+g^3)g)\\
 &= (\a_ig^2, \a_ig^3)+
(\a_ig^2, \a_ig^4)+(\a_ig^3, \a_i g^3)+(\a_ig^3, \a_i g^4)\\
&=-2+0+4-2=0
\end{split}
\]
and
\[
(\a_i, \a'_i g)= (\a_i, \a_i(g^2+g^3)g)= (\a_i, \a_i g^3)+(\a_i,
\a_ig^4)=0-2=-2.
\]
\eop

\begin{rem}\labtt{MminusU}
Since $M$ and $U$ are doubly even and since
$\frac{1}{\sqrt{2}}(U\cap M)\cong (A_1)^8$ and $(A_1)^*=\frac{1}2 A_1$, for any
$\b\in M \setminus (U\cap M) $,
\[ \b=\sum_{i=1}^4( \frac{b_i}2 \a_i +
\frac{b_i'}2 \a_i')\quad \text {where } b_i, b_i'\in \ZZ \text{ with
some } b_i, b_i' \text{ odd.}
\]
\end{rem}

\begin{lem}\labtt{glueaa4}
Let $\b\in M \setminus (U\cap M)$ with $(\b,\b)=4$. Then, one of the
following three cases holds.
\begin{enumerate}
\item[(i)] $|b_i|=1$ and  $b_i'=0$ for all
$i=1,2,3,4$; 

\item[(ii)]  $|b'_i|=1$ and $b_i=0$ for all
$i=1,2,3,4$; 

\item[(iii)] There exists a 3-set $\{i,j,k\} \subset \{1,2,3,4\}$ such that
$b_i^2=b_j^2=1$ and $b_i'^2=b_k'^2=1$.
\end{enumerate}
\end{lem}

\pf Let $\b=\sum_{i=1}^4( \frac{b_i}2 \a_i + \frac{b_i'}2 \a_i')\in
M \setminus (U\cap M)$ with $(\b,\b)=4$. Then we have $\sum_{i=1}^4 ( b_i^2
+ b_i'^2) =4$. Since no $|b_i|$ or $|b_i'|$ is greater than 1 (or else no $b_i$ or $b_i'$ is odd), $b_i, b_i'\in \{ -1, 0, 1\}$.
Moreover,  $(\b, \b g)=0$ or $-2$ since $\b\in M_4$.
By \eqref{inn},
\[
\begin{split}
(\b,\b g)
= &\left(\sum_{i=1}^4 ( \frac{b_i}2 \a_i + \frac{b_i'}2 \a_i'), \sum_{i=1}^4( \frac{b_i}2 \a_i g + \frac{b_i'}2 \a_i'g)\right )\\
= &\frac{1}4  \sum_{i=1}^4 \left( b_i^2 (-2) + 2 b_ib_i' (-2) \right
)=-\frac{1}2  \sum_{i=1}^4 \left( b_i^2  + 2 b_ib_i' \right
).
\end{split}
\]

If $(\b, \b g)=-2$, then
\[
\sum_{i=1}^4 \left( b_i^2  + 2 b_ib_i' \right)=4 =\sum_{i=1}^4 ( b_i^2
+ b_i'^2) .
\]
and hence we have \textbf{(*)}  $\sum_{i=1}^4 b_i'(b_i'- 2 b_i)=0$.

\medskip

Set $k_i:=b_i'(b_i'- 2 b_i)$. The values of $k_i$, for all $b_i, b_i'\in \{-1,0,1\}$, are listed in Table \ref{kbb}.

\begin{table}[bht]
\caption{ Values of $k_i$}
\begin{center}
\begin{tabular}{|c|c|c|c|c|c|c|c|}
\hline
   $b_i'$&  $0$& $-1$& $-1$& $-1$ &$1$&$1$&$1$  \cr \hline
   $b_i$& $-1,0,1$& $-1$& $0$& $1$ & $-1$& $0$& $1$ \cr \hline
   $k_i=b_i'(b_i'- 2 b_i)$&$0$& $-1$& $1$& $3$ &$3$&$1$&$-1$  \cr \hline
\end{tabular}
\end{center}
\label{kbb}
\end{table}%

Note that $k_i= 0, \pm 1$ or $3$ for all $i=1,2,3,4$. Therefore, up to the order of the indices, the values for $(k_1, k_2,k_3,k_4)$ are $(3, -1,-1,-1)$, $(1,-1, 1,-1)$, $(1,-1,0,0)$ or $(0,0,0,0)$.

However, for $(k_1, k_2,k_3,k_4)=(3, -1,-1,-1)$ or $(1,-1, 1,-1)$, $b_i'^2=b_i^2=1$ for all $i=1,2,3,4$ and then $\sum_{i=1}^4 ( b_i^2
+ b_i'^2)=8>4$. Therefore, $(k_1, k_2,k_3,k_4)=(1,-1,0,0)$ or $(0,0,0,0)$.

If $(k_1, k_2,k_3,k_4)=(1,-1,0,0)$, then we have, up to order,  $(b_1')^2=1$ (whence $k_1=1$), $b_1=0$, $b_2'=b_2=\pm 1$ (whence $k_2=-1$) and $b_3'=b_4'=0$. Since $\sum_{i=1}^4 ( b_i^2 + b_i'^2) =4$, $b_3^2+b_4^2=1$ and hence we have (iii).

If
$k_i=b_i'(b_i'- 2 b_i)=0$ for all $i=1,2,3,4$, then  $b_i'=0$ for all $i$ and we have (i). Note that  if $b_i'\neq 0$, $b_i'- 2 b_i\neq 0$.

Now assume $(\b, \b g)=0$. Then $\sum_{i=1}^4 \left( b_i^2  + 2 b_ib_i' \right)=0$. Note that this equation is the same as the above equation (*) in the case for $(\b, \b g)=-2$ if we replace $b_i$ by $b_i'$ and $b_i'$ by $-b_i$ for $i=1,2,3,4$. Thus, by the same argument as in the case for $(\b, \b g)=-2$, we have either $b_i^2  + 2 b_ib_i'=0$ for all $i=1,2,3,4$ or  $b_1^2  + 2 b_1b_1'=1, b_2^2  + 2 b_2b_2'=-1$, and $b_3=b_4=0$. In the first case, we have $b_i=0$ and $b_i'^2=1$ for all $i=1,2,3,4$, that means (ii) holds. For the later cases, we have $b_1=1, b_1'=0$, $b_2=1, b_2'=-1$, $b_3=b_4=0$, and $b_3'^2+b_4'^2=1$ and thus (iii) holds.
\eop

\begin{prop}\labtt{dih10glue}
By rearranging the indices if necessary, we have
\[
M= \mathrm{span}_\ZZ \left \{ {M\cap U, \frac{1}2
(\a_1+\a_2+\a_3+\a_4), \frac{1}2 (\a_1'+\a_2'+\a_3'+\a_4'), \atop
\frac{1}2 (\a_1+\a_2+\a_2'+\a_4'),\frac{1}2
(\a_1+\a_3+\a_2'+\a_3')}\right\}
\]
and
\[
N=Mg = \mathrm{span}_\ZZ \left \{ {N\cap U, \frac{1}2
(\b_1+\b_2+\b_3+\b_4), \frac{1}2 (\b_1'+\b_2'+\b_3'+\b_4'), \atop
\frac{1}2 (\b_1+\b_2+\b_2'+\b_4'),\frac{1}2
(\b_1+\b_3+\b_2'+\b_3')}\right\},
\]
where $\b_i=\a_i g$ and $\b_i'= \a_i'g$ for all $i=1,2,3,4$.
\end{prop}

\pf  By \refpp{glueaa4}, the norm $4$ vectors in $M\setminus (U\cap M)$ are of the form
\[
\frac{1}2 ( \pm \a_1\pm \a_2\pm \a_3\pm \a_4),\ \frac{1}2 ( \pm \a_1'\pm \a_2'\pm \a_3'\pm \a_4'),\ \text{ or } \frac{1}2 ( \pm \a_i\pm \a_j\pm \a_j'\pm \a_k'),
\]
where $i,j,k$ are distinct elements in $\{1,2,3,4\}$.

Since $M\cong EE_8$ and $U\cap M \cong (AA_1)^8$, the cosets of $M/(U\cap M)$
can be identified with the codewords of the Hamming [8,4,4] code $H_8$.

Let $\varphi: M/(U\cap M) \to H_8$ be an isomorphism of binary codes. For any
$\b\in M$, we denote the coset $\b+U\in M/(U\cap M)$ by $\bar{\b}$.
We shall also arrange the index set such that the first 4 coordinates correspond to the coefficient of $\half \a_1, \half\a_2, \half\a_3$ and $\half\a_4$ and the last 4 coordinates correspond to the coefficient of $\half\a_1',\half \a_2',\half \a_3'$ and $\half \a_4'$.

Since $(1, \dots , 1)\in H_8$, we have
\[
\frac{1}2 (\a_1+\a_2+\a_3+\a_4+\a_1'+\a_2'+\a_3'+\a_4')\in M.
\]

We shall also show that $\frac{1}2 (\a_1+\a_2+\a_3+\a_4)\in M$ and hence $\frac{1}2(\a_1'+\a_2'+\a_3'+\a_4')\in M$.

Since $M/(M\cap U)\cong H_8$, there exist $\b_1, \b_2,\b_3,\b_4\in M\setminus (U\cap M)$ such that $\varphi(\bar{\b}_1), \varphi(\bar{\b}_2),\varphi(\bar{\b}_3),\varphi(\bar{\b}_4)$ generates the Hamming code $H_8$. By \refpp{glueaa4}, their projections to the last 4 coordinates are all even and thus spans an even subcode of $\ZZ_2^4$, which has dimension $\leq 3$. Therefore, there exists $a_1, a_2, a_3, a_4\in \{0,1\}$,
not all zero such that
$\varphi(a_1\bar{\b}_1+a_2\bar{ \b}_2+a_3\bar{\b}_3+a_4\bar{\b}_4)$ projects to zero and so must equal $(11110000)$.  Therefore, $\frac{1}2 (\a_1+\a_2+\a_3+\a_4)\in M$.

Since $|M/(U\cap M)|=2^4$, there exists $\b'=\frac{1}2(\a_i+\a_j+\a_j'+\a_k')$ and $\b''=\frac{1}2(\a_m+\a_n+\a_n'+\a_\ell')$ such that
\[
M= \mathrm{span}_\ZZ \left \{ M\cap U, \frac{1}2
(\a_1+\a_2+\a_3+\a_4), \frac{1}2 (\a_1'+\a_2'+\a_3'+\a_4'), \b', \b''\right\}.
\]

Note that
\[
\b'+\b''=\frac{1}2 ((\a_i+\a_j+\a_m+\a_n) + (\a_j'+\a_k'+\a_n'+\a_\ell')).
\]
Let $A:=(\{i,j\}\cup \{m,n\})- (\{i,j\}\cap \{m,n\})$ and $A':=(\{j,k\}\cup \{n,\ell\})-(\{j,k\}\cap \{n,\ell \})$. We shall show that $|\{i,j\}\cap \{m,n\}|=|\{j,k\}\cap \{n,\ell\}|=1$ and $|A\cap A'|=1$.

Since $\varphi(\bar{\b'}+\bar{\b''})\in H_8$ but  $\varphi(\bar{\b'}+\bar{ \b''})\notin \mspan_{\ZZ_2}\{ (11110000), (00001111)\}$, by \refpp{glueaa4}, we have
\[
\b'+\b''\in  \frac{1}2 (\a_p +\a_q+ \a_p'+\a_r') + M\cap U,
\]
for some $p,q\in \{ i,j,m,n\}$, $p, r\in \{j,k,n,\ell\}$ such that $p,q, r$ are distinct.

That means
$\frac{1}2(\a_i+\a_j+\a_m+\a_n) \in \frac{1}2(\a_p +\a_q) + M\cap U$ and
  $\frac{1}2(\a_j'+\a_k'+\a_n'+\a_\ell') \in \frac{1}2(\a_p'+\a_r')+ M\cap U$.
  It implies that
$A=(\{i,j\}\cup \{m,n\})- (\{i,j\}\cap \{m,n\})=\{p, q\}$ and $A'=(\{j,k\}\cup \{n,\ell\})-(\{j,k\}\cap \{n,\ell \})=\{p, r\}$. Hence, $|\{i,j\}\cap \{m,n\}|=|\{j,k\}\cap \{n,\ell\}|=1$ and $|A\cap A'|=1$.

By rearranging the indices if necessary, we may assume
$ \b'=\frac{1}2 (\a_1+\a_2+\a_2'+\a_4')$,
$ \b''=\frac{1}2 (\a_1+\a_3+\a_2'+\a_3')$ and hence
\[
M= \mathrm{span}_\ZZ \left \{ {M\cap U, \frac{1}2
(\a_1+\a_2+\a_3+\a_4), \frac{1}2 (\a_1'+\a_2'+\a_3'+\a_4'), \atop
\frac{1}2 (\a_1+\a_2+\a_2'+\a_4'),\frac{1}2
(\a_1+\a_3+\a_2'+\a_3')}\right\}
\]

Now let $\b_i=\a_i g$ and $\b_i'= \a_i'g$ for all $i=1,2,3,4$. Then
\[
N=Mg = \mathrm{span}_\ZZ \left \{ {N\cap U, \frac{1}2
(\b_1+\b_2+\b_3+\b_4), \frac{1}2 (\b_1'+\b_2'+\b_3'+\b_4'), \atop
\frac{1}2 (\b_1+\b_2+\b_2'+\b_4'),\frac{1}2
(\b_1+\b_3+\b_2'+\b_3')}\right\},
\]
as desired.  \eop

\medskip
Next we shall show that $L$ contains a sublattice isomorphic to a tensor product $A_4 \otimes A_4$.

\begin{nota}\labtt{aa4inM}
Take 
\begin{align*}
\g_0:=&\ \frac{1}2(-\a_1+\a_2+ \a_2'-\a_4'),\\
\g_1:=&\ \a_1,\\
\g_2:=&\ -\frac{1}2(\a_1+\a_2+\a_3+\a_4), \\
\g_3:=&\ \a_3,\\
\g_4:=&\ -\frac{1}2(\a_3-\a_4+ \a_2'-\a_4')
\end{align*}
in $M$ (cf. \refpp{dih10glue}) and set $R:=\mathrm{span}_\ZZ\{ \g_1,\g_2,\g_3,\g_4\}$. Then $R \cong
AA_4$. Note that $\g_0=-(\g_1+\g_2+\g_3+\g_4)$.
\end{nota}

\begin{lem}\labtt{innerproduct}
For any $i=0, 1,2,3, 4$ and $j=1,2,3,4$, we have
\begin{enumerate}
\item[(i)]  $(\g_i, \g_i g)= (\g_i, \g_i g^4)=-2$  and $(\g_i, \g_i g^2)= (\g_i, \g_i g^3)=0$;

\item[(ii)]  $ (\g_{j-1}, \g_j g)= (\g_{j-1}, \g_j g^4)=1$  and \\ $(\g_{j-1}, \g_j g^2)= (\g_{j-1}, \g_j g^3)=0$;

\item[(iii)] $(\g_i, \g_j g^k)=0$ for any $k$ if $|i-j|> 1$

\item[(iv)] $ \g_i\DD\cong AA_4 $.
\end{enumerate}
\end{lem}

 \pf Straightforward. \eop

\begin{prop}
Let $T=R\DD$. Then $T\cong A_4\otimes A_4$.
\end{prop}

\pf By (iv) of \refpp{innerproduct}, $ \g_i\DD=\mathrm{span}_\ZZ\{ \g_i g^j\ |\ j=0,1,2,3,4\}\cong AA_4 $.

Let $\{e_1,e_2,e_3,e_4\}$ be a fundamental basis of $A_4$ and denote $e_0= -(e_1+e_2+e_3+e_4)$.
Now define a linear map $\varphi: T\to A_4\otimes A_4$ by  $\varphi (\g_i g^j) = e_i\otimes e_j$, for $i, j = 1,2,3,4$.
By the inner product formulas in \refpp{innerproduct},
\[
(\g_i g^j, \g_k g^\ell)= (\g_i, \g_k g^{\ell-j})=
\begin{cases}
4 &\text{ if } i=k, j=\ell,\\
-2 &\text{ if } i=k, |j-\ell|=1,\\
1 &\text{ if } |i-k|=1, |j-\ell|=1,\\
0 &\text{ if } |i-k|>1, |j-\ell|>1.
\end{cases}
\]
Hence, $(\g_i g^j, \g_k g^\ell)= (e_i\otimes e_j,e_k\otimes e_\ell)$ for all $i,j,k,\ell$ and $\varphi$ is an isometry.
\eop

\appendix

\section{General results about lattices}

\begin{lem}\labtt{pl} Let $p$ be a prime number, $f(x):=1+x+x^2+\cdots + x^{p-1}$.  Let
$L$ be a $\ZZ [x]$-module.  For $v\in L$, $pv \in L(x-1)+Lf(x)$.
\end{lem}
\pf
We may write $f(x)=\sum_{i=0}^{p-1}x^i =  \sum_{i=0}^{p-1}((x-1)+1)^i=(x-1)h(x)+p$, for some $h(x)\in \ZZ [x]$.  Then if $v\in L$, $pv=v(f(x)-(x-1)h(x))$.
\eop

\begin{lem}\labtt{f/fmfnelemabel}
Suppose that the four group $D$ acts on the  abelian group $A$.
If the fixed point subgroup of $D$ on $A$ is 0, then
$A/Tel(A,D)$ is an elementary abelian 2-group.
\end{lem}
\pf
Let $a\in A$ and let $r\in D$.    We claim that $a(r+1)$ is an eigenvector
for $D$.  It is clearly an eigenvector for $r$.   Take $s\in D$ so that $D=\la r, s\ra$.
Then $a(r+1)s=a(r+1)(s+1) - a(r+1) = -a(r+1)$ since $a(r+1)(s+1)$ is a fixed point.
So, $a(r+1)$ is an eigenvector for $D$.

To prove the lemma, we just calculate that
$a(1+r)+a(1+s)+a(1+rs)=2a+a(1+sr+s+rs)=2a$ since $a(1+sr+s+rs)$ is a fixed point.
\eop

\begin{lem}\labtt{prankdg}
Suppose that $X$ is a lattice of rank $n$ and $Y$ is a sublattice of
rank $m$.  Let $p$ be a prime number. Suppose that the $p$-rank of
$\dg X$ is $r$.
Then,
$(Y\cap p\dual X)/(Y\cap pX)$ has $p$-rank at least $r+m-n$.  In particular,
the $p$-rank of $\dg Y$ is at least $r+m-n$; and if $r+m>n$, then $p$ divides $det(Y)$.
\end{lem}
\pf
We may assume that $Y$ is a direct summand of $X$.  The quadratic space
$X/pX$ has dimension $n$ over $\FF_p$ and its radical $p\dual X/pX$ has dimension $r$.  The image of $Y$ in $X/pX$  is $Y+pX/pX$, and it has dimension $m$ since $Y$ is a direct summand of $X$.
Let $q$ be the quotient map $X/pX$ to $(X/pX)/(p\dual X/pX)\cong X/p\dual X\cong p^r$.  Then $dim(q(Y+pX/pX))\le n-r$, so that $dim(Ker(q)\cap (Y+pX/pX))\ge m-(n-r)=r+m-n$.

We note that $Ker(q)=p\dual X/pX$, so  the above proves
$r+m-n \le rank((Y\cap p\dual X)+pX/pX)=rank(Y\cap p\dual X)/(Y\cap p\dual X \cap pX))=rank((Y\cap p\dual X)/(Y\cap pX))
=rank((Y\cap p\dual X)/pY)
$.
Note that $(Y\cap p\dual X)/pY\cong (\frac 1p Y \cap \dual X)/Y \le \dual Y/Y$, which implies the inequality of the lemma.
\eop

\begin{lem}\labtt{rescaling2}
Suppose that $Y$ is an integral lattice such that
there exists an integer $r>0$ so that
$\dg Y$ contains a direct product of $rank(Y)$ cyclic groups of order $r$.
Then $\frac 1{\sqrt r} Y$ is an integral lattice.
\end{lem}
\pf
Let $Y< X< \dual Y$ be a sublattice such that $X/Y\cong (\ZZ_r)^{rank\, Y}$. Then
$x\in X$ if and only if $rx\in Y$. Let $y, y'\in Y$. Then
$(\frac 1{\sqrt r}y, \frac 1{\sqrt r}y')=(\frac{1}ry, y')\in (X, Y )\leq (Y^*, Y)= \ZZ$.
\eop

\begin{lem}\labtt{dualrescaling}
Suppose that $X$ is an integral lattice and that there is an integer $s\ge 1$ so that
$\frac 1{\sqrt s}X$ is an integral lattice.
Then the subgroup $s \dg {X}$ is isomorphic to $\dg {\frac 1{\sqrt s}X}$ and
 $\dg X / s\dg X$ is isomorphic to $s^{rank(X)}$.
\end{lem}
\pf
Study the diagram below, in which the horizontal arrows are multiplication by $\frac 1{\sqrt s}$.  The hypothesis implies that the finite abelian group $\dg X$ is a direct sum of cyclic groups, each of which has order divisible by $s$.
$$\begin{matrix}
\dual X & \longrightarrow & \frac 1{\sqrt s}\dual X \cr
| & & | \cr
s\dual X & \longrightarrow & \sqrt s \dual X \cr
| & &| \cr
X & \longrightarrow & \frac 1{\sqrt s}X \cr
\end{matrix}.$$

We take  a vector $a\in \RR \otimes X$ and note that $a \in \dual{(\frac 1{\sqrt s}X)}$ if and only if $(a, \frac 1{\sqrt s}X)\le \ZZ$ if and only if $(\frac 1{\sqrt s}a,X)\le \ZZ$ if and only if $\frac 1{\sqrt s}a \in \dual X$ if and only if $a \in \sqrt s \dual X$.  This proves the first statement.  The second statement follows because
$\dg X$ is a direct sum of $rank(X)$ cyclic groups, each of which has order divisible by $s$.
\eop

\begin{lem}\labtt{involonlattice}
Let $y$ be an order 2 isometry of a lattice $X$.  Then $X/Tel(X,y)$
is an elementary abelian 2-group. Suppose that $X/Tel(X,y)\cong
2^c$. Then we have $det(X^+(y)) det(X^-(y)) = 2^{2c}det(X)$ and for $\vep =
\pm$, the image of $X$ in $\dg {X^\vep (y)}$ is $2^c$.  In
particular, for $\vep = \pm$, $det(X^\vep (y))$ divides $2^c det(X)$
and is divisible by $2^c$.  Finally, $c\le rank(X^\vep (y))$, for
$\vep = \pm$, so that $c\le \half rank(X)$.
\end{lem}
\pf See \cite{gre8}. \eop

\begin{lem}\labtt{lifteigvec}
Suppose that $t$ is an involution acting on the abelian group $X$.  Suppose that $Y$ is a $t$-invariant subgroup of odd index so that $t$ acts on $X/Y$ as a scalar $c\in \{ \pm 1 \}$.
Then for every coset $x+Y$ of $Y$ in $X$, there exists $u\in x+X$ so that
$ut=cu$.
\end{lem}
\pf  First, assume that $c=1$.
Define $n:=\half (|X/Y|+1)$, then take $u:=nx(t+1)$.
This is fixed by $t$ and $u\equiv 2nx \equiv x (mod\, Y)$.

If $c=-1$, apply the previous argument to the involution $-t$.
\eop

\begin{lem}\labtt{teldec}
If $X$ and $Y$ are abelian groups with $|X:Y|$ odd, an involution
$r$ acts on $X$,  and $Y$ is $r$-invariant, then $X/Tel(X,r)\cong Y/Tel(Y,r)$.
\end{lem}
\pf
Since $X/Y$ has odd order, it is the direct sum of its two eigenspaces for the action of $r$.  
Use \refpp{lifteigvec} to show that $Y+Tel(X,r)=X$ and $Y\cap Tel(X,r)= Tel(Y,r)$.  \eop

\begin{lem}\labtt{halfsquare} Suppose that $X$ is an integral lattice which has rank $m\ge1$
and there exists a lattice $W$, so that $X\le W \le \dual X$ and
$W/X\cong 2^r$, for some integer $r\ge 1$.  Suppose further that
every nontrivial coset of $X$ in $W$ contains a vector with
noninteger norm. Then $r=1$.
\end{lem}
\pf Note that if $u+X$ is a nontrivial coset of $X$ in $W$, then
$(u, u)\in \half +\ZZ$.

Let $\phi : X \rightarrow Y$ be an isometry of lattices, extended
linearly to a map between duals.  Let $Z$ be the lattice between
$X\perp Y$ and $W \perp \phi (W)$ which is diagonal with respect to
$\phi$, i.e., is generated by  $X\perp Y$ and  all vectors of the
form $(x, x \phi )$, for $x\in W$.

Then $Z$ is an integral lattice.  In any integral lattice, the even sublattice has index 1 or 2.  Therefore, $r=1$ since the nontrivial cosets of $X\perp Y$ in $Z$ are odd.
\eop

\begin{lem}\labtt{mnfsummands}
Suppose that the  integral lattice $L$ has no vectors of norm  $2$ and that $L=M+N$, where $M\cong N\cong EE_8$. The sublattices $M, N, F=M\cap N$ are direct summands of $L=M+N$.
\end{lem}
\pf
Note that $L$ is the sum of even lattices, so is even.  Therefore, it has no vectors of norm 1 or 2.
Since $M$ by definition defines the  summand $S$ of negated vectors by
$t_M$, we get $S=M$ because $M\le S \le \half M$ and the minimum norm
of $L$ is 4.  A  similar statement holds for $N$.  The sublattice $F$
is therefore the sublattice of vectors fixed by both $t_M$ and $t_N$,
so it is clearly a direct summand of $L$.
\eop

\begin{lem} \labtt{commutatoreigen}
$D\cong Dih_6$, $\la g\ra=O_3(D)$ and $t$ is an involution in $D$.  Suppose that $D$ acts on the abelian group $A$, $3A=0$ and
$A(g-1)^2=0$.  Let $\vep = \pm 1$.
If $v\in A$ and $vt=\vep v$, then $v(g-1)t=-\vep v(g-1)$.
\end{lem}
\pf  Calculate $v(g-1)t=vt(g^{-1}-1)=\vep v (g^{-1}-1)$.  Since $(g-1)^2=0$, $g$   acts as 1 on the image of $g-1$, which is the image of $(g^{-1}-1)$.   Therefore,
$\vep v (g^{-1}-1)=\vep v (g^{-1}-1)g=\vep v(1-g)=-\vep v(g-1)$.
\eop

\begin{lem}\labtt{subquotientdg}
Suppose that $X$ is an integral lattice and $Y$ has finite index, $m$,  in $X$.   Then $\dg X$ is a subquotient of $\dg Y$ and
$|\dg X|m^2=|\dg Y|$.  The groups have isomorphic Sylow $p$-subgroups if $p$ is a prime which does not divide $m$.
\end{lem}
\pf
Straightforward.  \eop

\begin{lem}\labtt{dgxdgy}
Suppose that $X$ is a lattice, that $Y$ is a direct summand and $Z:=ann_X(Y)$.
Let $n:=|X:Y\perp Z|$.

(i)
The image of $X$ in $\dg Y$ has index dividing $(det(Y), det(X))$.  In particular, if $(det(Y), det(X))=1$, $X$ maps onto $\dg Y$

(ii)
Let $A:=ann_{\dual X}(Y)$.   Then $\dual X/(Y\perp A)\cong \dg Y$.

(iii)
There are epimorphisms of groups $\varphi_1 : \dg X \rightarrow \dual X/(X+A)$ and
$\varphi_2 : \dg Y \rightarrow \dual X/(X+A)$.

(iv)
We have isomorphisms $Ker (\varphi_1 ) \cong (X+A)/X$ and
 $Ker(\varphi_2) \cong  \psi(X)/Y \cong X/(Y\perp Z)$.  The latter is a group of order $n$.

In particular, $Im(\varphi_1 )\cong Im(\varphi_2)$ has order $\frac 1n det(Y)$.

(v)
If $p$ is a prime which does not divide $n$, then $O_p(\dg Y )$ injects into $O_p(\dg X)$.  This injection is an isomorphism onto if $(p, det(Z))=1$.
\end{lem}
\pf
(i) This is clear since the natural map $\psi : \dual X \rightarrow \dual Y$ is onto and $X$ has index $det(X)$ in $X^*$.

(ii) The natural map $\dual X \rightarrow \dual Y$ followed by the quotient map $\zeta : \dual Y \rightarrow \dual Y/Y$ has kernel $Y\perp A$.

(iii) Since $\dg X =\dual X/X$, we have the first epimorphism.   Since $X+A\ge Y+A$, existence of the second epimorphism follows from (ii).

(iv) First,   $Ker (\zeta \psi )=Y\perp A$ follows from (ii)
and the definitions  of $\psi$ and $\zeta$.
So, $Ker(\varphi_1 )\cong (X+A)/(Y+A)$.
Note that
$(X+A)/(Y+A)=(X+(Y+A))/(Y+A)\cong X/((Y+A)\cap X) = X/(Y+Z))$.  The latter quotient has order $n$. For the order statement, we use the formula $|Im(\psi )||Ker(\psi )|=|\dg Y|$.
For the second isomorphism, use $\dual Y \cong \dual X/A$
and
$\dg Y=\dual Y/Y \cong (\dual X/A)/((Y+A)/A)$ and note that in here the image of $X$ is $((X+A)/A)/((Y+A)/A) \cong (X+A)/(Y+A)$.

(v) Let $P$ be a Sylow
$p$-subgroup of $\dg Y$.
Then $P \cap Ker(\psi )=0$ since $(p,n)=1$.  Therefore $P$ injects into $Im(\psi )$.
The epimorphism $\varphi$ has kernel $A/Z=\dg Z$.  So, $P$ is isomorphic to a Sylow $p$-subgroup of $\dg X$ if $(p, det(Z))=1$.
\eop

\begin{lem} \labtt{elemabelsylow}
Suppose that $X$ is an integral lattice
and $E$ is an elementary abelian 2-group acting in $X$.
If $H$ is an orthogonal direct summand of $Tel(X,E)$ and $H$ is a direct summand of $X$, then
the odd order Sylow groups of $\dg H$ embed in $\dg X$.  In other notation, $O_{2'}(\dg H)$ embeds in $O_{2'}(\dg X)$
\end{lem}
\pf
Apply \refpp{dgxdgy} to $Y=H$, $n$ a divisor of $|X:Tel(X,E)|$, which is a power of 2.
\eop

\begin{lem}\labtt{dualxdualy}
Suppose that $X$ is a lattice and $Y$ is a direct summand  and $Z:=ann_X(Y)$.
Let $n:=|X:Y\perp Z|$.

Suppose that $v\in \dual X$ and that $v$ has order $m$ modulo $X$.
If $(m,det(Z))=1$, there exists $w\in Z$ so that $v-w\in \dual Y \cap \dual X$.   Therefore, the coset $v+X$ contains a representative in $\dual Y$.
Furthermore, any element of $\dual Y \cap (v+X)$ has order $m$ modulo $Y$.
\end{lem}
\pf
The image of $v$ in $\dg Z$ is zero, so
the restriction of the function $v$ to $Z$ is the same as taking a dot product with an element of $Z$.
In other words,  the projection of $v$ to $Z ^ *$ is already in $Z$. Thus, there exists $w\in Z$ so that $v-w \in ann_{\dual X}(Z) =  \dual Y$.

Define $u:=v-w$.   Then $mu =mv -mw \in X$.   Since $v-w \in \dual Y$, $mu \in X \cap \dual Y$, which is $Y$ since $Y$ is a direct summand of $X$.

Now consider an arbitrary $u \in (v+X) \cap \dual Y$.  We claim that its order modulo $Y$ is $m$.
There is $x \in X$ so that $u=x+v$.  Suppose $k>0$.  Then $ku=kx+kv$ is in $X$ if and only if $kv \in X$, i.e. if and only if $m$ divides $k$.
\eop

\begin{lem}\labtt{telx}
Let $D$ be a dihedral group of order $2n$, $n>2$ odd, and $Y$ a finitely generated free abelian group which is a $\ZZ [D]$-module, so that an element $1\ne g\in D$ of odd order acts with zero fixed point subgroup on $Y$.
Let $r$ be an involution of $D$ outside $Z(D)$.
Then
$Y/Tel(Y,r)$ is elementary abelian of order $2^{\half rank(Y)}$.
Consequently, $det(Tel(Y,r))=2^{rank(Y)}det(Y)$.
\end{lem}
\pf
The first statement follows since the odd order $g$ is inverted by $r$ and acts without fixed points on $Y$.  The second statement follows from the index formula for determinants.
\eop

\begin{lem}\labttr{sqrt3e8ine8}
Suppose that $X\cong E_8$ and that $Y$ is a sublattice such that $X/Y\cong 3^4$ and $Y\cong \sqrt 3 E_8$.  Then there exists an element $g$ of order 3 in $O(X)$ so that $X(g-1)=Y$.  In particular, $Y$ defines a partition on the set of 240 roots in $X$ where two roots are equivalent if and only if their difference lies in $Y$.
\end{lem}
\pf
Note that $Y$ has 80 nontrivial cosets in $X$ and the set $\Phi$ of roots has cardinality 240.
Let $x, y \in \Phi$ such that $x-y \in Y$ and $x, y$ are linearly independent.
Then $6\le (x-y,x-y)=4-2(x,y)$, whence $|(x, y)|\leq 2$.
Therefore $(x,y)=-1$, since $x, y$ are linearly independent roots.
Let $z$ be a third root which is congruent to $x$ and $y$ modulo $Y$.  Then
$(x,z)=-1=(y,z)$ by the preceding discussion.   Therefore the projection of $z$ to the span of $x, y$ must be $-x-y$, which is a root.  Therefore, $x+y+z=0$.

It follows that a nontrivial coset of $Y$ in $X$ contains at most three roots.  By counting, a nontrivial coset of $Y$ in $X$ contains exactly three roots.

Let $P:=span\{x, y\} \cong A_2$, $Q:=ann_X(P)$.   We claim that $Y\le P\perp Q$, which has index 3 in $X$.   Suppose the claim is not true.  Then the structure of $\dual P$ means that there exists $r\in Y$ so that $(r,x-y)$ is not divisible by 3.  Since $x-y\in Y$, we have a contradiction to
 $Y\cong \sqrt 3 E_8$.
 The claim implies that $g_P$, an automorphism of order 3 on $P$, extended to
$X$ by trivial action on $Q$, leaves $Y$ invariant (since it leaves invariant
any sublattice between  $P\perp Q$ and $\dual P \perp \dual Q$).  Moreover,
$X(g_P-1)=\ZZ (x-y)$.

Now take a root $x'$ which is in $Q$.
Let $y', z'$ be the other members of its equivalence class of $x'$.
We claim that these are also in $Q$.
We know that $x'-y'\in Y \le P\perp Q$, so
 $P':=span\{ x', y', z'\} \le P\perp Q$.
 Now, we claim that the projection of $P'$ to $P$ is 0.
 Suppose otherwise.   Then the projection of some root $u \in P'$ to $P$ is nonzero.   Therefore the projection is a root, i.e. $u \in P$.  But then $u$ is in the equivalence class of $x$ or $-x$ and so $P'=P$, a contradiction to $(x',P)=0$.

 We now have that the class of $x'$ spans a copy of $A_2$ in $Q$.  We may continue this procedure to get a sublattice $U=U_1\perp U_2\perp U_3 \perp U_4$ of $X$ such that $U_i \cong A_2$ for  of $X$ with the property that if
$g_i$ is an automorphism of order 3 on $U_i$ extended to $X$ by trivial action on $ann_X(U_i)$, then
each $U(g_i-1)\le Y$ and, by determinants, $g:=g_1g_2g_3g_4$ satisfies
$U(g-1)=Y$.
\eop

\begin{prop}\labtt{overa2timese8}
Suppose that $T\cong A_2\otimes E_8$.

(i)
Then
$\dg T\cong 3^8$ and the natural $\third \ZZ /\ZZ$-valued
quadratic form has maximal Witt index; in fact, there is a natural identification of quadratic spaces $\dg T$ with $E_8/3 E_8$, up to scaling.

(ii)
Define
$\mathcal O_k :=\{ X | T \le X \le \dual T, X \hbox{ is an integral lattice}, dim(X/T)=k \}$ (dimension here means over $\FF_3$).

(a) $\mathcal O_k$ is nonempty if and only if $0\le k \le 4$;

(b) $\mathcal O_k$ consists of even lattices
for each  $k$, $0\le k \le 4$;

(c) On $\dual T/T$, the action of $g$, the isometry of order 3 on $T$ corresponding to an order 3 symmetry of the $A_2$ tensor factor, is trivial.  Therefore any lattice between $T$ and $\dual T$ is $g$-invariant.

\smallskip

\smallskip

(iii)  the lattices in $\mathcal O_k $ embed in $E_8 \perp E_8$.  For a fixed $k$, the embeddings are unique up to the action of $O(E_8 \perp E_8 )$.

(iv)  the lattices in $\mathcal O_k $ are rootless if and only if $k=0$.
\end{prop}
\pf
(i) This follows since the quotient $\dual T / T $ is covered by $\third P$, where $P\cong \sqrt 6 E_8$ is $ann_T(E)$, where $E$ is one of the three $EE_8$-sublattices of $T$.

(ii) Observe that $X\in \mathcal O_k$ if and only if $X/T$ is a totally
singular subspace of $\dg T$.  This implies (a).  An integral lattice is even if
it contains an even sublattice of odd index.  This implies (b). For (c), note that $T^*/T$ is covered by $\frac{1}3P$ and $P(g-1)\leq T(g-1)^2 = 3T$. Hence $g$ acts trivially on $\dg T=T^*/T$. This means any lattice $Y$ such that $T\le Y \le T^*$ is $g$-invariant.

(iii) and (iv)   First, we take $X\in \mathcal O_4$ and prove $X\cong E_8 \perp E_8$.   Such an $X$ is even, unimodular and has rank 16, so is isometric to $HS_{16}$ or $E_8 \perp E_8$.  Since $X$ has a fixed point free automorphism of order 3, $X\cong E_8 \perp E_8$.  Such an automorphism fixes both direct summands.   Call these summands $X_1$ and $X_2$.  Define $Y_i:=X_i(g-1)$, for $i=1,2$.  Thus, $Y_i\cong \sqrt 3 E_8$.

The action of $g$ on $X_i\cong E_8$ is unique up to conjugacy, namely as a diagonally embedded cyclic group of order 3 in a natural $O(A_2)^4$ subgroup of $O(X_i)\cong Weyl(E_8)$ (this follows from the corresponding conjugacy result for $O^+(8,2)\cong Weyl(E_8)/Z(Weyl(E_8))$.

We consider how $T$ embeds in $X$.   Since $|X:T|=3^4$, $|X_i : T\cap X_i|$ divides $3^4$.
Since $T\cap X_i \ge X_i (g-1)$ and $T$ has no roots, \refpp{sqrt3e8ine8} implies that $T\cap X_i =  Y_i$, for $i=1,2$.

If $U\in \mathcal O_k$, $T\le U \le X$, then rootlessness of $U$ implies that
$U\cap X_i = Y_i$ for $i=1,2$.   Therefore, $U=T$, i.e. $k=0$.
\eop

\begin{lem}\labtt{scalarisometry}  Let $q$ be an odd prime power and let $(V,Q)$ be a finite dimensional quadratic space over $\FF_q$.
Let $c$ be a generator of $\FF_q^\times$.

If $dim(V)$ is odd, there exists $g \in GL(V)$ so that $gQ=c^2Q$.

If $dim(V)$ is even, there exists $g \in GL(V)$ so that $gQ=cQ$.
\end{lem}
\pf
The scalar transformation $c$ takes $Q$ to $c^2Q$.  This proves the result in case $dim(Q)$ is odd.
Now suppose that $dim(V)$ is even.

Suppose that $V$ has maximal Witt index.  Let $V=U\oplus U'$, where $U, U'$ are each totally singular.  We take $g$ to be $c$ on $U$ and 1 on $U'$.

Suppose that $V$ has nonmaximal Witt index.  The previous paragraph allows us to reduce the proof to the case $dim(V)=2$ with $V$ anisotropic. (One could also observe that if we write $V$ as the orthogonal direct sum of nonsingular 2-spaces, the result follows from the case $dim(V)=2$.) Then $V$ may be identified with $\FF_{q^2}$ and $Q$ with a scalar multiple of the norm map.  We then take $g$ to be multiplication by a scalar $b\in \FF_{q^2}$ such that $b^{q+1}=c$.
\eop

\section{Characterizations of lattices of small rank}

Some results in this section are in the literature.  We collect them here for convenience.

\begin{lem}\labtt{rank2}
Let $J$ be a rank 2 integral lattice.
If $det(J)\in \{1,2,3,4,5,6\}$, then $J$ contains a vector of norm 1 or 2.   If $det(J)\in \{1, 2\}$, $J$ is rectangular.   If $J$ is even, $J \cong A_1\perp A_1$ or $A_2$.
\end{lem}
\pf
The first two statements follows from values of the Hermite function.  Suppose that $J$ is even.  Then $J$ has a root, say $u$.  Then $ann_J(u)$ has determinant $\frac{1}2 det(J)$ or $2\, det(J)$.  If $ann_J(u)=\frac{1}2 det(J)$, then $J$ is an orthogonal direct sum $\ZZ u \perp \ZZ v$, for some vector $v\in J$.
For $det(J)$ to be at most 6 and $J$ to be even, $(v,v)=2$ and $J$ is the lattice $A_1\perp A_1$.
Now assume that
$ann_J(u)$ has determinant $2\, det(J)$, an integer at most 12.  Let $v$ be a basis for $ann_J(u)$.  Then $\half (u+v)\in J$ and so $\fourth (2+(v,v))\in 2\ZZ$ since $J$ is assumed to be even.  Therefore, $2+(v,v)\in 8\ZZ$.  Since $(v,v)\le 12$, $(v,v)=6$.  Therefore, $\half (u+v)$ is a root and we get $J \cong A_2$.
\eop

\begin{lem}\labtt{rank3}\label{rank 3}
Let $J$ be a rank 3 integral lattice.
If $det(J)\in \{1,2,3\}$, then $J$ is rectangular or $J$ is isometric to $\ZZ \perp A_2$.
If $det(J)=4$, $J$ is rectangular or is isometric to $A_3$. 
\end{lem}
\pf
If $J$ contains a unit vector, $J$ is orthogonally decomposable and we are done by \refpp{rank2}.
Now use the Hermite function: $H(3,2)=1.67989473\dots$, $H(3,3)=1.92299942\dots$ and $H(3,4)=2.11653473\dots $.  We therefore get an orthogonal decomposition unless possibly
$det(J)=4$ and $J$ contains no unit vector.  Assume that this is so.

If $\dg J$ is cyclic, the lattice $K=J+2J^*$ which is strictly between $J$ and $\dual J$ is integral and unimodular, so is isomorphic to $\ZZ^3$.
So, $J$ has index 2 in  $\ZZ^3$, and the result is easy to check.  If $\dg J\cong 2 \times 2$, we are done by a similar argument provided a nontrivial coset of $J$ in $\dual J$ contains a vector of integral norm.  If this fails to happen, we quote \refpp{halfsquare} to get a contradiction.
\eop

\begin{lem}\labtt{rank4det4}
Suppose that $X$ is an integral lattice which has rank 4 and determinant 4.
Then $X$ embeds with index 2 in $\ZZ^4$.
If $X$ is odd, $X$ is isometric to one of $2\ZZ \perp \ZZ^3, A_1\perp
A_1 \perp \ZZ^2, A_3\perp \ZZ$. If $X$ is even, $X\cong D_4$.
\end{lem}
\pf  Clearly, if $X$ embeds with index 2 in $\ZZ^4$, $X$ may be
thought of as the annihilator mod 2 of a vector $w\in \ZZ$ of the
form $(1,\dots, 1, 0, \dots , 0)$.  The isometry types for $X$
correspond to the cases where the weight of $w$ is 1, 2, 3 and 4. It
therefore suffices to demonstrate such an embedding.

First, assume that $\dg X$ is cyclic.   Then $X+2\dual X$ is an
integral lattice (since $(2x,2y)=(4x,y)$, for $x, y \in \dual X$)
and is unimodular, since it contains $X$ with index 2.  Then the classification of unimodular integral lattices of small rank implies $X+2\dual X\cong \ZZ^4$, and the conclusion is
clear.

Now, assume that $\dg X$ is elementary abelian. By
\refpp{halfsquare}, there is a nontrivial coset $u+X$ of $X$ in
$\dual X$ for which $(u,u)$ is an integer.  Therefore, the lattice
$X':=X+\ZZ u$ is integral and unimodular. By the classification of unimodular integral lattices,
$X'\cong \ZZ ^4$. \eop

\medskip

\begin{thm}\labttr{rank8} Let $L$ be a unimodular integral lattice of rank at most 8.  Then $L\cong \ZZ^n$ or $L\cong E_8$.
\end{thm}
\pf
This is a well-known classification.   The article \cite{gre8} has an elementary proof   
and discusses the history.  \eop

The next result is well known. The proof may be new.

\begin{prop}\labtt{uniquenesse6} Let $X$ be an integral lattice of
determinant 3 and rank at most 6.  Then $X$ is rectangular; or $X\cong
A_2 \perp \ZZ^m$, for some $m\le 4$; or $X\cong E_6$.
\end{prop}
\pf
Let $u\in \dual X \setminus X$.  Since $3u\in X$, $(u,u)\in \third \ZZ$.  Since $det(\dual X)=\third$, $(u, u)\in \third + \ZZ$ or $\tthird + \ZZ$.

Suppose $(u, u)\in \third + \ZZ$. Let $T\cong A_2$.    Then we quote  \refpp{duala2} to see that there is a unimodular lattice, $U$, which contains $X\perp T$ with index 3.

Suppose $U$ is even.  By the classification \refpp{rank8}, $U\cong E_8$.
A well-known property of $E_8$ is that all $A_2$-sublattices are in one orbit under the Weyl group.   Therefore, $X\cong E_6$.

If $U$ is not even, $U\cong \ZZ^n$, for some $n\le 8$. Any root in $\ZZ^n$ has the form $(\pm 1, \pm 1, 0,0,\dots , 0, 0)$.  It follows that every $A_2$ sublattice of $\ZZ^n$ is in one orbit under the isometry group $2\wr Sym_n$.  Therefore $X=ann_U(T)$ is rectangular.

Suppose $(u, u)\in \tthird + \ZZ$.  Then we consider a unimodular lattice $W$ which contains $X\perp \ZZ v$ with index 3, where $(v,v)=3$.  By the classification, $W\cong \ZZ^7$.  Any norm 3 vector in $\ZZ^n$ has the form $(\pm 1,\pm 1, \pm 1,0,0,\dots ,0)$ (up to coordinate permutation).  Therefore, $ann_W(v)$ must be isometric to $\ZZ^4 \perp A_2$.
\eop

\begin{lem}\labtt{rank4det5}  If $M$ is an even integral lattice of determinant
5 and rank 4, then $M\cong A_4$.
\end{lem}
\pf
Let $u \in \dual M$ so that $u +M$ generates
$\dg M$.
Then $(u,u)=\frac k5$, where $k$ is an integer.  Since $5u \in M$, $k$ is an
even integer.
Since $H(4,\frac 15 )=1.029593054\dots $,
a minimum norm vector in $\dual M$ does not lie in $M$, since $M$ is an even
lattice.   We may assume that $u$ achieves this minimum norm.  Thus,
$k\in \{2, 4\}$.

Suppose that $k=4$.
Then we may form $M \perp \ZZ 5v$, where $(v,v)=\frac 15$.
Define $w:=u+v$.   Thus, $P:=M+\ZZ w$ is a unimodular integral lattice.  By the
classification, $P\cong \ZZ^5$, so we identify $P$ with $\ZZ^5$.  Then
$M=ann_P(y)$ for some norm 5 vector $y$.   The only possibilities for such $y\in
P$ are $(2,1,0,0,0)$, $(1,1,1,1,1)$, up to monomial transformations.   Since $M$
is even, the latter possibility must hold and we get $M\cong A_4$.

Suppose that $k=2$.
We let $Q$ be the rank 2 lattice with Gram matrix $\begin{pmatrix}3&2\cr 2&3
\end{pmatrix}$.  So, $det(Q)=5$ and there is a generator $v\in \dual Q$ for
$\dual Q$ modulo $Q$ which has norm $\frac 35$.
We then form $M\perp Q$ and define $w:=u+w$.  Then $P:=M+Q+\ZZ w$ is an integral
lattice of rank 6 and determinant 1.  By the classification, $P\cong \ZZ^6$.  In
$P$, $M$ is the annihilator of a pair of norm 3 vector, say $y$ and $z$.  Each
corresponds in $\ZZ^6$ to some vector of shape
$(1,1,1,0,0,0)$, up to monomial transformation.    Since $M$ is even, the
6-tuples representing $y$ and $z$ must have supports which are disjoint 3-sets.
However, since $(y,z)=2$, by the Gram matrix, we have a contradiction. \eop

\begin{nota}\labttr{m(4,25)} We denote by $\mathcal{M}(4,25)$ an even integral lattice of rank 4 and determinant 25. By \refpp{rank4det25}, it is unique.
\end{nota}

\begin{lem}\labtt{rank4det25}
(i)
There exists a unique, up to isometry,  rank 4 even integral lattice whose
discriminant group has order 25.

(ii)
It is isometric to a glueing of the orthogonal direct sum $A_2\perp \sqrt 5 A_2$ by a glue vector of the shape $u+v$, where $u$ is in the dual of the first summand and $(u,u)=\frac 23$, and where $v$ is in the dual of the second summand and has norm $\frac {10}3$.

(iii)  The set of roots forms a system of type $A_2$; in particular, the lattice does not contain a pair of orthogonal roots.

(iv) The isometry group is isomorphic to $Sym_3\times Sym_3 \times 2$, where the first factor acts as the Weyl group on the first summand in (iii) and trivially on the second, the second factor acts as the Weyl group on the second summand of (iii) and trivially on the first, and where the third direct factor acts as $-1$ on the lattice.

(v) The isometry group  acts transitively on (a)  the six roots; (b) the 18 norm 4 vectors; (c) ordered pairs of norm 2 and norm 4 vectors which are orthogonal; (d) length 4 sequences of orthogonal vectors whose norms are 2, 4, 10, 20.

(vi)
An  orthogonal  direct sum of two such embeds as a sublattice of index $5^2$ in $E_8$.
\end{lem}
\pf
The construction of (ii) shows that such a lattice exists and it is easy to deduce (iii), (iv),  and (v).

We now prove (i).
Suppose that $L$ is such a lattice.
We observe that if the discriminant group were cyclic of order 25, the unique lattice strictly between $L$ and its dual would be even and unimodular.  Since $L$ has rank 4, this is impossible.   Therefore, the
discriminant group has shape $5^2$.

Since
$H(4,25)=3.44265186\dots$, $L$ contains a root, say $u$.  Define
$N:=ann_L(u)$.
Since $H(3,50)=4.91204199\dots ,$ $N$ contains a norm 2 or 4
element, say $v$.

Define $R:=\ZZ u \perp \ZZ v$ and  $P:=ann_L(R)$, a sublattice of
rank $2$ and determinant $2(v,v)\cdot 25$. Also, the Sylow 5-group
of $\dg P$ has exponent 5. Then $P\cong \sqrt 5 J$, where $J$ is an
even, integral lattice of rank 2 and $det(J)=2(v,v)$.  Since $det(J)$ is
even and the rank of the natural bilinear form on $J/2J$ is even, it follows
that $J=\sqrt{2}K$, for
an integral, positive definite lattice $K$.  We have $det(K)=\half
(v,v)\in \{1, 2\}$  and so $K$ is rectangular.  Also, $P\cong \sqrt {10}K$.
 So, $P$ has rectangular basis $w, x$ whose
norm sequence is $10, 5(v,v)$

Suppose that $(v,v)=2$.  Also, $L/(R\perp P) \cong 2\times 2$.  A
nontrivial coset of $R\perp P$ contains an element of the form
$\half y + \half z$, where $y \in span\{u, v\}$ and $z\in span\{w, x\}$. We
may furthermore arrange for $y=au+bv$, $z=cw+dx$, where $a, b, c, d
\in \{0,1\}$.   For the norm of $\half y +
\half z$ to be an even integer, we need $a=b=c=d=0$ or $a=b=c=d=1$.
This is incompatible with  $L/(R\perp P) \cong 2\times 2$.
Therefore, $(v,v)=4$.

We have  $L/(R\perp P)\cong 5^2$.
Therefore, $\frac 15 w$ and $\frac 15 x$ are in $\dual L$ but are not in $L$

Form the orthogonal sum $L\perp \ZZ y$, where $(y,y)=5$.  Define
$w:=\frac 15 x +\frac 15 y$.  Then $(w,w)=1$.  Also, $Q:=L+\ZZ w$ has rank 5, is integral and  contains
$L\perp \ZZ y$ with index 5, so has determinant 5.   Since $w$ is a unit vector,
$S:=ann_Q(w)$ has rank 4 and determinant 5, so $S\cong A_4$.
Therefore, $Q=\ZZ w \perp S$
and $L=ann_Q(y)$ for some $y\in Q$ of norm 5, where $y=e+f$, $e\in \ZZ w$, $f\in S$.
Since $S$ has no vectors of odd norm, $e\ne 0$ has odd norm.  Since $(y,y)=5$ and since $(e,e)$ is a perfect square, $(e,e)=1$ and $(f,f)=4$.   Since $O(A_4)$ acts transitively on norm $4$ vectors of $A_4$,  $f$ is uniquely determined up to the action of $O(S)$.   Therefore, the isometry type of $L$ is uniquely determined.

It remains to prove (vii).
 For one proof, use \refpp{1mod4}. Here is a second proof.  We may form an orthogonal direct sum of two such lattices and extend upwards by certain glue vectors.

Let $M_1$ and $M_2$ be two mutually orthogonal copies of $L$. Let $u, v, w,x$ be the orthogonal elements of $M_1$ of norm $2,4,10,20$ as defined in the proof of (i). Let $u', v', w', x'$ be the corresponding elements in $M_2$.  Set
\[
\g =\frac{1}5(w+x+x')\quad \text{ and } \quad \g'= \frac{1}5(x+w'+x').
\]
Their norm are both $2$. By computing the Gram matrix, it is easy to show that
$ E=\mspan_\ZZ\{ M_1, M_2, \g, \g'\}$ is integral and has determinant $1$. Thus, $E$ is even and so $E\cong E_8$.
\eop

\begin{lem}\labtt{1mod4}
 Let $p$ be a prime which is 1(mod 4).  Suppose that $M,  M'$ are lattices such that $\dg M$ and $\dg {M'}$  are elementary abelian $p$-groups which are isometric as quadratic spaces over $\FF_p$.  Let $\psi$ be such an isometry and let $c\in \FF_p$ be a square root of $-1$.    Then the overlattice $N$
of $M\perp M'$ spanned by the ``diagonal  cosets" $\{ \a + c \a \psi | \a \in \dg M \}$ is unimodular.  Also, $N$ is even if  $M$ and $M'$ are even.
\end{lem}
\pf The hypotheses imply that $N$ contains $M+M'$ with index
$|det(M)|$, so is unimodular. It is integral since the space of
diagonal cosets so indicated forms a maximal totally singular
subspace of the quadratic space $\dg M \perp \dg {M'}$.  The last
sentence follows since $|N:M\perp M'|$ is odd. \eop

\begin{lem}\labtt{rank4det125}  An even rank 4 lattice with discriminant group which is elementary abelian of order 125 is isometric to $\sqrt 5 \dual A_4$.
\end{lem}
\pf
Suppose that $L$ is such a lattice. Then  $\det(\sqrt{5}\dual L )=5$.
We may apply the result \refpp{rank4det5} to get  $\sqrt{5}\dual L \cong A_4$.
\eop

\begin{lem}\labttr{rank4det9}
An even integral lattice of rank 4 and determinant 9 is isometric to $A_2^2$.
\end{lem}
\pf
Let $M$ be such a lattice.     Since $H(4,9)=
2.66666666\dots $ and $H(3,18)=3.494321858\dots $, $M$ contains an
orthogonal pair of roots, $u, v$. Define $P:=\ZZ u \perp \ZZ v$.
The natural map $M\rightarrow   \dg P$ is onto since $(\det M, \det P)=1$.
Therefore, $Q:=ann_M(P)$ has determinant 36 and the image of $M$ in $\dg Q$ is
$2\times 2$.
Therefore, $Q$ has a rectangular basis $w, x$, each of norm 6 or with respective norms 2, 18.

We prove that 2, 18 does not occur.
Suppose that it does.  Then there is a sublattice $N$ isometric to $A_1^3$.  Since there are no even integer norm vectors in $\dual N \setminus N$, $N$ is a direct summand of $M$.  By coprimeness, the natural map of $M$ to $\dg N\cong 2^3$ is onto.   Then the natural map of $M$ to $\dg {\ZZ x}$ has image isomorphic to $2^3$.   Since $\dg {\ZZ x}$ is cyclic, we have a contradiction.

Since $M/(P\perp Q)\cong 2\times 2$ and $M$ is even, it is easy to see that $M$
is one of $M_1:=span\{P, Q, \half (u+w), \half (v+x)\}$ or
$M_2:=span\{P, Q, \half (u+x), \half (v+w)\}$.   These two
overlattices are isometric by the isometry defined by $u\mapsto u, v
\mapsto v, w\mapsto x, x \mapsto w$.
It is easy to see directly that they are isometric to $A_2^2$.  For
example, $M_1 =span\{u, x, \half (u+x)\} \perp  span\{v, w, \half
(v+w) \}$.
\eop

\section{Nonexistence of particular lattices}

\begin{lem}\labtt{no3times3}
Let $X\cong \ZZ^2$.
There is no sublattice of $X$ whose discriminant group is $3\times 3$.
\end{lem}
\pf Let $Y$ be such a sublattice.  Its index is 3.
Let $e, f$ be an orthonormal basis
of $X$. Then $Y$ contains $W:=span\{3e, 3f\}$ with index 3. Let
$v\in Y \setminus W$, so that $Y=W+\ZZ v$.  If $v$ is $e$ or $f$,
clearly $\dg Y$ is cyclic of order 9. We may therefore assume that
$v=e+f$ or $e-f$. Then $Y$ is spanned by $3e$ and $e\pm f$, and so
its Smith invariant sequence  1, 9. This final contradiction
completes the proof. \eop

\begin{lem}\labtt{norank4det3}  There does not exist an even rank 4 lattice of determinant 3.
\end{lem}
\pf
Let $L$ be such a lattice and let $u \in \dual L$ so that $u$ generates $\dual L$ modulo $L$.
Then $(u, u)=\frac k3$ for some integer $k>0$.  Since $L$ is even, $k$ is even.
Since $H(4,\frac 13)=1.169843567\dots <4/3$, we may assume that $k=2$.

We now form $L \perp \ZZ (3v)$, where $(v,v)=\third$.   Define $w:=u+v$.
The lattice $P:=L+\ZZ w$ is unimodular, so is isometric to $\ZZ^5$.  Since $det(M)=3$, $M=ann_P(y)$ for some vector $y$ of norm 3.   This forces $M$ to be isometric to $A_2 \perp \ZZ^2$, a contradiction to evenness.
\eop

\begin{coro} \labtt{norank4det27} There does not exist an even rank 4 lattice whose discriminant group is elementary abelian of order $3^3$.
\end{coro}
\pf
If $M$ is such a lattice, then $3\dual M$ has rank 4 and determinant 3.  Now use \refpp{norank4det3}.
\eop

\section{Properties of particular lattices}

\begin{lem}\labtt{smallrankssde8}
Suppose that $M\ne 0$ is an SSD sublattice of $E_8$ and that $rank(M)\le 4$.
Then $M$ contains a root, $r$, and
$ann_M(r)$ is an SSD sublattice of $E_8$ of rank $(rank(M)-1)$.
Also $M\cong A_1^{rank(M)}$ or $D_4$.
\end{lem}
\pf  Let $L:=E_8$.
Let $d:=det(M)$, a power of 2 and $k:=rank(M)$.
Note that $\dg M$ is elementary abelian of rank at most $k$.
If
$d=2^k$, then $\frac 1{\sqrt 2}M$ is unimodular, hence is isometric to $\ZZ^k$, and the conclusion holds.   So, we assume that
$d<2^k$.
For $n\le 4$ and $d|8$, it is straightforward to check that the Hermite function $H$ satisfies
$H(n,d)<4$.
Therefore, $M$ contains a root, say $r$.

Suppose that $M$ is a direct summand of $L$.
By \refpp{subrssd},
$N:=ann_M(r)$ is RSSD in $L$, hence is SSD in $L$ by \refpp{ssd=rssd} and we apply induction to conclude that
$N$ is an orthogonal sum of $A_1s$.
So $M$ contains $M'$ an orthogonal sum of $A_1$s with index 1 or 2.
Furthermore, $det(M')=2^k$.
If the index were 1, we would be done, so we assume the index is 2.
Since $d>1$, $d=2, 4$ or 8.
By the index formula for determinants, $2^2$ is a divisor of $d$.
Therefore, $d=4$ or 8.  However, if $d=8$, then
$det(M')=32$, which is impossible since $rank(M')\le 4$.
Therefore, $d=4$ and $rank(M')=4$.
It is trivial to deduce that $M\cong D_4$.

We now suppose that $M$ is not a direct summand of $L$.  Let $S$ be the direct summand of $L$ determined by $M$.  Then $S$ is SSD and the above analysis says $S$ is isometric to some $A_1^m$ or $D_4$.   The only SSD sublattices of $A_1^m$ are the orthogonal direct summands.   The only SSD sublattices of $D_4$ which are proper have determinant $2^4$ and so equal twice their duals and therefore are isometric to $A_1^4$.
\eop

\begin{lem}\labtt{ssde8}
Suppose that $M$ is an SSD sublattice of $E_8$.   Then $M$ is one of the sublattices in Table \ref{ssdsummand} and Table \ref{ssdnonsummand}.
\begin{table}[bht]
\caption{ \bf  SSD sublattices of $E_8$ which span direct summands}
\begin{center}
\begin{tabular}{|c|c|}
\hline
Rank & Type \\ \hline
$0$& $0$ \\ \hline
$1$& $A_1$ \\ \hline
$2$& $A_1\perp A_1$ \\ \hline
$3$& $A_1\perp A_1\perp A_1$ \\ \hline
$4$& $A_1\perp A_1\perp A_1\perp A_1$,  $D_4$ \\ \hline
$5$& $D_4\perp A_1$ \\ \hline
$6$& $D_6$ \\ \hline
$7$& $E_7$ \\ \hline
$8$& $E_8$ \\ \hline
\end{tabular}
\end{center}
\label{ssdsummand}
\end{table}%

\begin{table}[bht]
\caption{ \bf  SSD sublattices of $E_8$ which do not span direct summands}
\begin{center}
\begin{tabular}{|c|c|c|}
\hline
Rank & Type & contained in the summand\\ \hline
$4$& $A_1\perp A_1\perp A_1\perp A_1$ & $D_4$ \\ \hline
$5$& $A_1^{\perp 5}$ & $D_4\perp A_1$ \\ \hline
$6$& $A_1^{\perp 6}$, $D_4\perp A_1\perp A_1$ & $D_6$ \\ \hline
$7$& $A_1^{\perp 7}$, $D_4\perp A_1\perp A_1\perp A_1$, $D_6\perp A_1$ & $E_7$ \\ \hline
$8$& $A_1^{\perp 8}$, $D_4\perp A_1^{\perp 4}$, $D_4\perp D_4$,&$E_8$ \\
&  $D_6\perp A_1\perp A_1$, $E_7\perp A_1$ & \\ \hline
\end{tabular}
\end{center}
\label{ssdnonsummand}
\end{table}
\end{lem}

\pf  We may assume that $1\le rank(M)\le 7$.
First we show that $M$ contains a root.

If $rank(M)\le 4$, this follows from \refpp{smallrankssde8}.
If $rank(M)\ge 4$, then $N:=ann_L(M)$ has rank at most 4, so is isometric to one of
$A_1^k$ or $D_4$.

Suppose that $rank(N)=4$.  If
 $N\cong A_1^4$ and so $ann_L(N)\cong A_1^4$, which contains $M$ and whose only SSD sublattices are orthogonal direct summands, so $M=ann_L(N)$ and the result follows in this case.
If $N\cong D_4$, then $M\cong D_4$ or $A_1^4$ by an argument in the proof of
\refpp{smallrankssde8}.

We may therefore assume that $rank(N)\le 3$, whence $N\cong A_1^{rank(N)}$ and
$rank(M)\ge 5$.
Furthermore, we may assume that $rank(M) > rank(\dg M)$, or else we deduce that
$M\cong A_1^{rank(M)}$.  It follows that $det(M)$ is a proper divisor of 128.

Note that $\dg M$ has rank which is congruent to $rank(M)$ mod 2 (this follows
from the index determinant formula plus the fact that $\dg M$ is an elementary
abelian 2-group).
Therefore, since $rank(M)\le 7$,
$det(M)$ is a proper divisor of 64, i.e. is a divisor of 32.

For any $d\ge 2$, $H(n,d)$ is an increasing function of $n$ for $n\in [5,\infty )$.
For fixed $n$, $H(n,d)$ is increasing as a function of $d$.
Since $H(7,32)=3.888997243\dots $, we conclude that $M$ contains a root, say $r$.

Since $L/(M\perp N)$ is an elementary abelian 2-group by \refpp{involonlattice},
$M\perp N \ge 2L$.  Also, $r+2L$ contains a frame, $F$, a subset of 16 roots
which span an $A_1^8$-sublattice of $L$.  Since roots are orthogonally
indecomposable in $L$, $F=(F\cap M) \cup (F\cap N)$.
It follows that $M$ contains a sublattice $M'$ spanned by $F\cap M$, $M'\cong A_1^{rank(M)}$, and so  $M$ is generated by $M'$ and glue vectors of the form
$\half (a+b+c+d)$, where $a, b, c, d$ are linearly independent elements of $F\cap M$.
It is now straightforward to obtain the list in the conclusion by considering the cases of rank 5, 6 and 7 and subspaces of the binary length 8 Hamming code.
\eop

\begin{lem} \labtt{DD6inD6}Let $X\cong D_6$ and let $\mathfrak{S}=\{ Y\subset X|\ Y\cong DD_6\}$. Then $O(X)$ acts transitively on $\mathfrak{S}$.
\end{lem}

\pf Let $X\cong D_6$ and $R=2X^*$. Since $D_6^*/D_6\cong \ZZ_2\times \ZZ_2$, we have
$X\ge  R\ge  2X$ and $R/2X\cong D_6^*/D_6\cong \ZZ_2\times \ZZ_2$. Thus, the index of $R$ in $X$ is $2^6/2^2=2^4$.

Let $\bar{\ } :X \to X/2X$ be the natural projection. Then for any
$Y\in \mathfrak{S}$, $\bar{Y}$ is a totally isotropic subspace of
$\bar{X}$.   Note that $R/2X$ is the radical of $\bar{X}$ and thus $X/R\cong
2^4$ is nonsingular. Therefore, $\dim (Y+R)/R\leq 2$ and $Y/(Y\cap
R)\cong (Y+R)/R$ also has dimension $\leq 2$.

First we shall show that $Y\geq 2X$ and $\dim (Y+R)/R = 2$.
Consider the tower
\[
Y\geq Y\cap R \geq Y\cap 2X.
\]
Since $Y\cap R+2X$ is doubly even but $R$ is not, $R\neq Y\cap R+
2X$ and $(Y\cap R+ 2X)/ 2X \lvertneqq R/2X$. Thus $(Y\cap R+ 2X)/
2X\cong Y\cap R/ Y\cap 2X$ has dimension $\leq 1$ and hence
\[
|Y: Y\cap 2X|=|Y: Y\cap R|\cdot |Y\cap R: Y\cap 2X|\leq 2^3.
\]
However, $\det(Y)=2^8$ and $\det(2X)=2^{12}4=2^{14}$. Therefore,
$|Y: Y\cap 2X|\geq 2^3$ and hence $|Y: Y\cap 2X|=2^3$. This implies
$Y\cap 2X=2X$, i.e., $Y\geq 2X$. It also implies that $|Y+R:R|=|Y:
Y\cap R|=2^2$ and hence $(Y+R)/R$ is a maximal isotropic subspace of
$X/R\cong \ZZ_2^4$.

Let $2X\leq  M\leq X$ be such that $M/2X$ is maximal totally
isotropic subspace. Then  $M\geq R$ and $\drt{M}$ is an integral
lattice. Set $ M_{even}= \{ \a\in M|\, \half (\a,\a) \in 2\ZZ\}. $
Then $M_{even}$ is a sublattice of $M$ of index $1$ or $2$. If $Y$
is contained in such $M$, then $Y=M_{even}$. That means $Y$ is
uniquely determined by $M$.

Finally, we shall note that the Weyl group acts on $X/R$ as the symmetric group $Sym_6$.
Moreover, $Sym_6$ acts faithfully on $X/R\cong \ZZ_2^4$ and fixes the form $(\ ,\ )$, so it acts as $Sp(4,2)$. Thus it acts transitively on maximal totally isotropic subspace and we have the desired conclusion. \eop

\begin{lem}\labtt{sqrt2D} Let $X\cong D_{6}$ and let $Y\cong DD_6$ be a sublattice of $X$.
Then there exists a subset $\{\eta_1, \dots, \eta_{6}\}\subset X$ with $(\eta_i, \eta_j)=2\delta_{i,j}$ such that
\[
Y=\mspan_\ZZ\{
\eta_i \pm \eta_j|\ i,j= 1, \dots, 6 \}
\] and
\[
X=\mspan_\ZZ\left \{ \eta_1, \eta_2, \eta_4, \eta_6,
\half(-\eta_1+\eta_2-\eta_3+\eta_4), \half(-\eta_3+\eta_4-\eta_5+\eta_6)\right \}.
\]
\end{lem}

\pf We shall use the standard model for $D_{6}$, i.e.,
\[
D_{6}=\{ (x_1, x_2,\dots,  x_{6})\in \ZZ^{6}|\, x_1+\cdots +x_{6} \equiv 0 \mod 2\}.
\]
Let $\b_1=(1,1,0,0,0,0)$, $\b_2=(-1,1,0,0,0,0)$, $\b_3=(0,0,1,1,0,0)$, $\b_4=(0,0,-1,1,0,0)$, $\b_5=(0,0,0,0,1,1)$, and $\b_6=(0,0, 0,0,-1,1)$. Then, $(\b_i, \b_j)=2\delta_{i,j}$ and
\[
W=\mspan_\ZZ\{ \b_i \pm \b_j|\, i,j=1,2,3,4,5,6\}\cong DD_6.
\]
Note also that $\{(1,1,0,0,0,0)$, $(-1,1,0,0,0,0)$,  $(0,0,-1,1,0,0)$, $(0,0, 0,0,-1,1)$, $(-1,0,-1,0,0,0)$, $(0,0,-1,0,-1,0)\}$ forms a basis for $X$ (since their Gram matrix has determinant 4). By expressing them in $\b_1, \dots, \b_6$, we have
\[
X =\mspan_\ZZ\left \{ \b_1, \b_2, \b_4, \b_6,
\half(-\b_1+\b_2-\b_3+\b_4), \half(-\b_3+\b_4-\b_5+\b_6)\right \}.
\]

Let $Y\cong DD_6$ be a sublattice of $X$. Then by Lemma \ref{DD6inD6}, there exists $g\in O(X)$ such that
$Y= Wg$. Now set $\eta_i=\b_i g$ and we have the desired result.
 \eop

\begin{lem}\labtt{dd4ind4}
Let $X\cong D_4$ and let $Y\cong DD_4$ be a sublattice of $X$. Then
$Y=2 X^*$ and hence $ X\leq  \half Y$.
\end{lem}

\pf 
The radical of the form on $X/2X$ is $2\dual X/2X$.   If $W$ is any $A_2$-sublattice of $X$, its image in $X/2X$ complements $2\dual X/2X$.   Therefore, every element of $X \setminus 2\dual X$ has norm $2 (mod\, 4)$.  It follows that $Y\le 2\dual X$.   By determinants, $Y=2\dual X$. \eop

\begin{lem}\labtt{duala2}
Let $L$ be the $A_2$-lattice, with basis of roots $r, s$.  Let $g\in
O(L)$ and $|g|=3$. (i) Then $\dual L \cong \frac 1{\sqrt 3}L$ and
every nontrivial coset of $L$ in $\dual L$ has minimum norm
$\tthird$.  All norms in such a coset lie in
$\tthird + 2\ZZ$.  (ii) If $x$ is any root, $ann_L(x)=\ZZ (xg-xg^2)$ and
$xg-xg^2$ has norm 6.
\end{lem}
\pf  (i) The transformation $g:r\mapsto s, s\mapsto -r-s$ is an
isometry of order 3 and $h:=g-g^{-1}$ satisfies $h^2=-3$ and
$(xh,yh)=3(x,y)$ for all $x, y\in \QQ\otimes L$.
Furthermore, $g$ acts indecomposably on $L/3L\cong 3^2$. We have
$\third L = Lh^{-2}\ge Lh^{-1} \ge L$, with each containment having
index 3 (since $h^2=-3$). Since $\dual L$ lies strictly between $L$
and $\third L$ and is $g$-invariant, $\dual L=Lh^{-1}$.

Since $\dual L \cong \frac 1{\sqrt 3}L$, the minimum norm in $\dual
L $ is $\tthird$ by \refpp{duala2}.  The final statement follows since the six roots
$\pm r, \pm s, \pm (r+s)$ fall in two orbits under the action of
$\la g \ra$, the differences $rg^i-sg^i$ lie in $3\dual L $ and $r$
and $-r$ are not congruent modulo $3\dual L $.

(ii) The element $xg-xg^2$ has norm 6 and is clearly in $ann_L(x)$.
The sublattice $\ZZ x \perp \ZZ (xg-xg^2)$ has norm $2\cdot 6=12$,
so has index 3 in $L$, which has determinant 3.  Since $L$ is
indecomposable, $ann_L(x)$ is not properly larger than $\ZZ
(xg-xg^2)$. \eop

\begin{lem}\labtt{index3ina2}
Let $X\cong A_2$ and let $Y\le X, |X:Y|=3$. Then either $Y=3\dual X$
and its Smith invariant sequence $3, 9$; or $Y$ has Gram matrix
$\begin{pmatrix}2&-3\cr -3&18
\end{pmatrix}$, which has Smith invariant sequence $1, 27$.  In particular, such $Y$ has $\dg Y$ of rank 2 if and only if $Y=3\dual X$.
\end{lem}
\pf Let $r, s, t$ be roots in $X$ such that $r+s+t=0$.  Any two of
them form a basis for $X$.  The sublattices $span\{r,3s\}, span\{s,
3t\}, span\{t, 3r\}$ of index 3 are distinct (since their sets of
roots partition the six roots of $X$) and the index 3 sublattice
$3\dual X$ contains no roots.  Since there are just four sublattices
of index 3 in $X$, we have listed all four.  It is straightforward
to check the assertions about the Gram matrices.  Note that $3\dual
X \cong \sqrt 3 X$. \eop

\begin{prop}\labtt{det3ine8}
Suppose that $M$ is a sublattice of $L\cong E_8$, that $M$ is a
direct summand of $M$, that $M$ has  discriminant group which is elementary abelian of order $3^s$, for some $s$.
Then $M$ is 0, $L$, or is a natural $A_2$, $A_2\perp A_2$, or $E_6$
sublattice.  The respective values of $s$ are 0, 0, 1, 2 and 1. In
case $M$ is not a direct summand, the list of possibilities expands
to include $A_2\perp A_2\perp A_2$, $A_2\perp A_2\perp A_2\perp A_2$ and $A_2\perp E_6$ sublattices.
\end{prop}
\pf
One of $M$ and $ann_L(M)$ has rank at most 4 and the images of $L$ in their
discriminant group are isomorphic.  Therefore, $s\le 4$.   If $s$ were equal to
4, then both $M$ and $ann_L(M)$ would have rank $4$, and each would be
isometric to $\sqrt 3$ times some rank 4 integral unimodular lattice.   By
\refpp{rank8}, each would be isomorphic to $\sqrt 3 \ZZ^4$, which would
contradict their evenness.   Therefore, $s\le 3$.

The second statement is easy to derive from the first, which we
now prove.

We may replace $M$ by its annihilator in $L$ if necessary to assume
that $r:=rank(M)\le 4$. Since $M$ is even and $det(M)$ is a power of
3, $r$ is even.  We may assume that $r\ge 2$ and that $s\ge 1$. If
$r=2$, $M\cong A_2$ \refpp{rank2}. We therefore may and do assume
that $r=4$.

For $s\in \{1,3\}$, we quote \refpp{norank4det3} and \refpp{norank4det27} to see that there is no such $M$.
If $s=2$ we quote \refpp{rank4det9} to identify $M$.
\eop

\begin{prop}\labtt{det5ine8}
Suppose that $M$ is a sublattice of $L\cong E_8$, that $M$ is a
direct summand, that $M$ has discriminant group which is elementary abelian of order $5^s$, for some $s \le 4$.
Then $M$ is 0,  a natural $A_4$ sublattice, the rank 4 lattice $M(4,25)$
(cf. \refpp{m(4,25)}), the rank 4 lattice $\sqrt{5}\dual {A_4}\cong A_4(1)$
(cf. \refpp{a4(1)}) or $L$.  The respective
values of $s$ are 0,  1, 2, 3 and 0.
\end{prop}
\pf We may replace $M$ by its annihilator in $L$ if necessary to
assume that $r:=rank(M)\le 4$. Since $M$ is even and $det(M)$ is a
power of 5, $r$ is even.  We may assume that $r\ge 2$ and that $s\ge
1$. If $r=2$, $det(M)\equiv 3(mod \, 4)$ (consider the form of a Gram matrix
$\begin{pmatrix} a&b\cr b&d \end{pmatrix}$, which has even entries on the diagonal and odd determinant, whence $b$ is odd).   This is not possible
since $det(M)$ is a power of 5.

We therefore may and do assume that $r=4$.

Suppose that $s=4$, i.e., that $\dg M \cong 5^4$.
Then $M\cong \sqrt 5 J$, where $J$ is an integral lattice of
determinant 1.  Then $J\cong \ZZ^{rank(M)}$, which is not an even
lattice.  This is a contradiction since $M$ is even.
We conclude $s=rank(\dg M)\le 3$.  If $s=0$, $M$ is a rank 4 unimodular integral lattice, hence is odd by \refpp{rank8}, a contradiction.  Therefore,
$1\le s \le 3$.   The results
\refpp{rank4det5}, \refpp{rank4det25} and \refpp{rank4det125} identify $M$.
\eop

\begin{nota} \labttr{a4(1)}  The lattice $A_4(1)$ is defined by the Gram matrix
$$\begin{pmatrix} 4&-1&-1&-1\cr -1&4&-1&-1\cr -1&-1&4&-1\cr
-1&-1&-1&4
\end{pmatrix}.$$  It is spanned by vectors $v_1, \cdots, v_5$ which satisfy
$v_1+v_2+v_3+v_4+v_5=0$ and $(v_i,v_j)=-1+5\kron ij$.   Its isometry
group contains $Sym_5\times \la -1\ra$.
\end{nota}

\begin{lem}\labtt{norm4gorbit}
Suppose that $u$ is a norm 4 vector in an integral lattice $U$ where
$D$ acts as isometries so that $1+g+g^2+g^3+g^4$ acts as 0.   Then
one of three possibilities occurs.

(i) The unordered pair of scalars $(u,ug)=(u,g^4u)$ and
$(u,ug^2)=(u,ug^3)$ equals the unordered set $\{0,-2\}$; or

(ii) The unordered pair of scalars $(u,ug)=(u,g^4u)$ and
$(u,ug^2)=(u,ug^3)$ equals the unordered set $\{-3,1\}$; or

(iii)  $(u,ug)=(u,ug^2)=(u,ug^3)=(u,g^4u)=-1$.

The isometry types of the lattice $span\{u,ug, ug^2, ug^3,ug^4\}$
in these respective cases are $AA_4, A_4, A_4(1)$.
\end{lem}
\pf Straightforward. \eop

\begin{lem}\labtt{allabouta4(1)}
Let $X\cong A_4(1)$ \refpp{a4(1)}.  Then

(i) $X$ is rootless and contains exactly 10 elements of norm 4;

(ii) Suppose $u\in X$ has norm 4.   Then $ann_X(u)\cong \sqrt 5 A_3$.

(iii) $O(X)\cong 2\times Sym_5$; furthermore, if $X_4$ is the set of
norm 4 vectors and $\mathcal O$ is an orbit of a subgroup of order 5
in $O(X)$ on $X_4$, then the subgroup of $O(X)$ which preserves
$\mathcal O$ is a subgroup isomorphic to $Sym_5$, and is the
subgroup generated by  all reflections.

(iv)  Suppose that $Y\cong  A_4$ and that $g\in O(Y)$ has order 5,
then $Y(g-1)\cong X$.  Also, $O(Y)\cap O(Y(g-1))=N_{O(Y)}(\la g \ra )\cong 2\times 5{:}4$, where the right direct factor is a Frobenius group of order 20.

(v) $\dg{X}$ is elementary abelian.
\end{lem}
\pf (i) If the set of roots $R$ in
$X$ were nonempty, then $R$ would have an isometry of order 5.  Since the
rank of $R$ is at most 4, $R$ would be an $A_4$-system and so the sublattice of $X$ which $R$ generates would be $A_4$, which has determinant 5.  Since $det(X)=5^3$,
this is a contradiction.

By construction, $X$ has an cyclic group $Z$ of order 10 in $O(X)$
which has an orbit of 10 norm 4 vectors, which are denoted $\pm v_i$ in
\refpp{a4(1)}. Suppose that $w$ is a norm 4 vector outside the
previous orbit.  Let $g\in Z$ have order 5.  Then
$w+wg+wg^2+wg^3+wg^4=0$.   Therefore $0=(v,w+wg+wg^2+wg^3+wg^4)$,
which means that there exists an index $i$ so that $(w,vg^i)$ is
even. Since $vg^i$ and $w$ are linearly independent, $(vg^i,w)$ is
not $\pm 4$ and the sublattice $X'$ which  $vg^i$ and $w$ span has
rank 2.  Since $(vg^i,w)\in \{-2,0,2\}$, $X'\cong AA_1\perp AA_1$ or
$AA_2$.  This contradicts \refpp{prankdg}  (in that notation, $n=4$, $m=2$,
$p=5$,
$r=3$).

(ii) Let $K:=ann_X(u)$.   Since $(u,u)=4$ is relatively prime to
$det(X)$, the natural map $X\rightarrow \dg {\ZZ u}$ is onto.
Therefore $\ZZ u \perp K$ has index 4 in $X$.  Hence $\det(\ZZ u
\perp K)=4^2\cdot 5^3$  and $\det K=4\cdot 5^3$. Since $\dg{X}\leq
\dg{\ZZ u \perp K}=\dg{\ZZ u}\times \dg{K}$ and the Smith invariant
sequence of $X$ is 1,5,5,5, $\dg{K}$ contains an elementary group
$5^3$. Moreover, the image of $X$ in $\dg{K}$ isomorphic to the
image of $X$ in $\dg{\ZZ u}$, which is isomorphic to $\ZZ_4$.
Therefore,  $\dg{K}=\ZZ_4\times \ZZ_5^3$ by determinants. Hence, the
Smith invariant sequence for $K$ is $5, 5, 20$ and so $K\cong \sqrt
5 W$, for an integral lattice $W$ such that $\dg W=4$. Since $X$ is
even, $W$ is even.  We identify $W$ with $A_3$ by \refpp{rank3}.

(iii) We use the notation in the proof of (i).  By \refpp{norm4gorbit}, for any two
distinct vectors of the form $vg^i$, the inner product is $-1$, so
the symmetric group on the set of all $vg^i$ acts as isometries on
the $\ZZ$-free module spanned by them, and on the quotient of this
module by the $\ZZ$-span of $v+vg+vg^2+vg^3+vg^4$, which is
isometric to $X$.

Pairs of elements of norm 4 fall into classes according to their
inner products: $\pm 4, \pm 1$.  An orbit of an element of order 5
on $X_4$ gives pairs only with inner products $4, -1$ (since the sum of these five values is 0).  There are
two such orbits and an inner product between norm 4 vectors from
different orbits is one of  $-4, 1$.  The map $-1$ interchanges
these two orbits.  Therefore, the stabilizer of $\mathcal O$ has
index 2 in $O(X)$.  It contains the map which interchanges distinct
$vg^i$ and $vg^j$ and fixes other $vg^k$ in the orbit.  Such a map
is a reflection on the ambient vector space.  Since $Sym_5$ has just
two classes of involutions, it is clear that every reflection  in
$O(X) = \la -1\ra \times Stab_{O(X)}(\mathcal O )$ is contained in $
Stab_{O(X)}(\mathcal O )$.

(iv) We have $Y(g-1)=span_{\ZZ}\{ vg^i-vg^j | i, j \in \ZZ \}$.
By checking a Gram matrix, one sees that it is isometric to $X$.
We consider $O(Y)\cap O(Y(g-1))$, which clearly contains $N_{O(Y)}(\la g \ra)$.  We show that
this containment is equality.  We take for $Y$ the standard model, the set of coordinate sum 0 vectors in $\ZZ^5$.   Take $v\in Y(g-1)$, a norm 4 vector.  It has shape $(1,1,-1,-1,0)$ (up to reindexing).   The coordinate permutation $t$ which transposes the last two coordinates is not in $O(Y(g-1))$ (since $v(t-1)$ has norm 2).   Therefore $O(Y)$ does not stabilize $Y(g-1)$.    Since
$N_{O(Y)}(\la g \ra)$ is a maximal subgroup of $O(Y)$, it equals
$O(Y)\cap O(Y(g-1))$.

(v) Since $A_4(1)\cong \sqrt{5} A_4^*$ by \refpp{rank4det125},
$\dual{(A_4(1))}\cong \frac{1}{\sqrt{5}} A_4$. Thus,
$ 5 \dual{A_4(1)} < A_4(1)$ and $\dg{A_4(1)}$ is elementary abelian.
\eop

\begin{lem}\labtt{a41ine8}
Let $X\cong A_4(1)$ be a sublattice of $E_8$. If $X$ is a direct summand, then
$ann_{E_8}(X)\cong A_4(1)$.
\end{lem}

\pf.  Let $Y\cong E_8$ and let $X\cong A_4(1)$ be a sublattice of $Y$.
Since $X$ is a direct summand, the natural
map $Y \to  \dg{ X}$ is onto. Similarly,
the natural map from $Y \to \dg{ann_Y(X)} $ is also onto and  these two images
are isomorphic.
Thus,
$ \dg{ann_Y( X)}\cong \dg{X}\cong 5^3$. Hence,
$ann_Y( X) $ is isomorphic to $A_4(1)$ by \refpp{rank4det125}.\eop

\begin{rem}\labtt{a41ine82}
Note that $A_4(1)$ can be embedded into $E_8$ as a direct summand. Recall that
\[\begin{split}
&A_4(1)\cong \sqrt{5}A_4^*\\
&=\mspan_\ZZ\{\frac{1}{\sqrt{5}}(1,1,1,1,-4),  \frac{1}{\sqrt{5}}(1,1,1,-4,1), \frac{1}{\sqrt{5}}(1,1,-4,1,1), \frac{1}{\sqrt{5}}(1,-4,1,1,1)\}.
\end{split}
\]
Then,
\[
(A_4(1))^* \cong \frac{1}{\sqrt{5}}A_4= \left \{ \frac{1}{\sqrt{5}}(x_1, \dots, x_5)\left |\  \sum_{i=1}^5 x_i=0 \text { and } x_i\in \ZZ, i=1,\dots,5 \right . \right \}.
\]
Let
\[
Y=\mspan_\ZZ\left \{ {A_4(1) \perp A_4(1), \quad \frac{1}{\sqrt{5}}(1,-1, 0,0,0\, |\, 2,-2, 0,0,0)} \atop
{ \frac{1}{\sqrt{5}}(0,1,-1, 0,0\, |\, 0, 2,-2, 0,0),\  \frac{1}{\sqrt{5}}(0,0, 1,-1, 0\, |\, 0,0, 2,-2, 0,)} \right\}.
\]
Then $Y$ is a rank $8$ even lattice and $|Y: A_4(1)\perp A_4(1)|=5^3$. Thus $\det(Y)=1$ and $Y\cong E_8$.
Clearly, $A_4(1)$ is a direct summand by the construction.
\end{rem}

\begin{rem}\labttr{aboutinvolsoe6}
We have
$O(E_6)\cong Weyl(E_6) \times \la -1 \ra$.
Thus, outer involutions are negatives of inner involutions.  The next result does not treat inner and outer cases differently.
\end{rem}

\begin{lem}\labtt{involsoe6}
Let $t\in O(E_6)$ be an involution.
The negated sublattice for $t$ is either SSD (so occurs in the list for $E_8$ \refpp{ssde8}) or is
RSSD but not SSD and is isometric to one of $AA_2, AA_2\perp A_1, AA_2 \perp A_1\perp A_1, AA_2 \perp A_1\perp A_1 \perp A_1, A_5, A_5 \perp A_1, E_6$.  Moreover,
the isometry types of the RSSD sublattices determine them uniquely  up to the action of $O(E_6)$.
\end{lem}
\pf
Let $S$ be the negated sublattice and assume that it is not SSD.  Then the image of $E_6$ in $\dg S$ has index 3 and is an elementary abelian 2-group, so that $det(S)=2^a 3$,   where $a\le rank(S)$.  Note that $rank(S)\ge 2$.  Now, let $T:=ann_{E_6}(S)$, a sublattice of rank at most 4.  Since $det(S\perp T)=2^{2a}3$, $det(T)=2^a$ and the image of $E_6$ in $\dg T$ is all of $\dg T$.  Therefore, $T$ is SSD and we may find the isometry type of $T$ among the SSD sublattices of $E_8$.  As we search through SSD sublattices of rank at most 4 (all have the form $A_1^m$ or $D_4$), it is routine to determine the annihilators of their embeddings in $E_6$.
\eop

\begin{lem}\labtt{e6overaa2+4}
Suppose that $R\perp Q$ is an orthogonal direct sum with $Q\cong
AA_2$ and $R\cong D_4$. Let $\phi : \dg R \rightarrow \dg Q$ be any
monomorphism (recall that $\dg R\cong 2\times 2$ and $\dg Q\cong
2\times 2 \times 3$). Then the lattice $X$ which is between $R\perp
Q$ and $\dual R \perp \dual Q$ and which is the diagonal with
respect to $\phi$ is isometric to $E_6$. Furthermore, if $X$ is a
lattice isometric to $E_6$ which contains $R\perp Q$, then $X$ is
realized this way.
\end{lem}
\pf  Such $X$ have determinant 3. The cosets of order 2 for $D_4$ in
its dual have odd integer norms (the minimum is 1). The cosets of
order 2 for $AA_2$ in its dual have odd integer norms (the minimum is
1).  It follows that such $X$ above are even lattices.  By a
well-known characterization, $X\cong E_6$ (cf. \refpp{uniquenesse6}).

Conversely, suppose that $X$ is a lattice containing $R \perp Q$,
$X\cong E_6$.
Since $det(X)$ is odd, the image of the natural map
$X\rightarrow \dg R$ is onto.
Therefore,
$|X:R\perp Q|=4$.  The image of $X$ in $\dg Q$ is isomorphic to the image of $X$ in
$\dg R$.  The last statement follows.
 \eop

\begin{coro}\labtt{annue6}
(i) Let $Y$ be a sublattice of $X\cong E_6$ so that
$Y\cong D_4$.  Then $ann_X(Y)\cong AA_2$.

(ii) Let $U$ be a sublattice of $X\cong E_6$ so that
$U\cong AA_2$ and $X/(U\perp ann_X(U))$ is an elementary abelian 2-group.  Then $X/(U\perp ann_X(U))\cong 2^2$ and $ann_X(U)\cong D_4$.
\end{coro}
\pf
(i)  Let $Z:=ann_X(Y)$.  Since $(det(X),det(Y))=1$, the natural map of $X$ to
$\dg Y\cong 2\times 2$ is onto, so the natural map of $X$ to $\dg Z \cong
2\times 2 \times 3$ has image $2 \times 2$.  Since $rank(Z)=2$, this means
$\frac 1{\sqrt 2} Z$ is an integral lattice of determinant 2.  It is not
rectangular, or else there exists a root of $X$ whose annihilator contains $Y$,
whereas a root of $E_6$ has annihilator which is an $A_5$-sublattice, which does
not contain a $D_4$-sublattice (since an $A_5$ lattice does not contain an
$A_1^4$-sublattice).   Therefore, by \refpp{rank2},
$\frac 1{\sqrt 2} Z \cong A_2$.

(ii)  Use \refpp{rank4det4}.
\eop

\begin{nota}\labttr{xq}  We define two rank 4 lattices $X, Q$.   First, $X\cong A_1^2A_2$, $\dg X \cong 2^2\times 3$. Let $X$ have the decomposition into indecomposable summands $X=X_1\perp X_2\perp X_3$, where $X_1\cong X_2\cong A_1$ and $X_3\cong A_2$.
Let $\a_1\in X_1, \a_2\in X_2, \a_3, \a_4\in X_3$ be roots with $(\a_3,\a_4)=-1$.

We define $Q \cong ann_{E_6}(P)$, where $P$ is a sublattice of $E_6$ isometric to $A_1^2$.   Then $\dg Q \cong 2^2\times 3$ and $rank(Q)=4$.  Then $Q$ is not a root lattice (because in $E_6$, the annihilator of an $A_1$-sublattice is an $A_5$-sublattice; in an $A_5$-lattice, the annihilator of an $A_1$-sublattice is not a root lattice).  

We use the standard model for $E_6$, the annihilator in the standard model of $E_8$ of
$J:=span\{(1,-1,0,0,0,0,0,0),(0,1,-1,0,0,0,0,0)\}$.
So, $E_6$ is
the set of $E_8$ vectors with equal first three coordinates.

We  may take $P$ to be the span of $(0,0,0,1,1,0,0,0)$ and $(0,0,0,1,-1,0,0,0)$.   Therefore,
$Q=span\{u, Q_1,  w\}$, where $u=(2,2,2,0,0,0,0,0)$, $Q_1$ is the $D_3$-sublattice supported on the last three coordinates, and $w:=(1,1,1,0,0,1,1,1)$.
\end{nota}

\begin{lem}\labtt{aboutoq}
The action of $O(Q_1)\cong 2\times Sym_3$ extends to an action on $Q$.   This action is faithful on $Q/3\dual Q$.
\end{lem}
\pf
The action of $O(Q_1)\cong 2\times Sym_3$ extends to an action on $Q$ by letting reflections in roots of $Q_1$ act trivially on $u$ and by making the central involution of $O(Q_1)$ act as $-1$ on $Q$.  The induced action on $Q/3\dual Q$ is faithful since $Q_1$ maps onto $Q/3\dual Q$ (because $(3, det(Q_1))=1$) and the action on $Q_1/3Q_1$ is faithful.  In more detail, the action of $O_2(O(Q_1))\cong 2^3$ is by diagonal matrices and any normal subgroup of $O(Q_1)$ meets   $O_2(O(Q_1))$ nontrivially.
\eop

\begin{lem}\labtt{rank4inrank4}  We use notation \refpp{xq}.   Then
$X$ contains a sublattice $Y\cong \sqrt 3 Q$ and $X> Y> 3X$.
\end{lem}
\pf
We define $\b_1:=\a_1+\a_2+\a_3, \b_2:=-2\a_3-\a_4, \b_3:=\a_3+2\a_4$.   Then $Y_1:=span\{\b_1, \b_2, \b_3\}\cong \sqrt 3 D_3$.  The vector $\b_4:=3\a_1-3\a_2$ is orthogonal to $Y_1$ and has norm 36.  Finally, define
$\g:=\half \b_4 + \half(\b_1+2\b_2+3\b_3)=2\a_1-\a_2+2\a_4$.  Then  $Y:=span\{Y_1, \b_4, \g\}$ is the unique lattice containing $Y_1\perp \ZZ \b_4$ with index 2 whose intersection with
$\half Y_1$ is $Y_1$ and whose intersection with $\half \ZZ \b_4$ is $\ZZ \b_4$.  There is an analogous characterization for $Q$ and $\sqrt 3 Q$.  We conclude that $Y \cong \sqrt 3 Q$.

Moreover, by direct calculation, it is easy to show that
\[
\begin{split}
3\a_1&= \g+\b_1-\b_3, \quad  3\a_2= \g+\b_1-\b_3-\b_4,\\
3\a_3&=\b_1+2\b_3+\b_4-2\g, \quad 3\a_4=2\g -(\b_1+\b_2+\b_3+\b_4).
\end{split}
\]
Hence, $Y$ also contains $3X$.
\eop

\begin{lem}\labtt{noover}
Let $\mathcal M$ be the set of rank $n$ integral lattices. For $q\in \ZZ$,  let $\mathcal M(q)$ be the set of $X\in \mathcal M$ such that $X\le q \dual X$.
Suppose that $q$ is a prime, $X, Y \in \mathcal M(q)$, $Y\ge X$  and
$q$ divides $|Y:X|$.
Then $q^{n+2}$ divides $det(X)$.  In particular, if
$q^{n+2}$ does not divide $det(X)$, then  $X$ is not properly contained in a member of $\mathcal M(q)$.
\end{lem}
\pf
Use the index formula for determinants of lattices.
\eop

\begin{prop}\labtt{transxq}
For an integral lattice $K$, define $\tilde K:=K+2\dual K$.
Let $$\mathcal A:=\{ (R,S) | R\le S\le \RR^4, R\cong \sqrt 3 Q, S\cong
X\}$$  $$\mathcal B:=\{ (T,U) | T\le U\le \RR^4, T\cong \sqrt 3 X,
U\cong Q\}$$
$$\mathcal A':=\{ (R,S) | R\le S\le \RR^4, R\le 3\dual R,
det(R)=2^23^5, S\cong X\}$$
 $$\mathcal B':=\{ (T,U) | T\le U\le \RR^4, T\le 3T^*, det(T)=
2^23^5, U\cong Q\}.$$

Then
(i) $\mathcal A=\mathcal A'\ne \emptyset$ and $\mathcal B=\mathcal B'\ne \emptyset$;

(ii) the map $(T, U)\mapsto (\sqrt 3 U, \frac 1{\sqrt 3}T)$ gives a
bijection from $\mathcal B$ onto $\mathcal A$; furthermore if
$(T,U)\in \mathcal B$, then $T\ge 3\tilde U$ and if
$(R,S)\in \mathcal A$, then $R\ge 3\tilde  S$;

(iii)
$O(\RR^4)$ has one orbit on $\mathcal A$ and on $\mathcal B$.
\end{prop}

\pf
Clearly, $\mathcal A \subseteq \mathcal A'$ and $\mathcal B\subseteq
\mathcal B'$.
From \refpp{rank4inrank4},  $\mathcal A\neq \emptyset$ and  $\mathcal B\neq \emptyset$. Moreover, the formula in (ii) gives a bijection between $\mathcal A$ and $\mathcal B$.

Now, let $(E, F)$ be in $\mathcal A'  \cup \mathcal B '$.

We claim that $3\tilde F=F \cap 3\dual F$.
We prove this with the theory of modules over a PID.
Since $\dg F \cong 2^2 \times 3$, there exists a basis $a, b, c, d$ of $\dual F$
so that $a, b, 2c, 6d$ is a basis of $F$.   Then $a, b, 2c, 2d$ is a basis of $\tilde F$.    Since $3\dual F$ has basis $3a, 3b, 3c, 3d$, $F\cap 3\dual F$ has basis $3a, 3b, 6c, 6d$.   The claim follows.

Note that $F/3\tilde F$ is an elementary abelian 3-group of rank 3 and the claim implies that it is a nonsingular quadratic space.   Therefore, its totally singular subspaces have dimension at most 1.

We now study $E':=E+3\tilde F$, which maps onto a totally singular subspace of
$F/3\tilde F$.  Since totally singular subspaces have dimension at most 1, $|F:E'|$ is divisible by $3^2$ and so its determinant is
$|F:E'|^2det(F)=|F:E'|^2 2^2 3$.   However, $E'$ contains $E$, which has
determinant $2^23^5$.  We conclude that $E=E'$ has index $9$ in $F$.
Therefore, $E=E' \ge  3\tilde F$.
The remaining parts of (ii) follow.

Let $(T, U)\in \mathcal B$.   The action of $O(U)\cong Sym_4 \times 2$
on $U/3\tilde{U}$ is that of a monomial group with respect to a basis of
equal norm nonsingular vectors \refpp{aboutoq}.  It follows that the
action is transitive on maximal totally singular subspaces, of which
$T/3\tilde{U}$ is one.  This proves transitivity for $\mathcal B$.
Therefore  $\mathcal B=\mathcal B'$ and, using the $O(U)$-equivariant
bijection (ii),
$\mathcal A =\mathcal A'$.
 \eop

 \begin{coro}\labtt{transxq2} $\sqrt 3 Q$ does not embed in $Q$ and $\sqrt 3 X$ does not embed in $X$.
 \end{coro}
 \pf
 Use \refpp{transxq}, \refpp{rank4inrank4} and the fact that $X$ is not isometric to $Q$ ($X$ is a root lattice and $Q$ is not).
 \eop


\begin{lem}\labtt{e8overa2e6}
Suppose that $S\perp T$ is an orthogonal direct sum with $S\cong A_2,
T\cong E_6$.   The set of $E_8$ lattices which contain $S\perp T$ is in
bijection with
$\{ X|\, S\perp T\leq X\leq S^*\perp T^*, |X:S\perp T|=3, S^*\cap X=S, T^*\cap X=T\}$.
\end{lem}
\pf This is clear since any $E_8$ lattice containing $S\perp T$ lies
in $\dual S \perp \dual T$ and since the nontrivial cosets of $S$ in
$\dual S$ have norms in $\tthird +2\ZZ$ and the nontrivial cosets
of $T$ in $\dual T$ have norms in $\frac 43 + 2\ZZ$. \eop

\begin{lem}\labtt{aa1ind4}
Let $X\cong D_4$ and let $H\cong AA_1$ be a sublattice of $X$. Then the image of
the natural map $X^*$ to $H^*$ is $H^*=\frac{1}4 H$.
\end{lem}
\pf
A generator of $H$ has norm 4, so $H$ is a direct summand of $L$.  In general, if $W$ is a lattice and $Y$ is a direct summand of $W$, the natural map $\dual W \rightarrow \dual Y$ is onto.  The lemma follows.
\eop

\begin{lem}\labtt{aa2ind4} (i)
Up to the action of the root reflection group of $D_4$, there is a unique embedding of $AA_2$ sublattices.

(ii) We have transitivity of $O(D_4)$ on the set of $A_2$ sublattices and on the set of $AA_2$-sublattices.
In $D_4$, the annihilator of an $AA_2$ sublattice is an $A_2$ sublattice, and the annihilator of an $A_2$ sublattice is an $AA_2$-sublattice.
\end{lem}
\pf
(i)
Let $X\cong D_4$ and $Y\cong AA_2$.   Since every element of $Y$ has norm divisible by 4, $Y\le 2\dual X$.  Now let $s: =f-1$, where $f\in O_2(Weyl(X))$, $f^2=-1$.
Then $s^{-1}$ takes $2\dual {X}$ to $X$ and takes $Y$ to an $A_2$ sublattice of $X$.  Now use the well-known results that $A_2$ sublattices form one orbit under $Weyl(X)$ and $O(A_2)\cong Dih_{12}$ is induced on an $A_2$ sublattice of $D_4$ by its stabilizer in $Weyl(D_4)$.

(ii) We may take $Y:=span\{(-2,0,0,0),(1,1,1,1)\} \cong AA_2$.  Its annihilator is
$Z:=span\{(0,1,-1,0),(0,0,1,-1)\}\cong A_2$.  Trivially, $ann_X=(Z)=Y$.
\eop

\begin{lem}\labtt{gluee6}
Let $X\cong E_6$ and $Y, Z$ sublattices such that  $Z\cong D_4$ and $Y:=ann_X(Z)\cong AA_2$.  Define  $W:=2\dual Y$  (alternatively, $W$ may be characterized by the property that $Y\le W\le\dual Y$, $W/Y\cong 3$).
Then
$W \le \dual X$.
\end{lem}
\pf
By coprimeness, the natural map $X\rightarrow \dual Z$ is onto, and the image of $X$ in $\dg Z$ has order 4.   Therefore, the image of the natural map $X\rightarrow \dg Y$ has order 4 and so the image of the natural map $X\rightarrow \dual Y$ is $\half Y$.   The dual of $\half Y$ is $2\dual Y$, which contains Y with index 3 and satisfies $(X, 2\dual Y)\le \ZZ$.
\eop

\begin{lem}\labtt{duale6denom}
We have $(\dual {E_6}, \dual {E_6})=\third \ZZ$ and the norms of vectors in $\dual {E_6}\setminus E_6$ are in $\frac 43 +\ZZ$.
\end{lem}
\pf
This follows from the fact that $E_6$ has a sublattice of index 3 which is isometric to $A_2^3$ and the facts that $(\dual {A_2},\dual {A_2})=\tthird \ZZ$ and that a glue vector for $A_2^3$ in $E_6$ has nontrivial projection to the spaces spanned by each of the three summands.
\eop

\begin{hyp}\labttr{hypct}  $L$ is a rank 12 even integral lattice, $\dg L\cong 3^k$, for some integer $k$, $L$ is rootless and $\dual L$ contains no vector of norm $\frac 23$.
\end{hyp}

\begin{lem}\labtt{nonmaxwitt}
The quadratic space $\dg L$ in
\refpp{hypct}
has nonmaximal Witt index if $k$ is even.
\end{lem}
\pf
If the Witt index were maximal for $k$ is even, there would exist a lattice $M$ which satisfies $3L \le 3M \le L$ and $3M/3L$ is a totally singular space of dimension $\frac k2$.  Such an $M$ is even and unimodular.  A well-known theorem says that $rank(M)\in 8\ZZ$, a contradiction.  \eop

\begin{prop}\labtt{12+12inleech}
Let $L, L'$ be two lattices which satisfy hypothesis \refpp{hypct} for $k$ even, and which have the same determinant.
There exists an embedding of $L\perp L'$ into the Leech lattice.
\end{prop}
\pf
We form $L\perp L'$.
The quadratic spaces $\dg L$, $\dg L'$ have nonmaximal Witt index.

Let $g$ be a linear isomorphism  from $\dg L$ to $\dg {L'}$ which
takes the quadratic form on $\dg L$ to the negative of the quadratic form on $\dg {L'}$ \refpp{scalarisometry}.

Now, form the overlattice $J$ by gluing  from $\dg L$ to $\dg {L'}$ with $g$.
Clearly, $J$ has rank 24, is even and unimodular.  The famous characterization of the Leech lattices reduces the proof to showing that $J$ is rootless.

Suppose that $J$ has a root, $s$.  Write $s=r+r'$ as a sum of its projections to the rational spaces spanned by $L, L'$ respectively.  The norm of any element $x\in \dual L$  has the form $a/3$, where $a$ is an even integer at least 4.
 The norm of any element $x\in \dual {L'}$  has the form $b/3$, where $b$ is an even integer at least 4.   Therefore, we may assume that $r, r'$ have respective norms at least $\frac 43$.  Then $(s,s)\ge \frac 83 > 2$, a contradiction.
 \eop

\begin{lem}\labtt{hypctg}
Let $L$ be an even integral rootless lattice with $\dg L\cong 3^k$, for an integer $k$,  and an automorphism $g$ of order 3 without eigenvalue 1 such that $\dual L (g-1)\le L$.
Then $L$ satisfies hypothesis \refpp{hypct}.
\end{lem}
\pf
We need to show that if $v\in \dual L$, then $(v,v)\ge \frac 43$.  This follows since $v(g-1)\in L$, $(v(g-1),v(g-1))=3(v,v)$  and $L$ is rootless.
\eop

\begin{coro}\labtt{12+12inleech2} If
$L$, $L'$ satisfy hypotheses of \refpp{hypctg} and each of $L$, $L'$ is not properly contained in a rank 12 integral rootless lattice (such an overlattice satisfies \refpp{hypctg}), then $L\cong L'$ and $k=6$.
\end{coro}
\pf
Let $\Lambda$  be the Leech lattice.  We use results from \cite{Gr12} which analyze the elements of order 3 in $\Lambda $.

Take two copies $L_1, L_2$ of $L$.
We have by \refpp{12+12inleech}, an embedding of $L_1\perp L_2$ in $\Lambda$.  Identify $L_1\perp L_2$ with a sublattice of
$\Lambda$.

Since
$L_1$, $L_2$ are not properly contained in another lattice which satisfies \refpp{hypctg} and since $\Lambda$ is rootless, $L_1$ and $L_2$  are direct summands of $\Lambda$.
Since they are direct summands, $L_2=ann_{\Lambda } (L_1)$,
$L_1=ann_{\Lambda } (L_2)$ and
the natural maps of $\Lambda$ to $\dg {L_1}$ and $\dg {L_2}$ are onto.
The gluing construction shows that the automorphism $g$ of order 3 in $L$ as in \refpp{hypctg} extends to an automorphism of $\Lambda$ by given action on $L_2$
and trivial action on $L_1$.  Denote the extension by $g$.

We now do the same for $L'$, $g'$ in place of $L$, $g$.

From Theorem 10.35 of \cite{Gr12}, $g$ and $g'$ are conjugate in $O(\Lambda)$ and $det(L)=det(L')=3^6$.   A conjugating element takes the fixed point sublattice $L_1$ of $g$ to the fixed point sublattice $L_1'$ of $g'$.   Therefore, $L$ and $L'$ are isometric.
\eop

\begin{coro}\labtt{ctmaximal}  The Coxeter-Todd lattice is not properly contained in an integral, rootless lattice.
\end{coro}
\pf
Embed the Coxeter-Todd lattice $P$ in a lattice $Q$ satisfying the hypothesis of
\refpp{12+12inleech2}.    Since $det(P)=3^6=det(Q)$, $P=Q$.
\eop

\begin{lem}\labtt{trianglee6a2}
Let $X\cong E_8$, $P\le X, P\cong E_6$ and $Q:=ann_X(P)$.

(i) There exists a sublattice $R\cong A_2$ so that
$R\cap (P\cup Q)$ contains no roots.

(ii) If $r\in R$ is a root, then the orthogonal projection of $r$ to $P$ has norm $\frac 43$ and the projection to $Q$ has norm $\frac 23$.
\end{lem}
\pf
(i)
We may pass to a sublattice $Q_1\perp Q_2\perp Q_3 \perp Q$ of type $A_2^4$, where $P\ge Q_1\perp Q_2\perp Q_3$.   Then $X$ is described by a standard gluing with a tetracode, the subspace of $\FF_3^4$ spanned by $(0,1,1,1), (1,0,1,2)$, and elements $v_i$ of the dual of $Q_i$ ($Q_4:=Q$) where $v_i$ has norm $\frac 23$.  Then for example take
$R$ to be the span of
$ v_2+v_3+v_4, v_1+v_3-w$, where $w\in v_4+Q$ has norm
$\frac 23$ but $(w,v_4)=-\third$.  See \refpp{duala2}.

(ii) This follows since the norms in any nontrivial coset of $Q$ in $\dual Q$ is $\frac 23 +2\ZZ$.
\eop

\section{Values of the Hermite function}

\begin{nota}  Let $n$ and $d$ be positive integers. Define the Hermite function
\[
H(n,d):= \left( \frac{4}3\right)^{\frac{n-1}2} d^{(1/n)}.
\]
\end{nota}

\begin{thm}[Hermite: cf. proof in  \cite{kneser}, p. 83]
If a positive definite rank $n$ lattice has determinant $d$, it contains a nonzero vector of norm $\leq H(n,d)$.
\end{thm}

\begin{table}[bht]
\caption{ \bf  Values of the Hermite function $H(n,d)$; see \cite{kneser}, p.83. }
\begin{center}
\[
\begin{array}{|c|c|c|}
\hline
n & d & H(n,d)\\ \hline \hline
2& 2& 1.632993162\\ \hline
 2 & 3 & 2.000000000 \\ \hline
 2 & 4 & 2.309401077 \\ \hline
 2 & 5 & 2.581988897 \\ \hline
 2 & 6 & 2.828427125 \\ \hline
 2 & 7 & 3.055050464 \\ \hline
 2 & 8 & 3.265986324 \\ \hline
 2& 9 & 3.464101616 \\ \hline
 2& 12& 4.000000000 \\ \hline
 2 & 20 & 5.163977796 \\ \hline
 2 & 24& 5.656854249\\ \hline \hline
 5& 6& 2.543945033 \\ \hline
\end{array}
\,
\begin{array}{|c|c|c|}
\hline
 n & d & H(n,d)\\ \hline \hline
 3& 2& 1.679894733 \\ \hline
 3& 3& 1.922999426 \\ \hline
 3& 4& 2.116534735 \\ \hline
 3& 5& 2.279967929 \\ \hline
 3& 6& 2.422827457 \\ \hline
 3& 8& 2.666666666\\ \hline
 3& 10& 2.872579586 \\ \hline
 3& 12& 3.052571313 \\ \hline
 3& 16& 3.359789466 \\ \hline
 3& 24& 3.845998854\\ \hline
 3& 50& 4.912041997\\ \hline \hline
 6& 3&  2.465284531\\ \hline
\end{array}
\,
\begin{array}{|c|c|c|}
\hline
 n & d & H(n,d)\\ \hline \hline
  4& 2& 1.830904128\\ \hline
 4& 3& 2.026228495\\ \hline
 4& 4& 2.177324216\\ \hline
 4& 5& 2.302240057\\ \hline
 4& 6& 2.409605343\\ \hline
 4& 7& 2.504278443\\ \hline
 4& 8& 2.589289450\\ \hline
 4& 9& 2.666666668\\ \hline
 4& 12& 2.865519818\\ \hline
 4& 25& 3.442651865\\ \hline
 4& 125& 5.147965271\\ \hline \hline
7 & 32& 3.888997243\\ \hline
\end{array}
\]
\end{center}
\label{smallrank}
\end{table}%

\section{Embeddings of NREE8 pairs in the Leech lattice}

If $M, N$ is an NREE8 pair, then except for the case $\dih{4}{15}$,  $L=M+N$ can be
embedded in the Leech lattice $\Lambda$.
 In this
section, we shall describe such embeddings explicitly.

In the exceptional case $\dih{4}{15}$,
 $M\cap N\cong AA_1$ and $|t_Mt_N|=2$.  See  \refpp{rootlessdih415}.

\subsection{Leech lattice and its isometry group}
We shall recall some notations and review certain basic properties
of the Leech lattice $\Lambda$ and its isometry group $O(\Lambda)$, which is also known as $Co_0$, a perfect group of order $2^{22}\cdot 3^9\cdot 5^4\cdot 7^2\cdot 11\cdot 13\cdot 23$.

Let $\Omega=\{1,2,3,\dots, 24\}$ be a set of $24$ element and let
$\mathcal{G}$ be the extended Golay code of length $24$ indexed by $\Omega$. A subset
$S\subset \Omega$ is called a $\mathcal{G}$-set if
$S=\mathrm{supp}\,\a$ for some codeword $\a\in \mathcal{G}$. We
shall identify a $\mathcal{G}$-set with the corresponding codeword
in $\mathcal{G}$. A $\mathcal{G}$-set $\mathcal{O}$ is called an
\textit{octad} if $|\mathcal {O}|=8$ and is called a
\textit{dodecad} if $|\mathcal {O}|=12$. A \textit{sextet} is a
partition of $\Omega$ into six $4$-element sets of which the union
of any two forms a octad. Each $4$-element set in a sextet is called
a \textit{tetrad}.

For explicit calculations,  we shall use the notion of \textit{hexacode balance} to denote the codewords of the Golay code and the vectors in the Leech lattice. First we arrange  the set $\Omega$ into
a $4\times 6$ array such that the six columns forms a sextet. 

For each codeword in $\mathcal{G}$, $0$ and $1$ are marked by a
blanked and non-blanked space, respectively, at the corresponding
positions in the array. For example, $(1^8 0^{16})$ is denoted by
the array
\[
\begin{array}{|cc|cc|cc|}
\hline *& *&\ &\ &\ &\ \\
* & *&\ &\ &\ &\ \\
* & *&\ &\ &\ &\ \\
* & *&\ &\ &\ &\ \\ \hline
\end{array}
\]
\medskip

The following is a standard construction of the Leech lattice.
\begin{de}[Standard Leech lattice \cite{cs,Gr12}]\labtt{gen}
Let $e_{i}: =\frac{1}{\sqrt{8}}\left( 0,\dots ,4,\dots ,0\right) $ for $i\in \Omega$. Then $(e_i,e_j)=2\delta_{i,j}$. Denote $e_{X}:=\sum_{i\in X}e_{i}$ for $X\in \mathcal{G}$. The \textsl{standard Leech lattice} $\Lambda $ is a lattice of rank 24 generated by
the vectors:
\begin{eqnarray*}
&&\frac{1}{2}e_{X}\, ,\quad \text{where } X\text{ is a generator of the Golay code }%
\mathcal{G}; \\
&&\frac{1}{4}e_{\Omega }-e_{1}\,\,; \\
&&e_{i}\pm e_{j}\,,\text{ }i,j \in \Omega.
\end{eqnarray*}
\end{de}

\begin{rem}
By arranging the set $\Omega$ into a $4\times 6$ array, every vector in the Leech lattice $\Lambda$ can be written as the form
\[
X=\frac{1}{\sqrt{8}}\left[ X_1 X_2 X_3 X_4 X_5 X_6\right],\quad \text{
juxtaposition of column vectors}.
\]
For example,
\[
\frac{1}{\sqrt{8}}\,
\begin{array}{|rr|rr|rr|}
 \hline 2  & 2& 0  & 0 & 0 & 0 \\
 2 &  2 &  0&  0& 0& 0\\
 2 &  2& 0& 0&  0 & 0 \\
 2 &  2& 0&  0& 0 & 0\\ \hline
\end{array}
\]
denotes the vector $\displaystyle \frac1{2}\, e_A$, where $A$ is the
codeword 
\[
\begin{array}{|cc|cc|cc|}
\hline *& *&\ &\ &\ &\ \\
* & *&\ &\ &\ &\ \\
* & *&\ &\ &\ &\ \\
* & *&\ &\ &\ &\ \\ \hline
\end{array}\, .
\]
\end{rem}

\begin{de}
A set of vectors $\{\pm \b_1, \dots,\pm \b_{24}\} \subset \Lambda$
is called a \textsl{frame} of $\Lambda$ if $( \b_i, \b_j)= 8
\delta_{i,j}$ for all $i,j\in \{1, \dots,24\}$. For example,  $\{\pm
2e_1, \dots, \pm 2e_{24}\}$ is a frame and we call it the {\it standard
frame}.
\end{de}

Next, we shall recall some basic facts about the involutions in
$O(\Lambda)$.

Let $\mathcal{F}=\{\pm\b_1, \dots, \pm\b_{24}\}$ be a frame. For
any subset $S\subset \Omega$, we can define an isometry
$\varepsilon_S^{\mathcal{F}}: \mathbb{R}^{24} \to \mathbb{R}^{24}$
by $\varepsilon_S^{\mathcal{F}}(\b_i)= -\b_i$ if $i\in S$ and
$\varepsilon_S^{\mathcal{F}}(\b_i)= \b_i$ if $i\notin S$.
The involutions in
$O(\Lambda)$ can be characterized as follows:

\begin{thm}[\cite{cs,Gr12}]\labtt{invoinLambda}
There are exactly $4$ conjugacy classes of involutions in
$O(\Lambda)$. They correspond to the involutions
$\varepsilon_S^{\mathcal{F}}$, where $\mathcal{F}$ is a frame and
$S\in \mathcal{G}$ is an octad, the complement of an octad, a
dodecad, or the set $\Omega$. Moreover, the eigen-sublattice
$\{v\in\Lambda\mid \varepsilon_S^{\mathcal{F}}(v)=-v\}$ is
isomorphic to $EE_8$, $BW_{16}$, $DD_{12}^+$ and $\Lambda$, respectively,
where ${BW}_{16}$ is the Barnes-Wall lattice of rank $16$.
\end{thm}

\subsection{Standard $EE_8$s in the Leech lattice}
We shall describe some standard $EE_8$s in the Leech lattice in this subsection.

\subsubsection{$EE_8$ corresponding to octads in different frames}

Let $\mathcal{F}=\{\pm \b_1, \dots,\pm \b_{24}\} \subset \Lambda$ be
a frame and denote $\a_i:=\b_i/2$. For any octad $\mathcal{O}$,
denote
\[
E_\mathcal{F}(\mathcal{O})=\mathrm{span}\left\{ \a_i \pm \a_j, i,j
\in \mathcal{O}, \frac{1}2\sum_{i\in \mathcal{O}} \a_i\right\}.
\]
Then $E_\mathcal{F}(\mathcal{O})$ is a sublattice of $\Lambda$
isomorphic to $EE_8$. If $\{\pm 2e_1, \dots, \pm
2e_{24}\}$ is the standard frame, we shall simply denote
$E_\mathcal{F}(\mathcal{O})$ by $E(\mathcal{O})$.

\medskip

Next we shall consider another frame. Let
\[
A= \frac{1}2 \left [
\begin{matrix}
-1 & 1& 1& 1\\
1 & -1& 1& 1\\
1 & 1& -1& 1\\
1 & 1& 1& -1
\end{matrix}\right].
\]

\begin{nota}\labtt{defofxi}
Define a linear map $\xi:\Lambda \to \Lambda$  by  $ X\xi = AXD$, where
\[
X=\frac{1}{\sqrt{8}}\left[ X_1 X_2 X_3 X_4 X_5 X_6\right]
\]
is a vector in the Leech lattice $\Lambda$ and $D$ is the diagonal matrix
\[
\begin{pmatrix}
-1&0&0&0&0&0\\
0& 1&0&0&0&0\\
0& 0& 1&0&0&0\\
0&0&0 &1&0&0\\
0&0&0&0&1&0\\
0&0&0&0&0& 1\\
\end{pmatrix}.
\]
Recall that $\xi$ defines an isometry of $\Lambda$ (cf. \cite[p. 288]{cs} and \cite[p. 97]{Gr12}).
\end{nota}

Let $\mathcal{F}:=\{\pm 2e_1, \dots, \pm 2e_{24}\}$ be the standard frame. Then
$\mathcal{F}_\xi= \{\pm 2e_1, \dots, \pm 2e_{24}\} \xi $ is also a
frame. In this case, $E(\mathcal{O})\xi=
E_{\mathcal{F}_\xi}(\mathcal{O})$ is also isomorphic to $EE_8$ for
any octad $\mathcal{O}$. Note that if $\mathcal{F}$ is a frame and $g\in O(\Lambda)$, then $\mathcal{F}g$ is also frame.

\subsubsection{$EE_8$ associated to an even permutation in octad stabilizer}

The subgroup of $Sym_{\Omega}$ which fixes $\mathcal{G}$ setwise is the Mathieu group $M_{24}$, which is a simple group of order $2^{10}\cdot 3^3\cdot 5\cdot 7\cdot 11\cdot 23$.
Recall that $M_{24}$ is transitive on octads. The
stabilizer of an octad is the group $2^4{:}Alt_8 \cong AGL(4,2)$ and
it acts as the alternating group $Alt_8$ on the octad. If we fix a
particular point outside the octad, then every even permutation on
the octad can be extended to a unique element of $M_{24}$ which
fixes the point.

Let $\sigma=(ij)(k\ell)\in Sym(\mathcal{O})$ be a product of $2$
disjoint transpositions on the standard octad $\mathcal{O}$.
Then $\sigma$ determines a sextet which contains
$\{i,j,k,\ell\}$ as a tetrad and $\sigma$  extends uniquely
to an element $\tilde{\sigma}$ which fixes a particular point
outside the octad. Note that $\tilde{\sigma}$ fixes $2$ tetrads and
keeps the other $4$ tetrads invariant. Moreover, $\tilde{\sigma}$ has a rank $8$ $(-1)$-eigenlattice which we call $E$, and that $E$ is isometric to $EE_8$.

Take $\tilde{\sigma}$  to be the involution (UP6) listed in
\cite[pp. 49--52]{Gr12}.
\[
\scalebox{1.0} 
{
\begin{pspicture}(0,-1.0)(3.0,1.0)
\definecolor{color0c}{rgb}{0.5019607843137255,0.5019607843137255,0.5019607843137255}
\rput(0.0,-1.0){\psgrid[gridwidth=0.028222222,subgridwidth=0.014111111,gridlabels=0.0pt,subgriddiv=0,subgridcolor=color0c](0,0)(0,0)(3,2)}
\psdots[dotsize=0.12](0.24,0.76) \psdots[dotsize=0.12](0.74,0.76)
\psdots[dotsize=0.12](0.24,0.26) \psdots[dotsize=0.12](0.74,0.26)
\psdots[dotsize=0.12](1.24,0.76) \psdots[dotsize=0.12](1.74,0.74)
\psdots[dotsize=0.12](1.24,0.26) \psdots[dotsize=0.12](1.74,0.24)
\psdots[dotsize=0.12](2.24,0.74) \psdots[dotsize=0.12](2.24,0.26)
\psdots[dotsize=0.12](2.74,0.76) \psdots[dotsize=0.12](2.74,0.26)
\psdots[dotsize=0.12](0.24,-0.24) \psdots[dotsize=0.12](0.72,-0.24)
\psdots[dotsize=0.12](0.24,-0.74) \psdots[dotsize=0.12](0.74,-0.74)
\psdots[dotsize=0.12](1.24,-0.24) \psdots[dotsize=0.12](1.74,-0.24)
\psdots[dotsize=0.12](1.24,-0.74) \psdots[dotsize=0.12](1.74,-0.74)
\psdots[dotsize=0.12](2.24,-0.26) \psdots[dotsize=0.12](2.24,-0.76)
\psdots[dotsize=0.12](2.74,-0.24) \psdots[dotsize=0.12](2.72,-0.74)
\psline[linewidth=0.04cm](0.24,0.76)(0.72,0.76)
\psline[linewidth=0.04cm](0.72,-0.22)(0.74,-0.7)
\psline[linewidth=0.04cm](1.24,0.28)(1.7,0.72)
\psline[linewidth=0.04cm](1.72,-0.24)(2.22,0.24)
\psline[linewidth=0.04cm](2.22,-0.76)(2.72,-0.26)
\pscustom[linewidth=0.04] {
\newpath
\moveto(1.24,-0.22) \lineto(1.24,-0.17)
\curveto(1.24,-0.145)(1.25,-0.1)(1.26,-0.08)
\curveto(1.27,-0.06)(1.285,-0.015)(1.29,0.01)
\curveto(1.295,0.035)(1.31,0.075)(1.32,0.09)
\curveto(1.33,0.105)(1.355,0.135)(1.37,0.15)
\curveto(1.385,0.165)(1.41,0.195)(1.42,0.21)
\curveto(1.43,0.225)(1.455,0.25)(1.47,0.26)
\curveto(1.485,0.27)(1.515,0.29)(1.53,0.3)
\curveto(1.545,0.31)(1.57,0.34)(1.58,0.36)
\curveto(1.59,0.38)(1.615,0.41)(1.63,0.42)
\curveto(1.645,0.43)(1.67,0.455)(1.68,0.47)
\curveto(1.69,0.485)(1.715,0.51)(1.73,0.52)
\curveto(1.745,0.53)(1.78,0.55)(1.8,0.56)
\curveto(1.82,0.57)(1.855,0.59)(1.87,0.6)
\curveto(1.885,0.61)(1.915,0.63)(1.93,0.64)
\curveto(1.945,0.65)(1.975,0.67)(1.99,0.68)
\curveto(2.005,0.69)(2.045,0.71)(2.07,0.72)
\curveto(2.095,0.73)(2.145,0.74)(2.17,0.74)
\curveto(2.195,0.74)(2.225,0.74)(2.24,0.74) }
\pscustom[linewidth=0.04] {
\newpath
\moveto(1.24,-0.68) \lineto(1.25,-0.64)
\curveto(1.255,-0.62)(1.265,-0.58)(1.27,-0.56)
\curveto(1.275,-0.54)(1.29,-0.505)(1.3,-0.49)
\curveto(1.31,-0.475)(1.335,-0.44)(1.35,-0.42)
\curveto(1.365,-0.4)(1.385,-0.36)(1.39,-0.34)
\curveto(1.395,-0.32)(1.415,-0.29)(1.43,-0.28)
\curveto(1.445,-0.27)(1.475,-0.24)(1.49,-0.22)
\curveto(1.505,-0.2)(1.53,-0.165)(1.54,-0.15)
\curveto(1.55,-0.135)(1.57,-0.105)(1.58,-0.09)
\curveto(1.59,-0.075)(1.61,-0.045)(1.62,-0.03)
\curveto(1.63,-0.015)(1.655,0.01)(1.67,0.02)
\curveto(1.685,0.03)(1.715,0.055)(1.73,0.07)
\curveto(1.745,0.085)(1.775,0.115)(1.79,0.13)
\curveto(1.805,0.145)(1.83,0.175)(1.84,0.19)
\curveto(1.85,0.205)(1.875,0.235)(1.89,0.25)
\curveto(1.905,0.265)(1.935,0.295)(1.95,0.31)
\curveto(1.965,0.325)(1.995,0.355)(2.01,0.37)
\curveto(2.025,0.385)(2.06,0.41)(2.08,0.42)
\curveto(2.1,0.43)(2.135,0.45)(2.15,0.46)
\curveto(2.165,0.47)(2.2,0.49)(2.22,0.5)
\curveto(2.24,0.51)(2.28,0.535)(2.3,0.55)
\curveto(2.32,0.565)(2.365,0.585)(2.39,0.59)
\curveto(2.415,0.595)(2.455,0.615)(2.47,0.63)
\curveto(2.485,0.645)(2.515,0.675)(2.53,0.69)
\curveto(2.545,0.705)(2.585,0.725)(2.61,0.73)
\curveto(2.635,0.735)(2.68,0.745)(2.7,0.75) }
\pscustom[linewidth=0.04] {
\newpath
\moveto(1.74,-0.7) \lineto(1.76,-0.66)
\curveto(1.77,-0.64)(1.79,-0.6)(1.8,-0.58)
\curveto(1.81,-0.56)(1.845,-0.515)(1.87,-0.49)
\curveto(1.895,-0.465)(1.935,-0.425)(1.95,-0.41)
\curveto(1.965,-0.395)(1.985,-0.36)(1.99,-0.34)
\curveto(1.995,-0.32)(2.01,-0.285)(2.02,-0.27)
\curveto(2.03,-0.255)(2.05,-0.225)(2.06,-0.21)
\curveto(2.07,-0.195)(2.095,-0.165)(2.11,-0.15)
\curveto(2.125,-0.135)(2.15,-0.105)(2.16,-0.09)
\curveto(2.17,-0.075)(2.2,-0.05)(2.22,-0.04)
\curveto(2.24,-0.03)(2.275,-0.01)(2.29,0.0)
\curveto(2.305,0.01)(2.335,0.03)(2.35,0.04)
\curveto(2.365,0.05)(2.395,0.07)(2.41,0.08)
\curveto(2.425,0.09)(2.455,0.115)(2.47,0.13)
\curveto(2.485,0.145)(2.52,0.175)(2.54,0.19)
\curveto(2.56,0.205)(2.605,0.225)(2.63,0.23)
\curveto(2.655,0.235)(2.685,0.24)(2.7,0.24) }
\usefont{T1}{ptm}{m}{n}
\rput(5.0,-0.05){\small (UP 6)}
\end{pspicture}
}  
\]
Then $\tilde{\sigma}$ stabilizes the octad
\[
\begin{array}{|cc|cc|cc|}
\hline  * & * &  &  &  & \\
* & * &\ &\ &\ &\ \\
* & * &\ &\ &\ &\ \\
* & * &\ &\ &\ &\ \\ \hline
\end{array}
\]
and  determines as above the sublattice
\[
E=\mathrm{span}_\ZZ\left \{ \pm\a_i\pm \a_j, \frac{1}2
\sum_{i=1}^{8} \epsilon_i \a_i\right \},
\]
where

\begin{align*}
&\a_1=\frac{1}{\sqrt{8}}\,
\begin{array}{|cc|cc|cc|}
\hline\ 2 & -2&\ 0  & \ 0 &\ 0 &\ 0 \\
\ 0 & \ 0&\ 0&\ 0&\ 0&\ 0\\
\ 0 &  \ 2&\ 0&\ 0&\ 0 &\ 0 \\
\ 0 & -2&\ 0&\ 0&\ 0 &\ 0 \\ \hline
\end{array}, &&
 \a_2=\frac{1}{\sqrt{8}}\,
\begin{array}{|cc|cc|cc|}
\hline -2 &\ 2&\ 0  &\ 0 &\ 0 &\ 0 \\
\ 0 & \ 0&\ 0&\ 0&\ 0& 0\\
\ 0 &\ 2&\ 0&\ 0&\ 0 &\ 0 \\
\ 0 &  -2&\ 0&\ 0&\ 0 &\ 0 \\ \hline
\end{array},\\
&\a_3=\frac{1}{\sqrt{8}}\,
\begin{array}{|cc|cc|cc|}
\hline \ 0 &\  0&\ 0  &\ 2 & \ 0 &\ 0 \\
\ 0 &\ 0& -2&\ 0&\ 0&\ 0\\
\ 0 &\ 0&\ 0&\ 0&\ 0 & -2\\
\ 0 &\ 0&\ 0&\ 0&\ 2 &\ 0 \\ \hline
\end{array}, &
&\a_4=\frac{1}{\sqrt{8}}\,
\begin{array}{|cc|cc|cc|}
\hline \ 0 &\ 0&\ 0  & -2 &\ 0 &\ 0 \\
\ 0 &\ 0&\ 2&\ 0&\ 0&\ 0\\
\ 0 &\ 0&\ 0&\ 0&\ 0 & -2\\
 \ 0 &\ 0&\ 0&\ 0&\ 2 &\ 0 \\ \hline
\end{array},
\end{align*}
\begin{align*}
&\a_5=\frac{1}{\sqrt{8}}\,
\begin{array}{|cc|cc|cc|}
\hline \ 0 &\ 0&\ 0  &\ 0& - 2 &\ 0 \\
\ 0 &\ 0&\ 0&\ 0&\ 0& 2\\
\ 0 &\ 0& \ 2&\ 0&\ 0 &\ 0\\
\ 0 &\ 0&\ 0& -2&\ 0 &\ 0 \\ \hline
\end{array}, &
&\a_6=\frac{1}{\sqrt{8}}\,
\begin{array}{|cc|cc|cc|}
\hline\ 0 &\ 0&\ 0  &\ 0&\ 2 &\ 0 \\
\ 0 &\ 0&\ 0&\ 0&\ 0&\ 2\\
\ 0 &\ 0& -2&\ 0&\ 0 &\ 0\\
\ 0 &\ 0& \ 0& -2&\ 0 &\ 0 \\ \hline
\end{array},\\
&\a_7=\frac{1}{\sqrt{8}}\,
\begin{array}{|cc|cc|cc|}
\hline \ 0 &\ 0&\ 0  &\ 0& \ 0 & 2 \\
\ 0 &\ 0&\ 0&\ 0& -2& \ 0\\
\ 0 &\ 0& \ 0&\ 2&\ 0 &\ 0\\
\ 0 &\ 0& -2& \ 0&\ 0 &\ 0 \\ \hline
\end{array}, &
&\a_8=\frac{1}{\sqrt{8}}\,
\begin{array}{|cc|cc|cc|}
\hline \ 0 &\ 0&\ 0  &\ 0& \ 0 & -2 \\
\ 0 &\ 0&\ 0&\ 0& -2& \ 0\\
\ 0 &\ 0& \ 0&\ 2&\ 0 &\ 0\\
\ 0 &\ 0&\ 2& \ 0&\ 0 &\ 0 \\ \hline
\end{array},
\end{align*}
$i,j=1,\dots,8$ and $\epsilon_i=\pm 1$ such that $\prod_{i=1}^8
\epsilon_i=1$.  Then $E$ is a sublattice in $\Lambda$ which is
isomorphic to $EE_8$. By our construction, it is also clear
that $\tilde{\sigma}$ acts as $-1$ on $E$ and $1$ on $ann_\Lambda(E)$.

Recall that $\tilde{\sigma}$ is acting on $\Lambda$ from the right according to  our convention.

\subsection{$EE_8$ pairs in the Leech lattice}\label{sec:F3}
In this subsection, we shall describe certain NREE8 pairs $M,N$ explicitly inside the Leech lattice. By using the uniqueness theorem (cf. Theorem \ref{uniquenessofl}), we know that our examples are actually isomorphic to the lattices in Table \ref{NREE8SUM}. It turns out that
except for $\dih{4}{15}$, all lattices in Table \ref{NREE8SUM} can be embedded into the Leech lattice.

We shall note that the lattice $L=M+N$ is uniquely determined (up to isometry) by the rank of $L$ and the order of the dihedral group $D:=\la t_M, t_N\ra$ except for $\dih{8}{16,0}$ and $\dih{8}{16, DD_4}$. Extra information about $ann_M(N)$ is needed to distinguish them.

\medskip

Let $M$ and $N$ be $EE_8$ sublattices of the Leech lattice $\Lambda$. Let
$t:=t_M$ and $u:=t_N$ be the involutions of $\Lambda$ such that $t$
and $u$ act on $M$ and $N$ as $-1$ and act as $1$ on $M^\perp$ and
$N^\perp$, respectively. Set $g:=ut$ and $D:=\la t,u\ra $, the dihedral
group generated by $t$ and $u$.

\medskip

\begin{nota}\labtt{octad}
In this subsection, $\mathcal{O}$, $\mathcal{O}'$, $\mathcal{O}''$, etc denote some arbitrary octads while $\mathcal{O}_1$, $\mathcal{O}_2$, $\mathcal{O}_3$, and  $\mathcal{O}_4$ denote the fixed octads given as follows.
\[
\begin{split}
\mathcal{O}_1=
\begin{array}{|cc|cc|cc|}
\hline  * & * &  &  &  & \\
* & * &\ &\ &\ &\ \\
* & * &\ &\ &\ &\ \\
* & * &\ &\ &\ &\ \\ \hline
\end{array}\ , \qquad
&
\mathcal{O}_2=\begin{array}{|cc|cc|cc|}
\hline \ & *&* &* &* &* \\
* &\ &\ &\ &\ &\ \\
* &\ &\ &\ &\ &\ \\
* &\ &\ &\ &\ &\ \\ \hline
\end{array}\ ,\\
\mathcal{O}_3= \begin{array}{|cc|cc|cc|}
\hline\  &\ &\ &\ & * & * \\
 \ &\ &\ &\ & * & * \\
 \ &\ &\ &\ &* & * \\
 \ &\  &\ &\ & * & * \\ \hline
\end{array}\ , \qquad
&
\mathcal{O}_4=  \begin{array}{|cc|cc|cc|}
\hline \ && *&* &* &* \\
 &\ & *&* &* &* \\
 &\ & &\ &\ &\ \\
 &\ & &\ &\ &\ \\ \hline
\end{array}.
\end{split}
\]
\end{nota}

\begin{rem}
All the Gram matrices in this subsection are computed by multiplying the matrix $A$ by its transpose $A^t$, where $A$ is the matrix whose rows form an ordered basis given in each case. The Smith invariants sequences are computed using the command \textit{ismith} in Maple 8.
\end{rem}

\subsubsection{$|g|=2$.} In this case, $M\cap
N \cong 0, AA_1, AA_1\perp AA_1$ or $DD_4$.

\medskip

\noindent {\bf Case: $\dih{4}{15}$}: This case does not embed into
$\Lambda$.

If $M\cap N\cong AA_1$, then $L=M+N \cong \dih{4}{15}$ contains a
sublattice isometric to $AA_1\perp EE_7\perp EE_7$, which cannot be
embedded in the Leech lattice $\Lambda$ because the $(-1)$-eigenlattice of the involution $g=t_Mt_N$ has rank $14$ but there is no such involution in $O(\Lambda)$ (cf. Theorem \ref{invoinLambda}).

\medskip

\begin{nota}\labtt{octad2}
Let $\mathcal{O}=\{i_1, \dots, i_8\}$ and $\mathcal{O}'=\{j_1, \dots, j_8\}$ be $2$ distinct octads and denote $M := E(\mathcal{O})$ and $N:= E(\mathcal{O}')$.  Since the Golay code $\mathcal{G}$ is a type $II$ code (doubly even) and the minimal
norm of $\mathcal{G}$ is $8$, $|\mathcal{O}\cap \mathcal{O}'|$ is
either $0$, $2$, or $4$.
\end{nota}
\medskip

\noindent $\bf \dih{4}{16}$

When $|\mathcal{O}\cap \mathcal{O}'|=0$, clearly $M\cap N=0$ and $M+N\cong EE_8\perp EE_8$.

\medskip

\noindent $\bf \dih{4}{14}$

Suppose $\mathcal{O}\cap \mathcal{O}'=\{i_1, i_2\}=\{j_1, j_2\}$. Then $|\mathcal{O}\cap \mathcal{O}'|=2$ and  $F=M\cap N=\mspan_\ZZ\{e_{i_1}+e_{i_2}, e_{i_1}-e_{i_2}\} \cong AA_1\perp AA_1$. In this case, $ann_{M}(F)\cong ann_N(F)\cong DD_6$ and $L$ contains a sublattice of type $AA_1\perp AA_1\perp DD_6\perp DD_6$ which has index $2^4$ in $L$. Note that $L$ is of rank $14$. By computing the Gram matrices, it is easy to check that
\[
\begin{split}
& \{\frac{1}2(e_{i_1}+\cdots+e_{i_8})\}\cup   \{ -e_{i_{k}}+ e_{i_{k-1}}|\, 7 \geq  k \geq 3\}\cup \{-e_{i_2}+e_{i_1}, -e_{i_1}-e_{i_2}\}
\end{split}
\]
is a basis of $M$ and
\[
\begin{split}
\{-e_{i_2}+e_{i_1}, -e_{i_1}-e_{i_2}\}\cup \{  e_{j_{k-1}}- e_{j_k}|\, 3\leq  k \leq 7\}\cup\{\frac{1}2(e_{j_1}+\cdots+e_{j_8})\}
\end{split}
\]
is a basis of $N$. Thus,
\[
\begin{split}
& \{\frac{1}2(e_{i_1}+\cdots+e_{i_8})\}\cup   \{ -e_{i_{k}}+ e_{i_{k-1}}|\, 8 \geq  k \geq 3\}\cup \{-e_{i_2}+e_{i_1}, -e_{i_1}-e_{i_2}\}\\
& \cup \{  e_{j_{k-1}}- e_{j_k}|\, 3\leq  k \leq 8\}\cup\{\frac{1}2(e_{j_1}+\cdots+e_{j_8})\}
\end{split}
\]
is a basis of $L$ and the Gram matrix of $L$ is given by

\[
\left[
\begin{array}{cccccccccccccc}
 4&0&0&0&0&0&0&-2&1&0&0&0&0&1\\
 0&4&-2&0&0&0&0&0&0&0&0&0&0&0\\
 0&-2&4&-2&0&0&0&0&0&0&0&0&0&0\\
 0&0&-2&4&-2&0&0&0&0&0&0&0&0&0\\
 0&0&0&-2&4&-2&0&0&0&0&0&0&0&0\\
 0&0&0&0&-2&4&-2&-2&2&0&0&0&0&1\\
 0&0&0&0&0&-2&4&0&-2&0&0&0&0&0\\
 -2&0&0&0&0&-2&0&4&-2&0&0&0&0&-2\\
 1&0&0&0&0&2&-2&-2&4&-2&0&0&0&0\\
 0&0&0&0&0&0&0&0&-2&4&-2&0&0&0\\
 0&0&0&0&0&0&0&0&0&-2&4&-2&0&0\\
 0&0&0&0&0&0&0&0&0&0&-2&4&-2&0\\
 0&0&0&0&0&0&0&0&0&0&0&-2&4&0\\
 1&0&0&0&0&1&0&-2&0&0&0&0&0&4
 \end{array}
 \right] .
 \]
The Smith invariant sequence is $ 1111 2222 2222 44$.

\medskip
\noindent $\bf \dih{4}{12}$

Suppose $\mathcal{O}\cap \mathcal{O}'=\{i_1,i_2,i_3,i_4\}=\{j_1,j_2,j_3,j_4\}$ (cf. Notation \ref{octad2}). Then $|\mathcal{O}\cap \mathcal{O}'|=4$ and $F= M\cap N= \mspan_\ZZ\{e_{i_k}\pm e_{i_l}|\, 1\leq k< l \leq 4  \}\cong
DD_4$. Thus, $ann_M(F)\cong ann_N(F)\cong DD_4$. In this case, $L$ is of rank $12$ and it contains a sublattice of type $DD_4\perp DD_4\perp DD_4$ which has index $2^4$ in $L$. Note that
$
\{e_{i_1}+ e_{i_2}, e_{i_1}- e_{i_2}, e_{i_2}- e_{i_3}, e_{i_3}- e_{i_4}\}
$
is a basis of $F=M\cap N$. A check of Gram matrices also shows that
$$
\{e_{i_1}+ e_{i_2}, e_{i_1}- e_{i_2}, e_{i_2}- e_{i_3}, e_{i_3}- e_{i_4}\}
\cup \{e_{i_4}- e_{i_5}, e_{i_5}- e_{i_6}, e_{i_6}- e_{i_7}, \frac{-1}2(e_{i_1}+\cdots+e_{i_8}) \}
$$
is a basis of $M$ and
$$
\{e_{i_1}+ e_{i_2}, e_{i_1}- e_{i_2}, e_{i_2}- e_{i_3}, e_{i_3}- e_{i_4}\}  \cup \{e_{j_4}- e_{j_5}, e_{j_5}- e_{j_6}, e_{j_6}- e_{j_7}, \frac{-1}2(e_{j_1}+\cdots+e_{j_8})\}
$$
is a basis of $N$. Therefore, $L=M+N$ has a basis
\[
\begin{split}
&\{e_{i_1}+ e_{i_2}, e_{i_1}- e_{i_2}, e_{i_2}- e_{i_3}, e_{i_3}- e_{i_4}\}\\
\cup &\{e_{i_4}- e_{i_5}, e_{i_5}- e_{i_6}, e_{i_6}- e_{i_7}, \frac{-1}2(e_{i_1}+\cdots+e_{i_8}) \}\\
 \cup &\{e_{j_4}- e_{j_5}, e_{j_5}- e_{j_6}, e_{j_6}- e_{j_7}, \frac{-1}2(e_{j_1}+\cdots+e_{j_8})\}
\end{split}
\]
and the Gram matrix of $L$ is given by

\[
\left[
\begin{array}{cccccccccccc}
 4&0&2&0&0&0&0&-2&0&0&0&-2\\
 0&4&-2&0&0&0&0&0&0&0&0&0\\
 2&-2&4&-2&0&0&0&0&0&0&0&0\\
 0&0&-2&4&-2&0&0&0&-2&0&0&0\\
 0&0&0&-2&4&-2&0&0&2&0&0&-1\\
 0&0&0&0&-2&4&-2&0&0&0&0&0\\
 0&0&0&0&0&-2&4&0&0&0&0&0\\
 -2&0&0&0&0&0&0&4&-1&0&0&2\\
 0&0&0&-2&2&0&0&-1&4&-2&0&0\\
 0&0&0&0&0&0&0&0&-2&4&-2&0\\
 0&0&0&0&0&0&0&0&0&-2&4&0\\
 -2&0&0&0&-1&0&0&2&0&0&0&4
 \end{array}
 \right]
\]
whose Smith invariant sequence is $1111 2222 22 44$.

\medskip

\subsubsection{$ |g|=3$.}\labtt{g=3} In this case, $M\cap N =0$ or
$AA_2$.

\medskip
\noindent $\bf \dih{6}{16}$

\begin{nota}\labtt{basisofEO1}
Let $M:=E(\mathcal{O}_1)\cong EE_8$, where $\mathcal{O}_1$ is the
octad  described in Notation \refpp{octad}.  We choose a basis
$\{\b_1,\dots, \b_8\}$ of $M$, where

\begin{align*}
\b_1= \frac{1}{\sqrt{8}}\,
\begin{array}{|cc|cc|cc|}
\hline 4 & -4& 0  & 0 & 0 & 0 \\
0 &  0&  0&  0& 0& 0\\
 0 & 0&  0& 0&  0 & 0 \\
 0 &  0& 0&  0& 0 & 0\\ \hline
\end{array}, &&
 \b_2= \frac{1}{\sqrt{8}}\,
 \begin{array}{|rr|rr|rr|}
\hline 0  & 4& 0  & 0 & 0 & 0 \\
      -4 &  0&  0&  0& 0& 0\\
       0 & 0&  0& 0&  0 & 0 \\
       0 &  0& 0&  0& 0 & 0\\ \hline
\end{array}, \\
\b_3= \frac{1}{\sqrt{8}}\,
 \begin{array}{|rr|rr|rr|}
\hline 0  & 0& 0  & 0 & 0 & 0 \\
       4 &  -4&  0&  0& 0& 0\\
       0 & 0&  0& 0&  0 & 0 \\
       0 &  0& 0&  0& 0 & 0\\ \hline
\end{array}, &&
\b_4=\frac{1}{\sqrt{8}}\,
\begin{array}{|rr|rr|rr|}
\hline 0  & 0& 0  & 0 & 0 & 0 \\
       0&  4 &  0&  0& 0& 0\\
       -4 & 0&  0& 0&  0 & 0 \\
       0 &  0& 0&  0& 0 & 0\\ \hline
\end{array},
\\
\b_5=\frac{1}{\sqrt{8}}\,
\begin{array}{|rr|rr|rr|}
\hline 0  & 0& 0  & 0 & 0 & 0 \\
       0&  0 &  0&  0& 0& 0\\
       4 & -4&  0& 0&  0 & 0 \\
       0 &  0& 0&  0& 0 & 0\\ \hline
\end{array},  &&
\b_6=\frac{1}{\sqrt{8}}\,
\begin{array}{|rr|rr|rr|}
\hline 0  & 0& 0  & 0 & 0 & 0 \\
       0 & 0&  0&  0& 0& 0\\
       0 & 4&  0& 0&  0 & 0 \\
       -4&  0& 0&  0& 0 & 0\\ \hline
\end{array}, \\
\b_7=\frac{1}{\sqrt{8}}\,
\begin{array}{|rr|rr|rr|}
\hline -4  & -4& 0  & 0 & 0 & 0 \\
       0&  0 &  0&  0& 0& 0\\
       0 & 0&  0& 0&  0 & 0 \\
       0 &  0& 0&  0& 0 & 0\\ \hline
\end{array}, &&
\b_8=\frac{1}{\sqrt{8}}\,
\begin{array}{|rr|rr|rr|}
\hline\ 2  &\ 2& 0  & 0 & 0 & 0 \\
       2 &  2 &  0&  0& 0& 0\\
       2 &  2& 0& 0&  0 & 0 \\
       2 &  2& 0&  0& 0 & 0\\ \hline
\end{array}.
\end{align*}
\end{nota}

Let $N$ be the lattice generated by the vectors
\begin{align*}
&\a_1=\frac{1}{\sqrt{8}}\,
\begin{array}{|rr|rr|rr|}
\hline\ \ 2 & -2&\ \ 2  & -2 &\ \ 2 &\ \, 0 \\
\ 0 & \ 0&\ 0&\ 0&\ 0&  2\\
\ 0 &  \ 0&\ 0&\ 0&\ 0 & 2 \\
\ 0 & \ 0&\ 0&\ 0&\ 0 &\ \, 2\\ \hline
\end{array}, &&
 \a_2=\frac{1}{\sqrt{8}}\,
\begin{array}{|rr|rr|rr|}
\hline 0 &\ 2&\ 0&\ 2&  0 &\ 2 \\
      -2 &\ 0& -2&\ 0&\ 0 &\ -2 \\
     \ 0 &\ 0&\ 0&\ 0&\ 2 &\ -2\\
     \ 0 &\ 0&\ 0&\ 0&\ 0 &\ 0 \\ \hline
\end{array},\\
&\a_3=\frac{1}{\sqrt{8}}\,
\begin{array}{|rr|rr|rr|}
\hline \ 0 &\  0&\ 0  &\ 0 &\ 0 & -2 \\
     \ \ 2 & -2 & \ \,2 & -2 & 2  &\ 0\\
       \ 0 &\ 0 &\ 0  &\ 0 & -2  &  0\\
       \ 0 &\ 0 &\ 0  &\ 0 & 2  &\ 0 \\ \hline
\end{array}, &
&\a_4=\frac{1}{\sqrt{8}}\,
\begin{array}{|rr|rr|rr|}
\hline  \ 0 &\ 0  &   \ 0&    0&\ 0 &\ 2 \\
        \ 0 &\ \ 2&   \ 0&\ \ 2&  -2&\ \ 2\\
         -2 &\   0&    -2&\   0&\ 0 & \ 0\\
        \ 0 &  \ 0&   \ 0& \  0&\ 0 & -2 \\ \hline
\end{array},\\
&\a_5=\frac{1}{\sqrt{8}}\,
\begin{array}{|rr|rr|rr|}
\hline   \ 0 &\ 0&  \ 0&\ 0& \ 0 &-2 \\
         \ 0 &\ 0&  \ 0&\ 0&\ 2& 0\\
       \ \ 2 & -2&\ \ 2& -2&\ 2 &0\\
         \ 0 &\ 0&\ 0&  0& -2 & 0 \\ \hline
\end{array}, &
&\a_6=\frac{1}{\sqrt{8}}\,
\begin{array}{|rr|rr|rr|}
\hline \ 0 &  \ 0&\ 0 &  \ 0& 0   &\ 2 \\
       \ 0 &  \ 0& \ 0&  \ 0& 0   & -2\\
       \ 0 &\ \ 2&   0&\ \ 2&\ -2   &\ 2\\
        -2 &  \ 0&  -2&    0&\ 0 &\ \ 0 \\ \hline
\end{array},\\
&\a_7=\frac{1}{\sqrt{8}}\,
\begin{array}{|rr|rr|rr|}
\hline     -2 &  -2&  -2& -2& \ \ 0 & -2 \\
           \ 0 &\ 0& \ 0&\ 0& -2& \ 0\\
           \ 0 &\ 0& \ 0&\ 0& -2 &\ 0\\
           \ 0 &\ 0&   0&\ 0& -2 &\ 0 \\ \hline
\end{array}, &
&\a_8=\frac{1}{\sqrt{8}}\,
\begin{array}{|rr|rr|rr|}
\hline \ \ 1 &\ \ \, 1&\ \ 1  &\ \ \,1& \ -3 & 1 \\
       \ \ 1 &\ \ 1&\ \ 1  &\ \ 1& \  \ 1 & 1 \\
       \ \ 1 &\ \ 1&\ \ 1  &\ \ 1& \  \ 1 & 1 \\
       \ \ 1 &\ \ 1&\ \ 1  &\ \ 1& \  \ 1 & 1 \\ \hline
\end{array}.
\end{align*}
By checking the inner products, it is easy to shows that $N\cong EE_8$. Note that $\a_1, \dots, \a_7 $ are supported on octads and thus $N\leq \Lambda$ by \refpp{gen}.

In this case, $M\cap N=0$.  Then $\{\b_1, \b_2, \dots, \b_8, \a_1,\dots, \a_8\}$ is a basis of $L=M+N$ and the Gram matrix of $L=M+N$ is given by

\[
\left[
\begin{array}{rrrr rrrr| rrrr rrrr}
4&-2&0&0&0&0&0&0&2&-1&0&0&0&0&0&0\\
-2&4&-2&0&0&0&-2&0&-1&2&-1&0&0&0&-1&0\\
0&-2&4&-2&0&0&0&0&0&-1&2&-1&0&0&0&0\\
0&0&-2&4&-2&0&0&0&0&0&-1&2&-1&0&0&0\\
0&0&0&-2&4&-2&0&0&0&0&0&-1&2&-1&0&0\\
0&0&0&0&-2&4&0&0&0&0&0&0&-1&2&0&0\\
0&-2&0&0&0&0&4&-2&0&-1&0&0&0&0&2&-1\\
0&0&0&0&0&0&-2&4&0&0&0&0&0&0&-1&2\\ \hline
2&-1&0&0&0&0&0&0&4&-2&0&0&0&0&0&0\\
-1&2&-1&0&0&0&-1&0&-2&4&-2&0&0&0&-2&0\\
0&-1&2&-1&0&0&0&0&0&-2&4&-2&0&0&0&0\\
0&0&-1&2&-1&0&0&0&0&0&-2&4&-2&0&0&0\\
0&0&0&-1&2&-1&0&0&0&0&0&-2&4&-2&0&0\\
0&0&0&0&-1&2&0&0&0&0&0&0&-2&4&0&0\\
0&-1&0&0&0&0&2&-1&0&-2&0&0&0&0&4&-2\\
0&0&0&0&0&0&-1&2&0&0&0&0&0&0&-2&4
\end{array}
\right] \]

\medskip

By looking at the Gram matrix, it is clear that $L=M+N\cong A_2\otimes E_8$. The Smith invariant
sequence is $ 1111 1111 3333 3333$.

\noindent $\bf \dih{6}{12}$

Let $M:=E(\mathcal{O}_2)$ and $N:=M\xi $, where  $\mathcal{O}_2$ is
the  octad described in Notation \refpp{octad} and $\xi$ is the
isometry defined in Notation \refpp{defofxi}.

\begin{nota}\labtt{baisiofEO2}
Set
\begin{align*}
\g_1= \frac{1}{\sqrt{8}}\,
\begin{array}{|rr|rr|rr|}
\hline 0 & 0& 0  & 0 & 0 & 0 \\
 0 &  0&  0&  0& 0& 0\\
 -4 & 0&  0& 0&  0 & 0 \\
 4 &  0& 0&  0& 0 & 0\\ \hline
\end{array}, &&
 \g_2= \frac{1}{\sqrt{8}}\,
 \begin{array}{|rr|rr|rr|}
\hline 0  & 0& 0  & 0 & 0 & 0 \\
      -4 &  0&  0&  0& 0& 0\\
       4 & 0&  0& 0&  0 & 0 \\
       0 &  0& 0&  0& 0 & 0\\ \hline
\end{array}, \\
\g_3= \frac{1}{\sqrt{8}}\,
 \begin{array}{|rr|rr|rr|}
\hline 0  & -4& 0  & 0 & 0 & 0 \\
       4 &  0&  0&  0& 0& 0\\
       0 & 0&  0& 0&  0 & 0 \\
       0 &  0& 0&  0& 0 & 0\\ \hline
\end{array}, &&
\g_4=\frac{1}{\sqrt{8}}\,
\begin{array}{|rr|rr|rr|}
\hline 0  & 4& -4  & 0 & 0 & 0 \\
       0 &  0 &  0&  0& 0& 0\\
       0 & 0&  0& 0&  0 & 0 \\
       0 &  0& 0&  0& 0 & 0\\ \hline
\end{array},
\end{align*}
\begin{align*}
\g_5=\frac{1}{\sqrt{8}}\,
\begin{array}{|rr|rr|rr|}
\hline 0  & 0& 4  & -4 & 0 & 0 \\
       0 & 0&  0&  0& 0& 0\\
       0 & 0&  0& 0&  0 & 0 \\
       0 &  0& 0&  0& 0 & 0\\ \hline
\end{array}, &&
\g_6=\frac{1}{\sqrt{8}}\,
\begin{array}{|rr|rr|rr|}
\hline 0  & 0& 0  & 4 & -4 & 0 \\
       0&  0 &  0&  0& 0& 0\\
       0 & 0&  0& 0&  0 & 0 \\
       0 &  0& 0&  0& 0 & 0\\ \hline
\end{array}, \\
\g_7=\frac{1}{\sqrt{8}}\,
\begin{array}{|rr|rr|rr|}
\hline 0  & 0& 0  & 0 & 0 & 0 \\
       0&  0 &  0&  0& 0& 0\\
       -4 & 0&  0& 0&  0 & 0 \\
       -4 &  0& 0&  0& 0 & 0\\ \hline
\end{array}, &&
\g_8=\frac{1}{\sqrt{8}}\,
\begin{array}{|rr|rr|rr|}
\hline\ 0  &\ 2& 2  & 2 & 2 & 2 \\
       2 &  0 &  0&  0& 0& 0\\
       2 &  0& 0& 0&  0 & 0 \\
       2 &  0& 0&  0& 0 & 0\\ \hline
\end{array}.
\end{align*}
Then $\{\g_1, \dots,\g_8\}$ is a basis of $M$ and $\{\g_1\xi ,
\dots,\g_8\xi \}$ is a basis of $N$.
\end{nota}

By the definition of $\xi$, it is easy to show that $\g_1\xi=-\g_1,
\g_2\xi=-\g_2$. Moreover, for any $\a\in M=E(\mathcal{O}_2)$, $\a
\xi$ is supported on $\mathcal{O}_2$ if and only if $\a\in
\mspan_\ZZ\{ \g_1, \g_2\}$.  Hence, $F=M\cap N=\mspan_\ZZ\{ \g_1,
\g_2\} \cong AA_2$. Then $ann_M(F)\cong ann_N(F)\cong EE_6$ and
$L=M+N$ is of rank $14$.

Note that
$\{\g_1, \g_2, \g_3, \dots, \g_8\}$ is a basis of $M$ and $\{\g_1, \g_2, \g_3\xi, \dots, \g_8\xi \} $ is a basis of $N$. Therefore,
$$\{\g_1, \g_2\} \cup \{ \g_3, \dots, \g_8\}\cup \{\g_3\xi, \dots, \g_8\xi \} $$
is a basis of $L$ and the Gram matrix of $L$ is given by
\[
\left[
\begin{array}{cccccccccccccc}
  4&-2&0&0&0&0&0&0&0&0&0&0&0&0\\
  -2&4&-2&0&0&0&-2&0&-2&0&0&0&2&0\\
  0&-2&4&-2&0&0&0&0&0&1&0&0&-2&0\\
  0&0&-2&4&-2&0&0&0&1&-2&1&0&0&0\\
  0&0&0&-2&4&-2&0&0&0&1&-2&1&0&0\\
  0&0&0&0&-2&4&0&0&0&0&1&-2&0&0\\
  0&-2&0&0&0&0&4&-2&2&0&0&0&0&1\\
  0&0&0&0&0&0&-2&4&0&0&0&0&-1&-2\\
  0&-2&0&1&0&0&2&0&4&-2&0&0&0&0\\
  0&0&1&-2&1&0&0&0&-2&4&-2&0&0&0\\
  0&0&0&1&-2&1&0&0&0&-2&4&-2&0&0\\
  0&0&0&0&1&-2&0&0&0&0&-2&4&0&0\\
  0&2&-2&0&0&0&0&-1&0&0&0&0&4&2\\
  0&0&0&0&0&0&1&-2&0&0&0&0&2&4
  \end{array}
  \right]
\]
whose Smith invariant sequence is $1111\ 11111 \ 333 66$.

Recall that  $ann_M(F)+ann_N(F) \cong A_2\otimes E_6$ (cf. \refpp{tensorwitha2}) and thus $L$ contains a sublattice isometric to $AA_2 \perp (A_2\otimes E_6)$.

\medskip

\subsubsection{$|g|=4$.}
In this case, $M\cap N =0$ or
$AA_1$.  There are 2 subcases for $M\cap N=0$.

\medskip
\noindent $\bf \dih{8}{16, 0}$

Let $M:=E(\mathcal{O}_1)$, where $\mathcal{O}_1$ is the octad
as described in Notation \refpp{octad}.

Take $\{\b_1,\dots,  \b_8\}$ as defined in Notation
\refpp{basisofEO1}. Then it is a basis of $M=E(\mathcal{O}_1)$. Let
$N$ be the $EE_8$ sublattice generated by
\[
\a_1=\frac{1}{\sqrt{8}}
\begin{array}{|rr|rr|rr|}
\hline
     0&   0 &     0 &    0&   \ 0&   \ 0  \\
    -2&  -2 &    -2 &   -2&\ \, 0&\ \, 0\\
 \ \,0&\ \,0 &\ \, 0 &\ \, 0&\ \, 0&\ \, 0\\
    -2&  -2 &    -2 &   -2&\ \, 0&\ \, 0\\
  \hline
\end{array}
\qquad 
\a_2=\frac{1}{\sqrt{8}}
\begin{array}{|rr|rr|rr|}
\hline 0& \ \,0  &0 &\ 0 &\ 0 &\ 0  \\
 \ \ 4 &\ \ 0 &\ \, 4 &\ \, 0&\ \, 0&\ \, 0\\
 \ \,0&\ \, 0 &\ \,0&\ \,0&\ \,0&\ \,0\\
\ \,0&\ \,0 &\ \, 0 &\ \,0&\ \,0&\ \,0\\ \hline
\end{array}
\]
\[
\a_3=\frac{1}{\sqrt{8}}
\begin{array}{|rr|rr|rr|}
\hline -4& \ \ 0  &-4 &\ \ 0 &\ \ 0 &\ \,0  \\
 \ \,0 &\ \,0 &\ \, 0 &\ \, 0&\ \, 0&\ \, 0\\
 \ \,0&\ \, 0 &\ \,0&\ \,0&\ \,0&\ \,0\\
\ \,0&\ \,0 &\ \, 0 &\ \,0&\ \,0&\ \,0\\ \hline
\end{array}
\qquad 
\a_4=\frac{1}{\sqrt{8}}
\begin{array}{|rr|rr|rr|}
\hline
    2& \ \,0 &   2 &\ \,0&\ \, 0&\ \, 0  \\
   -2& \ \ 0 &  -2 &\ \,0&\ \, 0&\ \, 0\\
\ \,2&\ \, 0&\ \,2 &\ \,0&\ \, 0&\ \, 0\\
\ \,2&\ \,0 &\ \,2 &\ \,0&\ \, 0&\ \, 0\\ \hline
\end{array}
\]
\[
\a_5=\frac{1}{\sqrt{8}}
\begin{array}{|rr|rr|rr|}
\hline
     0&\ \,0 &     0 &   \ 0&   \ 0&   \ 0  \\
 \ \,0&\ \,0 &\ \, 0 &\ \, 0&\ \, 0&\ \, 0\\
    -4&\ \ 0 &    -4 &\ \ 0&\ \ 0&\ \, 0\\
 \ \,0&\ \,0 &\ \, 0 &\ \, 0&\ \, 0&\ \, 0\\ \hline
\end{array}
\qquad
\a_6=\frac{1}{\sqrt{8}}
\begin{array}{|rr|rr|rr|}
\hline
    0& \ \,2 &   0 &\ \,2&\ \, 0&\ \, 0  \\
    0& \ \,2 &  0 &\ \, 2&\ \, 0&\ \, 0\\
\ \,2&\ \, 0&\ \,2 &\ \,0&\ \, 0&\ \, 0\\
  - 2&\ \,0 &   -2 &\ \,0&\ \, 0&\ \, 0\\ \hline
\end{array}
\]
\[
\a_7=\frac{1}{\sqrt{8}}
\begin{array}{|rr|rr|rr|}
\hline
     0&   -4 &\ \  0 &  -4&\ \ 0   &\  \ 0  \\
 \ \,0&\ \,0 &\ \, 0 &\ \, 0&\ \, 0&\ \, 0\\
     0&\ \,0 &     0 &\ \, 0&\ \, 0&\ \, 0\\
 \ \,0&\ \,0 &\ \, 0 &\ \, 0&\ \, 0&\ \, 0\\ \hline
\end{array}
\qquad
\a_8=\frac{1}{\sqrt{8}}
\begin{array}{|rr|rr|rr|}
\hline
    0& \ \,2 &   0 &\ \,2&\ \, 0&\ \, 0  \\
    0&    -2 &   0 &   -2&\ \, 0&\ \, 0\\
\ \,0&\ \, 2 &\ \,0&\ \,2&\ \, 0&\ \, 0\\
    0&\ \, 2 &   0 &\ \,2&\ \, 0&\ \, 0\\ \hline
\end{array}
\]

In this case,  $M\cap N=0$ and $ann_N(M)= ann_M(N)=0$. Moreover, the set
$\{\b_1, \dots, \b_8, \a_1,\dots \a_8\}$ forms a basis of $L=M+N$ and
the Gram matrix of $L$ is

\[
\left[
\begin{array}{cccccccccccccccc}
4&0&-2&0&0&0&0&-2&1&-2&2&0&0&0&0&0\\
0&4&-2&0&0&0&0&0&1&-2&-2&2&0&0&0&0\\
-2&-2&4&-2&0&0&0&0&-1&2&0&-2&2&-1&0&0\\
0&0&-2&4&-2&0&0&0&1&0&0&0&-2&2&0&0\\
0&0&0&-2&4&-2&0&0&-1&0&0&1&0&-2&2&-1\\
0&0&0&0&-2&4&-2&0&1&0&0&0&0&0&-2&2\\
0&0&0&0&0&-2&4&0&-1&0&0&0&0&1&0&-2\\
-2&0&0&0&0&0&0&4&-2&1&-1&1&-1&1&-1&1\\
1&1&-1&1&-1&1&-1&-2&4&-2&0&0&0&0&0&0\\
-2&-2&2&0&0&0&0&1&-2&4&0&-2&0&0&0&0\\
2&-2&0&0&0&0&0&-1&0&0&4&-2&0&0&0&0\\
0&2&-2&0&1&0&0&1&0&-2&-2&4&-2&0&0&0\\
0&0&2&-2&0&0&0&-1&0&0&0&-2&4&-2&0&0\\
0&0&-1&2&-2&0&1&1&0&0&0&0&-2&4&-2&0\\
0&0&0&0&2&-2&0&-1&0&0&0&0&0&-2&4&-2\\
0&0&0&0&-1&2&-2&1&0&0&0&0&0&0&-2&4
\end{array}
\right].
\]
The Smith invariant sequence is $ 1111 1111 2222 2222$.

It is clear that $L\leq ann_\Lambda(E(\mathcal{O}_3))$ (see \refpp{octad} for the definition of $\mathcal{O}_3$).
On the other hand, $\det(L)=2^{8}=\det(ann_\Lambda(E(\mathcal{O}_3)))$. Hence, $L=ann_\Lambda(E(\mathcal{O}_3))$ is isomorphic to $BW_{16}$ (cf. Section \refpp{dih8FJ}).

\medskip

\noindent $\bf \dih{8}{16, DD_4}$

Define $M:=E(\mathcal{O}_2)\xi$ and  $N:=E(\mathcal{O}_4)$,  where
$\mathcal{O}_2$ and $\mathcal{O}_4$ are defined as in Notation
\refpp{octad}. We shall use the set $\{\g_1\xi, \dots, \g_8\xi\}$
defined in Notation \refpp{baisiofEO2} as a basis of $M$ and the set
$\{\a_1, \dots, \a_8\}$ as a basis of $N$, where

\[
\begin{split}
&\a_1=\frac{1}{\sqrt{8}}

\right].
\]
Therefore, no vector in $M=E(\mathcal{O}_2)\xi$ can be supported on
$\mathcal{O}_4$ and hence $M\cap N=0$. Moreover,  we have
\[
\begin{split}
ann_M(N)&=\{ \g\xi\in M|\, (\gamma\xi, \a_i)=0, \text{ for all } i=1, \dots, 8\}\\
&=\{ \g\xi\in M|\, supp\,\gamma\cap \mathcal{O}_4=\emptyset \}\\
& =\mspan_\ZZ\{\g_1\xi, \g_2\xi, \g_3\xi, \g_7\xi \} \cong DD_4.
\end{split}
\]

Set
\[
\begin{split}
\delta_1= \begin{array}{|rr|rr|rr|}
\hline 0 & 0& -1 & 0 &0 &0 \\
 0 & 0&  1 & 0 &0 &0 \\
 0 & 0& 1 & 0 &0 &0 \\
 0 & 0& 1 & 0 &0 &0 \\ \hline
\end{array}, \quad
\delta_2= \begin{array}{|rr|rr|rr|}
\hline 0& 0 & 0& -1 & 0 &0  \\
 0& 0 & 0&  1 & 0 &0  \\
 0&  0 & 0& 1 & 0 &0  \\
 0&  0 & 0& 1 & 0 &0  \\ \hline
\end{array}, \\
\delta_3= \begin{array}{|rr|rr|rr|}
\hline 0& 0 & 0&0& -1  &0 \\
 0& 0 & 0&  0& 1 &0  \\
 0& 0 & 0& 0& 1 &0  \\
 0& 0 & 0& 0& 1 &0 \\ \hline
\end{array}, \quad
\delta_4= \begin{array}{|rr|rr|rr|}
\hline 0& 0& 0& 0 & 0& -1  \\
 0& 0 &0& 0 & 0&  1   \\
 0& 0 &0&  0 & 0& 1   \\
 0& 0& 0&  0 & 0& 1  \\ \hline
\end{array}. \\
\end{split}
\]

Then
\[
\begin{split}
ann_N(M)= &\{ \a\in E(\mathcal{O}_4)| (\a, \g_i\xi )=0 \text{ for all } i=1, \dots,8\}\\
= &\{ \a\in E(\mathcal{O}_4)| (\a, \delta_i)=0 \text{ for all } i=1,2,3,4\}\\
= & \frac{1}{\sqrt{8}} \mspan_\ZZ
\left\{  {
\begin{array}{|rr|rr|rr|}
\hline   0& 0& 4& 0& 0& 0 \\
  0& 0& 4& 0& 0& 0\\
  0& 0& 0& 0& 0& 0\\
  0& 0& 0& 0& 0& 0\\
  \hline
  \end{array},\
\begin{array}{|rr|rr|rr|}
\hline   0& 0& 0& 4& 0& 0 \\
  0& 0& 0& 4& 0& 0\\
  0& 0& 0& 0& 0& 0\\
  0& 0& 0& 0& 0& 0\\
  \hline
  \end{array}, }
  \atop
{
\begin{array}{|rr|rr|rr|}
\hline
  0& 0& 0& 0& 4& 0 \\
  0& 0& 0& 0& 4& 0\\
  0& 0& 0& 0& 0& 0\\
  0& 0& 0& 0& 0& 0 \\
  \hline
  \end{array}, \
\begin{array}{|rr|rr|rr|}
\hline   0& 0&-2&-2&-2&-2\\
   0& 0&-2&-2&-2&-2\\
   0& 0& 0& 0& 0& 0\\
   0& 0& 0& 0& 0& 0\\
   \hline
\end{array}
  }
\right\}\\
\cong &DD_4
\end{split}
\]
In this case,  $L=M+N$ is of rank $16$  and $\{\a_1, \dots,\a_8, \g_1\xi,\dots, \g_8\xi\} $ is a basis of $L$. The Gram matrix of $L$ is given by
\[
\left[
\begin{array}{cc cc cc cc |cccccccc}
4&-2&0&0&0&0&0&0&0&0&0&2&-2&0&0&-1\\
-2&4&-2&0&0&0&2&0&0&0&0&-1&0&1&0&0\\
0&-2&4&-2&0&0&0&0&0&0&0&0&2&-2&0&1\\
0&0&-2&4&-2&0&0&0&0&0&0&0&-1&0&1&-1\\
0&0&0&-2&4&-2&0&0&0&0&0&0&0&2&-2&1\\
0&0&0&0&-2&4&0&0&0&0&0&0&0&-1&0&-1\\
0&2&0&0&0&0&4&-2&0&0&0&0&0&0&0&0\\
0&0&0&0&0&0&-2&4&0&0&0&0&0&0&0&0\\ \hline
0&0&0&0&0&0&0&0&4&-2&0&0&0&0&0&0\\
0&0&0&0&0&0&0&0&-2&4&-2&0&0&0&0&0\\
0&0&0&0&0&0&0&0&0&-2&4&-2&0&0&0&0\\
2&-1&0&0&0&0&0&0&0&0&-2&4&-2&0&0&0\\
-2&0&2&-1&0&0&0&0&0&0&0&-2&4&-2&0&2\\
0&1&-2&0&2&-1&0&0&0&0&0&0&-2&4&-2&0\\
0&0&0&1&-2&0&0&0&0&0&0&0&0&-2&4&0\\
-1&0&1&-1&1&-1&0&0&0&0&0&0&2&0&0&4
\end{array} \right]
\]
whose Smith invariant sequence is $111111112222
4444$.

A check of the Gram matrices also shows that
\[
ann_M(ann_M(N))=\mspan_\ZZ\{ \g_5\xi,\ \g_6\xi,\ \g_7'\xi,\ \g_8'\xi\}\cong DD_4
\]
and
\[
ann_N(ann_N(M))=\mspan\ZZ\{ \a_1,\ \a_3,\ \a_5,\
\a_8'\}\cong DD_4,
\]
where
\[
\g_7'=\frac{1}{\sqrt{8}}\,
\begin{array}{|rr|rr|rr|}
\hline  0& 0& 0& 0& 4&-4\\
  0& 0& 0& 0& 0& 0 \\
  0& 0& 0& 0& 0& 0\\
  0& 0& 0& 0& 0& 0\\
  \hline
  \end{array},\quad
\g_8'=\frac{1}{\sqrt{8}}\,
\begin{array}{|rr|rr|rr|}
\hline  0& 0& 0& 0& 4& 4\\
  0& 0& 0& 0& 0& 0 \\
  0& 0& 0& 0& 0& 0\\
  0& 0& 0& 0& 0& 0\\
  \hline
  \end{array},
\]
\[
  \a_8'=\frac{1}{\sqrt{8}}\,
  \begin{array}{|rr|rr|rr|}
\hline   0& 0&-2&-2&-2&-2\\
   0& 0& 2& 2& 2& 2\\
   0& 0& 0& 0& 0& 0\\
   0& 0& 0& 0& 0& 0\\
   \hline
\end{array}.
\]

Let $K=ann_M(ann_M(N))+ann_N(ann_N(M))$. Then $K$ is generated by
\[
\g_5\xi,\ \g_6\xi,\ \g_7'\xi,\ \g_8'\xi,\ \a_1,\ \a_3,\ \a_5,\
\a_8',
\]
The determinant of $K$ is $2^8$ and thus it is also isometric to
$EE_8$.

\bigskip

\noindent $\bf \dih{8}{15}$

When $M\cap N\cong AA_1$, this is the only possible case.

Let $\sigma_1$ and $\sigma_2$ be the involutions given as follows.
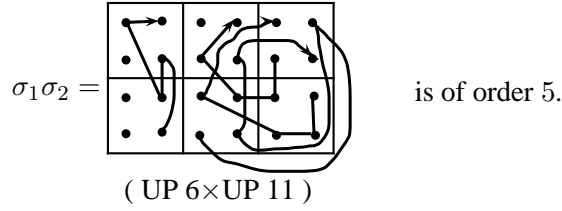
\[
\scalebox{1.0} 
{
%
}
\]
\vspace{0.5cm}

Let $M$ and $N$ be the $EE_8$ lattices corresponding to
$\sigma_1$ and $\sigma_2$, respectively.  Then,
\[
M\cap N =\mspan_\ZZ \left \{ \frac{1}{\sqrt{8}}
\begin{array}{|rr|rr|rr|}
\hline 0& 0& 0& 0& 0& 0\\
0& 0& 0& 0& 0& 0\\
0& 0&-2& 2&-2& 2\\
0& 0& 2&-2&-2& 2 \\ \hline
\end{array}
\right \}
\cong AA_1.
\]

Let
\[
\a_1=\frac{1}{\sqrt{8}}\,
\begin{array}{|rr|rr|rr|}
\hline 0& 0& 0& 0& 0& 0\\
0& 0& 0& 0& 0& 0\\
0& 0& 0& 0& 0&-4\\
0& 0& 0& 0& 4& 0 \\ \hline
\end{array}, \qquad
\a_2=\frac{1}{\sqrt{8}}
\begin{array}{|rr|rr|rr|}
\hline 0& 0& 0& 0& 0& 0\\
0& 0& 0& 0& 0& 0\\
0& 0&-2& 2&-2& 2\\
0& 0& 2&-2&-2& 2\\
\hline
\end{array},
\]
\[
\a_3=\frac{1}{\sqrt{8}} \,
\begin{array}{|rr|rr|rr|}  \hline
0& 0& 0& 0& 0& 0\\
0& 0& 0& 0& 0& 0\\
0& 0& 4& 0& 0& 0\\
0& 0&-4& 0& 0& 0
\\ \hline \end{array}, \qquad
\a_4=\frac{1}{\sqrt{8}} \, \begin{array}{|rr|rr|rr|}  \hline
0& 0& 0& 0&-2& 2\\
0& 0& 0& 0& 2&-2\\
0& 0&-2&-2& 0& 0\\
0& 0& 2& 2& 0& 0
\\ \hline \end{array},
\]
\[
\a_5=\frac{1}{\sqrt{8}} \, \begin{array}{|rr|rr|rr|}  \hline
0& 0& 0& 0& 4&-4\\
0& 0& 0& 0& 0& 0\\
0& 0& 0& 0& 0& 0\\
0& 0& 0& 0& 0& 0
\\ \hline \end{array}, \qquad
\a_6=\frac{1}{\sqrt{8}}
\,
\begin{array}{|rr|rr|rr|}  \hline
0& 0& 0& 0&-2& 2\\
0& 0& 0& 0&-2& 2\\
-2& 2& 0& 0& 0& 0\\
2&-2& 0& 0& 0& 0
\\ \hline \end{array},
\]
\[
\a_7=\frac{1}{\sqrt{8}} \,
\begin{array}{|rr|rr|rr|}  \hline
0& 0& 0& 0& 0& 0\\
  0& 0& 0& 0& 0& 0\\
  4& 0& 0& 0& 0& 0\\
  -4& 0& 0& 0& 0& 0
  \\ \hline \end{array}, \qquad
\a_8=\frac{1}{\sqrt{8}} \,
\begin{array}{|rr|rr|rr|}  \hline
  0& 0& 0& 0& 0& 0\\
  0& 0& 0& 0&-2& 2\\
  0&-2&-2& 0& 2& 0\\
  0& 2& 2& 0& 0&-2
  \\ \hline \end{array},
\]
\[
\a_1'=  \frac{1}{\sqrt{8}} \,
  \begin{array}{|rr|rr|rr|}  \hline
  0& 0& 0& 0& 0& 0\\
  0& 0& 0& 0& 0& 0\\
  0& 0& 0& 0& 0& 0\\
  0& 0&-4& 4& 0& 0
  \\ \hline \end{array},
\]
\[
\a_3'= \frac{1}{\sqrt{8}} \,
  \begin{array}{|rr|rr|rr|}  \hline
  0& 0& 0& 0& 0& 0\\
  0& 0& 0& 0& 0& 0\\
  0& 0& 4& 0& 0&-4\\
  0& 0& 0& 0& 0& 0
  \\ \hline \end{array}, \quad
\a_4'= \frac{1}{\sqrt{8}} \,
  \begin{array}{|rr|rr|rr|}  \hline
  0& 0& 0& 0& 0& 0\\
  0& 0&-2& 2& 2&-2\\
  0& 0&-2&-2& 2& 2\\
  0& 0& 0& 0& 0& 0
  \\ \hline \end{array},
\]
\[
\a_5'=\frac{1}{\sqrt{8}} \,
  \begin{array}{|rr|rr|rr|}  \hline
  0& 0& 0& 0& 0& 0\\
  0& 0& 4& 0&-4& 0\\
  0& 0& 0& 0& 0& 0\\
  0& 0& 0& 0& 0& 0
  \\ \hline \end{array},
\quad
\a_6'=\frac{1}{\sqrt{8}} \,
  \begin{array}{|rr|rr|rr|}  \hline
  2& 0& 0& 0& 0& 0\\
  0&-2&-2&-2& 2& 2\\
  0& 2& 0& 0& 0& 0\\
  0&-2& 0& 0& 0& 0
  \\ \hline \end{array},
\]
\[
\a_7'=
  \frac{1}{\sqrt{8}} \,
  \begin{array}{|rr|rr|rr|}  \hline
  0& 0& 0& 0& 0& 0\\
  0& 4& 0& 0& 0& 0\\
  0&-4& 0& 0& 0& 0\\
   0& 0& 0& 0& 0& 0
   \\ \hline \end{array},
\quad
\a_8'=
   \frac{1}{\sqrt{8}} \,
   \begin{array}{|rr|rr|rr|}  \hline
   -2& 0& 0& 0& 0& 0\\
   0& 0& 0&-2& 0& 2\\
   0& 0&-2& 0& 0& 2\\
   0& 2& 0& 0& 2&-2
   \\ \hline \end{array}.
\]

Note that $M\cap N=\mspan_\ZZ\{\a_2\}$. A check of Gram matrices also shows that
$\{\a_1,\a_2,\a_3, \dots,\a_8\}$ is a basis of $M$ and $\{\a_1',\a_2,\a_3', \dots,\a_8'\}$ is
a basis of $N$. Thus, $\{\a_1,\a_2,\a_3, \dots,\a_8\}\cup \{\a_1',\a_3', \dots,\a_8'\}$ is a basis of $L=M+N$ and the Gram matrix of $L$ is given by

\[
\left[
\begin{array}{ccccccccccccccc}
4&-2&0&0&0&0&0&0&0&2&-1&0&0&0&0\\
-2&4&-2&0&0&0&0&0&-2&-2&0&0&0&0&0\\
0&-2&4&-2&0&0&0&-2&2&2&-1&0&0&0&-1\\
0&0&-2&4&-2&0&0&0&0&-1&2&-1&0&0&0\\
0&0&0&-2&4&-2&0&0&0&0&0&0&0&0&0\\
0&0&0&0&-2&4&-2&0&0&0&-1&1&1&-1&0\\
0&0&0&0&0&-2&4&0&0&0&0&0&0&0&0\\
0&0&-2&0&0&0&0&4&-1&-1&0&1&-1&1&2\\
0&-2&2&0&0&0&0&-1&4&0&0&0&0&0&0\\
2&-2&2&-1&0&0&0&-1&0&4&-2&0&0&0&-2\\
-1&0&-1&2&0&-1&0&0&0&-2&4&-2&0&0&0\\
0&0&0&-1&0&1&0&1&0&0&-2&4&-2&0&0\\
0&0&0&0&0&1&0&-1&0&0&0&-2&4&-2&0\\
0&0&0&0&0&-1&0&1&0&0&0&0&-2&4&0\\
0&0&-1&0&0&0&0&2&0&-2&0&0&0&0&4
\end{array}
\right] \]

The Smith invariant sequence for $L$  is $ 1111 1 111 11 44 444$.

\medskip

\subsubsection{$ |g|=5$.} In this case, $M\cap N=0$ and $ann_M(N)=
ann_N(M)=0$.

\medskip
\noindent $\bf \dih{10}{16}$
Let $\sigma_1$ and $\sigma_2$ be the involutions given as follows:
\[
\scalebox{1.0} 
{
%
}
\]

\medskip

Let $M$ and $N$ be the $EE_8$ lattices corresponding to
$\sigma_1$ and $\sigma_2$, respectively. Then, $ M\cap N =0$ and
$ann_M(N)=ann_N(M)=0$. The lattice $L=M+N$ is of rank $16$. By Gram matrices, it is easy to check
\[
\a_1=\frac{1}{\sqrt{8}}\,
\begin{array}{|rr|rr|rr|}  \hline
 4&-4& 0& 0& 0& 0\\
0& 0& 0& 0& 0& 0\\
0& 0& 0& 0& 0& 0\\
0& 0& 0& 0& 0& 0\\ \hline \end{array}, \quad
\a_2=\frac{1}{\sqrt{8}}\,
\begin{array}{|rr|rr|rr|}  \hline
-2& 2& 0&-2& 0& 0\\
0& 0& 2& 0& 0& 0\\
0& 2& 0& 0& 0&-2\\
0&-2& 0& 0& 2& 0\\ \hline \end{array},
\]
\[
\a_3=\frac{1}{\sqrt{8}}\,
\begin{array}{|rr|rr|rr|}  \hline
 0& 0& 0& 4& 0& 0\\
0& 0&-4& 0& 0& 0\\
0& 0& 0& 0& 0& 0\\
0& 0& 0& 0& 0& 0\\ \hline \end{array}, \quad
\a_4=\frac{1}{\sqrt{8}}\,
\begin{array}{|rr|rr|rr|}  \hline
0& 0& 0&-2&-2& 0\\
0& 0& 2& 0& 0&-2\\
0& 0& 2& 0& 0& 2\\
0& 0& 0& 2&-2& 0\\ \hline \end{array},
\]
\[
\a_5=\frac{1}{\sqrt{8}}\,
\begin{array}{|rr|rr|rr|}  \hline
 0& 0& 0& 0& 4& 0\\
0& 0& 0& 0& 0& 0\\
0& 0&-4& 0& 0& 0\\
0& 0& 0& 0& 0& 0\\ \hline \end{array}, \quad
\a_6=\frac{1}{\sqrt{8}}\,
\begin{array}{|rr|rr|rr|}  \hline
0& 0& 0& 0&-2&-2\\
0& 0& 0& 0&-2& 2\\
0& 0& 2& 2& 0& 0\\
0& 0& 2&-2& 0& 0\\ \hline \end{array},
\]
\[
\a_7=\frac{1}{\sqrt{8}}\,
\begin{array}{|rr|rr|rr|}  \hline
0& 0& 0& 0& 0& 0\\
0& 0& 0& 0& 4& 0\\
0& 0& 0&-4& 0& 0\\
0& 0& 0& 0& 0& 0\\ \hline \end{array}, \quad
\a_8=\frac{1}{\sqrt{8}}\,
\begin{array}{|rr|rr|rr|}  \hline
 0& 0& 0& 0&-2& 2\\
0& 0& 0& 0& 0& 0\\
0&-2& 2& 0& 0&-2\\
0& 2&-2& 0& 2& 0\\ \hline \end{array},
\]
form a basis of $M$ and
\[
\a_1'=\frac{1}{\sqrt{8}}\,
\begin{array}{|rr|rr|rr|}  \hline
4& 0& 0& 0& 0& 0\\
0& 0& 0& 0& 0& 0\\
0& 0& 0& 0& 0& 0\\
0&-4& 0& 0& 0& 0\\ \hline \end{array}, \quad
\a_2'=\frac{1}{\sqrt{8}}\,
\begin{array}{|rr|rr|rr|}  \hline
-2& 0& 0& 0& 0& 0\\
0& 2&-2&-2& 2& 2\\
0&-2& 0& 0& 0& 0\\
0& 2& 0& 0& 0& 0\\ \hline \end{array},
\]
\[
\a_3'=\frac{1}{\sqrt{8}}\,
\begin{array}{|rr|rr|rr|}  \hline
0& 0& 0& 0& 0& 0\\
0& 0& 4& 0&-4& 0\\
0& 0& 0& 0& 0& 0\\
0& 0& 0& 0& 0& 0\\ \hline \end{array}, \quad
\a_4'=\frac{1}{\sqrt{8}}\,
\begin{array}{|rr|rr|rr|}  \hline
0& 0& 0& 0& 0& 0\\
0& 0&-2& 2& 2&-2\\
0& 0&-2&-2& 2& 2\\
0& 0& 0& 0& 0& 0\\ \hline \end{array},
\]
\[
\a_5'=\frac{1}{\sqrt{8}}\,
\begin{array}{|rr|rr|rr|}  \hline
0& 0& 0& 0& 0& 0\\
0& 0& 0& 0& 0& 0\\
0& 0& 4& 0& 0&-4\\
0& 0& 0& 0& 0& 0\\ \hline \end{array}, \quad
\a_6'=\frac{1}{\sqrt{8}}\,
\begin{array}{|rr|rr|rr|}  \hline
0& 0& 0& 0& 0& 0\\
0& 0& 0& 0& 0& 0\\
0& 0&-2& 2&-2& 2\\
0& 0&-2& 2&-2& 2\\ \hline \end{array},
\]
\[
\a_7'=\frac{1}{\sqrt{8}}\,
\begin{array}{|rr|rr|rr|}  \hline
0& 0& 0& 0& 0& 0\\
0& 0& 0& 0& 0& 0\\
0& 0& 0& 0& 0& 0\\
0& 0& 0& 0& 4&-4\\ \hline \end{array}, \quad
\a_8'=\frac{1}{\sqrt{8}}\,
\begin{array}{|rr|rr|rr|}  \hline
 0& 0& 0& 0& 0& 0\\
0&-2& 0&-2& 0& 2\\
0& 2&-2& 0& 0& 2\\
0& 0& 2&-2& 0& 0\\ \hline \end{array}
\]
form a basis of $N$. In addition, $\{\a_1, \dots,\a_8, \a_1', \dots, \a_8'\}$ is a basis of $L=M+N$.
The Gram matrix of $L$ is then given by

\[
\left[
\begin {array}{cccccccccccccccc}
4&-2&0&0&0&0&0&0&2&-1&0&0&0&0&0&0\\
 -2&4&-2&0&0&0&0&0&0&-1&1&-1&1&-1&1&0\\
 0&-2&4&-2&0&0&0&0&0&1&-2&1&0&0&0&0\\
 0&0&-2&4&-2&0&0&0&0&-1&1&0&0&1&-1&-1\\
 0&0&0&-2&4&-2&0&-2&0&0&0&1&-2&1&0&1\\
 0&0&0&0&-2&4&-2&0&0&0&1&-2&1&-1&0&1\\
 0&0&0&0&0&-2&4&0&0&1&-2&2&0&-1&0&0\\
 0&0&0&0&-2&0&0&4&-1&1&0&-1&2&-1&1&-2\\
 2&0&0&0&0&0&0&-1&4&-2&0&0&0&0&0&0\\
 -1&-1&1&-1&0&0&1&1&-2&4&-2&0&0&0&0&0\\
 0&1&-2&1&0&1&-2&0&0&-2&4&-2&0&0&0&0\\
 0&-1&1&0&1&-2&2&-1&0&0&-2&4&-2&0&0&0\\
 0&1&0&0&-2&1&0&2&0&0&0&-2&4&-2&0&-2\\
 0&-1&0&1&1&-1&-1&-1&0&0&0&0&-2&4&-2&0\\
 0&1&0&-1&0&0&0&1&0&0&0&0&0&-2&4&0\\
 0&0&0&-1&1&1&0&-2&0&0&0&0&-2&0&0&4
\end{array}
\right]
\]

The Smith invariant sequence is $1111 1111 1111 5555$.

\medskip
\subsubsection{$ |g|=6$.} In this case, $M\cap N=0$, and
$ann_N(M)\cong ann_M(N)\cong AA_2$.

\medskip
\noindent $\bf \dih{12}{16}$

Let $\sigma_1$ and $\sigma_2$ be the involutions given as follows:

\[
\scalebox{1.0} 
{
%
}
\]

\bigskip

Let $M$ and $N$ be the $EE_8$ lattices corresponding to
$\sigma_1$ and $\sigma_2$, respectively. Then, $ M\cap N =0$. Moreover, we have
\[
\begin{split}
ann_M(N)&=\frac{1}{\sqrt{8}} \mspan_\ZZ\left \{
\,
\begin{array}{|rr|rr|rr|}  \hline
0& 0& 0& 0& 0& 0\\
0& 0& 0& 0& 0& 0\\
4& 0& 0& 0& 0& 0\\
-4& 0& 0& 0& 0& 0
\\ \hline \end{array}
,
\
\begin{array}{|rr|rr|rr|}  \hline
0& 0& 0& 0& 2& -2\\
0& 0& 0& 0& 0& 0\\
-2& 0&2& 0& 2& 0\\
2& 0& -2& 0& 0&-2
\\ \hline \end{array}
 \right\}\\
& \cong AA_2
\end{split}
\]
and
\[
\begin{split}
ann_N(M)&=\frac{1}{\sqrt{8}} \mspan_\ZZ\left \{
\,
\begin{array}{|rr|rr|rr|}  \hline
0& 0& 0& 0& 0& 0\\
0& 0& 0& 0& 0& 0\\
0& 4& 0& 0& 0& 0\\
0& -4& 0& 0& 0& 0
\\ \hline \end{array}
,
\
\begin{array}{|rr|rr|rr|}  \hline
0&  0&  0& 0& 0& 0\\
0&  0&  0& 0& 2& -2\\
0& -2&  0& -2& 0& 2\\
0& -2&  0&  2& -2& 0
\\ \hline \end{array}
 \right\}\\
& \cong AA_2.
\end{split}
\]
In this case, $L= M+N$ is of rank $16$. By Gram matrices, it is easy to check
\[
\a_1=\frac{1}{\sqrt{8}}\,
\begin{array}{|rr|rr|rr|}  \hline
4&-4& 0& 0& 0& 0\\
0& 0& 0& 0& 0& 0\\
0& 0& 0& 0& 0& 0\\
0& 0& 0& 0& 0& 0
\\ \hline \end{array}
, \quad
\a_2=\frac{1}{\sqrt{8}}\,
\begin{array}{|rr|rr|rr|}  \hline
-2 & 2& 0&-2& 0& 0\\
 0 & 0& 2& 0& 0& 0\\
 2 & 0& 0& 0&-2& 0\\
-2 & 0& 0& 0& 0& 2
 \\ \hline \end{array}
 ,
\]
\[
\a_3=\frac{1}{\sqrt{8}}
\,
\begin{array}{|rr|rr|rr|}  \hline
 0& 0& 0& 4& 0& 0\\
 0& 0&-4& 0& 0& 0\\
 0& 0& 0& 0& 0& 0\\
 0& 0& 0& 0& 0& 0
 \\ \hline \end{array}
, \quad
\a_4=\frac{1}{\sqrt{8}}\,
\begin{array}{|rr|rr|rr|}  \hline
  0& 0& 0&-2& 2& 0\\
  0& 0& 2& 0& 0& 2\\
  0& 0& 0&-2& 2& 0\\
  0& 0&-2& 0& 0&-2
\\ \hline \end{array} ,
\]
\[
\a_5=\frac{1}{\sqrt{8}}\,
\begin{array}{|rr|rr|rr|}  \hline
 0& 0& 0& 0&-4& 0\\
 0& 0& 0& 0& 0& 0\\
 0& 0& 0& 0& 0& 0\\
 0& 0& 4& 0& 0& 0
\\ \hline \end{array} , \quad
\a_6=\frac{1}{\sqrt{8}}
\, \begin{array}{|rr|rr|rr|}  \hline
0& 0& 0& 0& 2& 2\\
0& 0& 0& 0& 2&-2\\
0& 0&-2& 2& 0& 0\\
0& 0&-2&-2& 0& 0
\\ \hline \end{array} ,
\]
\[
\a_7=\frac{1}{\sqrt{8}}\,
\begin{array}{|rr|rr|rr|}  \hline
0& 0& 0& 0& 0& 0\\
0& 0& 0& 0&-4& 0\\
0& 0& 0& 0& 0& 0\\
0& 0& 0& 4& 0& 0
\\ \hline \end{array} , \quad
\a_8=\frac{1}{\sqrt{8}}\,
\begin{array}{|rr|rr|rr|}  \hline
 0& 0& 0& 0& 2&-2\\
 0& 0& 0& 0& 0& 0\\
-2& 0& 2& 0&-2& 0\\
 2& 0&-2& 0& 0& 2
\\ \hline \end{array} , \]
form a basis of $M$ and
\[
\a_1'=\frac{1}{\sqrt{8}}\,
\begin{array}{|rr|rr|rr|}  \hline
 4& 0& 0& 0& 0& 0\\
 0&-4& 0& 0& 0& 0\\
 0& 0& 0& 0& 0& 0\\
 0& 0& 0& 0& 0& 0
\\ \hline \end{array} , \quad
\a_2'=\frac{1}{\sqrt{8}}\,
\begin{array}{|rr|rr|rr|}  \hline
-2& 0& 0& 0& 0& 0\\
 0& 2&-2& 2&-2& 2\\
 0& 2& 0& 0& 0& 0\\
 0&-2& 0& 0& 0& 0
\\ \hline \end{array} , \]
\[
\a_3'=\frac{1}{\sqrt{8}}\,
\begin{array}{|rr|rr|rr|}  \hline
0& 0& 0& 0& 0& 0\\
0& 0& 4&-4& 0& 0\\
0& 0& 0& 0& 0& 0\\
0& 0& 0& 0& 0& 0
\\ \hline \end{array} , \quad
\a_4'=\frac{1}{\sqrt{8}}\,
\begin{array}{|rr|rr|rr|}  \hline
 0& 0& 0& 0& 0& 0\\
 0& 0&-2& 2& 2&-2\\
 0& 0&-2&-2& 2& 2\\
 0& 0& 0& 0& 0& 0
\\ \hline \end{array}
, \]
\[
\a_5'=\frac{1}{\sqrt{8}}\,
\begin{array}{|rr|rr|rr|}  \hline
0& 0& 0& 0& 0& 0\\
0& 0& 0& 0& 0& 0\\
0& 0& 4& 0&-4& 0\\
0& 0& 0& 0& 0& 0
\\ \hline \end{array} , \quad
\a_6'=\frac{1}{\sqrt{8}}\,
\begin{array}{|rr|rr|rr|}  \hline
0& 0& 0& 0& 0& 0\\
0& 0& 0& 0& 0& 0\\
0& 0&-2& 2& 2&-2\\
0& 0&-2&-2& 2& 2
\\ \hline \end{array} , \]
\[
\a_7'=\frac{1}{\sqrt{8}}\,
\begin{array}{|rr|rr|rr|}  \hline
0& 0& 0& 0& 0& 0\\
0& 0& 0& 0& 0& 0\\
0& 0& 0& 0& 0& 0\\
0& 0& 0& 4&-4& 0
\\ \hline \end{array} , \quad
\a_8'=\frac{1}{\sqrt{8}}\,
\begin{array}{|rr|rr|rr|}  \hline
0& 0& 0& 0& 0& 0\\
0& 0& 0& 0&-2& 2\\
0&-2&-2& 0& 2& 0\\
0& 2& 2& 0& 0&-2
\\ \hline \end{array} .
\]
form a basis for $N$. Note that $\{\a_1,\dots,\a_8, \a_1',\cdots,\a_8'\}$ is a  basis of $L$ and the Gram matrix of $L$ is given by
\[
\left[
\begin {array}{cccccccc|cccccccc}
  4&-2&0&0&0&0&0&0&2&-1&0&0&0&0&0&0\\
 -2&4&-2&0&0&0&0&0&-1&0&1&-1&1&0&0&-1\\
  0&-2&4&-2&0&0&0&0&0&1&-2&1&0&0&0&0\\
  0&0&-2&4&-2&0&0&0&0&0&1&0&-1&0&0&1\\
  0&0&0&-2&4&-2&0&-2&0&0&0&0&0&-1&0&1\\
  0&0&0&0&-2&4&-2&0&0&-1&0&1&-1&2&-1&-1\\
  0&0&0&0&0&-2&4&0&0&1&0&-1&0&-1&2&1\\
  0&0&0&0&-2&0&0&4&0&0&0&-1&2&0&0&-2\\  \hline
  2&-1&0&0&0&0&0&0&4&-2&0&0&0&0&0&0\\
  -1&0&1&0&0&-1&1&0&-2&4&-2&0&0&0&0&0\\
  0&1&-2&1&0&0&0&0&0&-2&4&-2&0&0&0&0\\
  0&-1&1&0&0&1&-1&-1&0&0&-2&4&-2&0&0&0\\
  0&1&0&-1&0&-1&0&2&0&0&0&-2&4&-2&0&-2\\
  0&0&0&0&-1&2&-1&0&0&0&0&0&-2&4&-2&0\\
  0&0&0&0&0&-1&2&0&0&0&0&0&0&-2&4&0\\
 0&-1&0&1&1&-1&1&-2&0&0&0&0&-2&0&0&4
\end {array} \right]
\]
\medskip

The Smith invariant sequence is $1111 1111 1111 6666$.

\medskip

Now let $g=(\sigma_1\sigma_2)^2$. Then
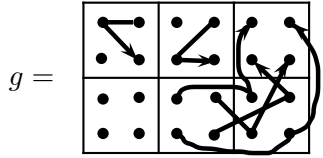
\[
\scalebox{1} 
{
%
}
\]
Note that $Mg$ is also isometric to $EE_8$ and it has a basis
\[
\a_1 g =\frac{1}{\sqrt{8}}\,
\begin{array}{|rr|rr|rr|}  \hline
-4& 0& 0& 0& 0& 0\\
0& 4& 0& 0& 0& 0\\
0& 0& 0& 0& 0& 0\\
0& 0& 0& 0& 0& 0
\\ \hline \end{array}
, \quad
\a_2 g=\frac{1}{\sqrt{8}}\,
\begin{array}{|rr|rr|rr|}  \hline
 2 & 0& 0& 0&-2& 2\\
 0 &-2&-2& 2& 0& 0\\
 2 & 0& 0& 0& 0& 0\\
-2 & 0& 0& 0& 0& 0
 \\ \hline \end{array}
 ,
\]
\[
\a_3 g=\frac{1}{\sqrt{8}}
\,
\begin{array}{|rr|rr|rr|}  \hline
 0& 0& 0& 0& 0& 0\\
 0& 0& 4&-4& 0& 0\\
 0& 0& 0& 0& 0& 0\\
 0& 0& 0& 0& 0& 0
 \\ \hline \end{array}
, \quad
\a_4 g=\frac{1}{\sqrt{8}}\,
\begin{array}{|rr|rr|rr|}  \hline
  0& 0& 0& 0& 2&-2\\
  0& 0&-2& 2& 0& 2\\
  0& 0& 2& 2& 2& 0\\
  0& 0& 0& 0&-2&-2
\\ \hline \end{array} ,
\]
\[
\a_5 g=\frac{1}{\sqrt{8}}\,
\begin{array}{|rr|rr|rr|}  \hline
 0& 0& 0& 0& 0& 0\\
 0& 0& 0& 0& 0& 0\\
 0& 0&-4& 0& 0& 0\\
 0& 0& 0& 0& 0& 4
\\ \hline \end{array} , \quad
\a_6 g=\frac{1}{\sqrt{8}}
\, \begin{array}{|rr|rr|rr|}  \hline
0& 0& 0& 0& 0& 0\\
0& 0& 0& 0& 0& 0\\
0& 0& 2&-2&-2&-2\\
0& 0& 2& 2& 2&-2
\\ \hline \end{array} ,
\]
\[
\a_7 g=\frac{1}{\sqrt{8}}\,
\begin{array}{|rr|rr|rr|}  \hline
0& 0& 0& 0& 0& 0\\
0& 0& 0& 0& 0& 0\\
0& 0& 0& 0& 0&-4\\
0& 0& 0& 4& 0& 0
\\ \hline \end{array} , \quad
\a_8 g=\frac{1}{\sqrt{8}}\,
\begin{array}{|rr|rr|rr|}  \hline
 0& 0& 0& 0&-2& 2\\
 0& 0& 0& 0& 0& 0\\
-2& 0& 2& 0& 2& 0\\
 2& 0&-2& 0& 0&-2
\\ \hline \end{array} . \]

Hence we have
\[
\begin{split}
 Mg\cap N &= \frac{1}{\sqrt{8}} \mspan_\ZZ\left\{
{%
}
\]
and it commutes with $t_N$. In this case, $t_{Mg}$ and $t_N$ generates a dihedral group of order $4$ and $Mg +N$ is isometric to the lattice $\dih{4}{12}$.

{}

\end{document}